\documentclass[10pt]{amsart}
\usepackage[dvips]{graphicx}
\RequirePackage{amssymb}
\RequirePackage[T1]{fontenc}

\newtheorem{theorem}{Theorem}[section]
\newtheorem{lemma}[theorem]{Lemma}
\newtheorem{conjecture}[theorem]{Conjecture} 
\newtheorem{proposition}[theorem]{Proposition}
 
\theoremstyle{definition}

\newtheorem{remark}[theorem]{Remark}

\newcommand{\Tr}{\text{Tr}}

\newcommand{\nc}{\newcommand}
\newcommand{\on}{\operatorname}

\nc{\cD}{{\mathcal D}}

\newcommand{\kk}{{\bf k}}

\newcommand{\QQ}{{\mathbb Q}}
\newcommand{\CC}{{\mathbb C}}
\newcommand{\RR}{{\mathbb{R}}}
\newcommand{\NN}{{\mathbb N}}

\renewcommand{\Tr}{{\mathfrak T}}
\nc{\SL}{{\mathfrak sl}}
\nc{\gt}{{\mathfrak gt}}
\nc{\grt}{\mathfrak{grt}}
\nc{\krv}{\mathfrak{krv}}
\nc{\kv}{\mathfrak{kv}}
\renewcommand{\t}{{\mathfrak{t}}}
\nc{\vt}{{\ov{\t}}}

\renewcommand{\u}{{\mathfrak{u}}}
\newcommand{\pb}{{\mathfrak{pb}}}

\nc{\G}{{\mathfrak{g}}}
\renewcommand{\u}{{\mathfrak{u}}}

\newcommand{\f}{{\mathfrak{f}}}
\newcommand{\tder}{{\mathfrak{tder}}}

\nc{\wh}{\widehat}\nc{\wt}{\widetilde} 
\renewcommand{\i}{\on{i}}

\newcommand{\ov}{\overline}

\newcommand{\nn}{{\bf n}}

\newcommand{\Mg}{{\on{Mg}}}
\newcommand{\Sg}{{\on{Sg}}}
\newcommand{\ben}{\begin{enumerate}}
\newcommand{\een}{\end{enumerate}}

\newcommand{\cC}{{\mathcal C}}

\newcommand{\lbr}{[\![}
\newcommand{\rbr}{]\!]}

\begin{document}

\title[Drinfeld associators and solutions of the Kashiwara--Vergne equations]
{Drinfeld associators, braid groups and explicit solutions of the 
Kashiwara--Vergne equations}

\author{A. Alekseev}
\address{Section de math\'ematiques, Universit\'e de Gen\`eve, 2-4 rue du 
Li\`evre, c.p. 64, 1211 Gen\`eve 4, Switzerland}
\email{alekseev@math.unige.ch}

\author{B. Enriquez}
\address{IRMA (CNRS UMR7501) et Universit\'e de Strasbourg, 
7 rue Ren\'e Descartes, 67084 Strasbourg cedex, France }
\email{enriquez@math.u-strasbg.fr}

\author{C. Torossian}
\address{Institut Mathématiques de Jussieu, Universit\'e Paris 7, CNRS;
Case 7012, 2 place Jussieu, 75005 Paris, France}
\email{torossian@math.jussieu.fr}

\maketitle

\begin{abstract} 
The Kashiwara--Vergne (KV) conjecture states the existence of solutions 
of a pair of equations related with the Campbell--Baker--Hausdorff
series. It was solved by Meinrenken and the first author over ${\mathbb R}$, 
and in a formal version, by the first and last authors over a field of 
characteristic 0. In this paper, we give a simple and explicit
formula for a map from the set of Drinfeld associators 
to the set of solutions of the formal KV equations. Both 
sets are torsors under the actions of prounipotent groups, and we show 
that this map is a morphism of torsors. When specialized to the KZ associator, 
our construction yields a solution over $\RR$ of the original KV conjecture. 
\end{abstract}

\section*{Introduction}

In \cite{KV}, M. Kashiwara and M. Vergne formulated a conjecture on the form of the 
Campbell--Baker--Hausdorff (CBH) series. This conjecture triggered the 
work of several authors (for a review see \cite{T2}). In particular, 
Kashiwara--Vergne settled it for solvable Lie algebras (\cite{KV}), 
Rouvi\`ere gave a proof for ${\mathfrak{sl}}_2$ (\cite{R}), and 
Vergne (\cite{V}) and Alekseev--Meinrenken 
(\cite{AM1}) proved it for quadratic Lie algebras. All these constructions
lead to explicit rational formulas for solutions of the KV conjecture. 
The general case was settled in the positive by Alekseev--Meinrenken 
(\cite{AM2}) using Kontsevich's deformation quantization theory
and results in \cite{T1}. The corresponding solution is defined over 
${\mathbb R}$, and expresses as an infinite series where coefficients are
combinations of Kontsevich integrals on configuration spaces
and integrals over simplices. The values of most of these coefficients
remain unknown.

Later, the first and last authors
gave another proof (\cite{AT}), based on Drinfeld's theory of associators. 
In that paper, the Kashiwara--Vergne (KV) conjecture was reformulated as the 
problem of constructing special automorphisms of the free Lie 
algebra with two generators with coboundary Jacobian (see Section 
\ref{sec:result}); 
the authors also showed that each associator gives rise to an affine line of 
such automorphisms. The solution is defined as a nonabelian cochain
with coboundary equal to the associator. Such a construction is 
inspired by the theory of quantization of Lie bialgebras, and the existence
problem is solved by showing that obstructions vanish in all degrees. 

The purpose of the present work is to give a direct construction of the map 
$M_1(\kk)\to \on{SolKV}(\kk)$, $\Phi\mapsto \mu_{\Phi}$ from associators 
to solutions of the 
KV equations (we work over a field $\kk$ of characteristic $0$). 
Namely, for $\Phi\in M_1(\kk)$, $\mu_\Phi$ is the automorphism 
of the topologically free Lie algebra generated by $x,y$ given by 
\begin{equation} \label{0}
\mu_\Phi : x\mapsto \Phi(x,-x-y)x\Phi(x,-x-y)^{-1}, \quad 
y \mapsto e^{-(x+y)/2}\Phi(y,-x-y)y\Phi(y,-x-y)^{-1}e^{(x+y)/2}. 
\end{equation}
Our main result (Theorem \ref{thm:main}) is the identity 
\begin{equation}\label{Phimumu}
\Phi(t_{12},t_{23})\circ\mu_{\Phi}^{12,3}\circ\mu_{\Phi}^{1,2} = 
\mu_{\Phi}^{1,23}\circ\mu_{\Phi}^{2,3}. 
\end{equation} 
This identity implies that the Jacobian of $\mu_\Phi$ is a 
cocycle, and therefore a coboundary according to cohomology 
computations in \cite{AT}; it can then be expressed using the $\Gamma$-function 
$\Gamma_{\Phi}$ of $\Phi$ (see \cite{DT,E}). Identity (\ref{Phimumu}) also implies 
that $\mu_\Phi$ is special, i.e., satisfies 
\begin{equation} \label{damien}
\mu_{\Phi}(\on{log}(e^{x}e^{y}))=x+y  
\end{equation}
(see Subsection \ref{X:Y} and also Proposition 7.4 
in \cite{AT}); we also give a direct proof of (\ref{damien}) 
based on the hexagon and duality identities satisfied by $\Phi$.  
The conjunction of (\ref{damien}) and of the fact that the 
Jacobian of $\mu_\Phi$ is a coboundary actually means that 
$\mu_\Phi$ is a solution of the KV equations introduced in 
\cite{AT}. 

The affine line of solutions of the KV equations 
attached in \cite{AT} to $\Phi$ then takes the form 
$\{\on{Inn}(e^{s(x+y)}) \circ \mu_\Phi, s\in \kk\}$, 
where $\on{Inn}(g) = (u\mapsto gug^{-1})$. It remains
an open question whether all  the solutions of the KV
equation are of this form. 

The strategy for proving (\ref{Phimumu}) is as follows. 
For each associator $\Phi$ and each parenthesization $O$ of a 
word in $n$ identical letters (the letter is $\bullet$), 
Drinfeld and Bar-Natan define  
an isomorphism $\tilde\mu_{\Phi}^{O}
:\on{PB}_n(\kk)\to \on{exp}(\hat\t_n)$ from the 
prounipotent completion of the pure braid group with $n$ 
strands to the group associated with the holonomy Lie algebra. 
Note that $\on{PB}_n$
contains the free group $\on{F}_{n-1}$ as a normal subgroup, 
while $\t_n$ contains the free Lie algebra $\f_{n-1}$
as an ideal; we show that the above isomorphisms restrict to isomorphisms 
$\mu_{\Phi}^{O}:\on{F}_{n-1}(\kk)\to \on{exp}(\hat\f_{n-1})$
(in the case of the left parenthesization, this was proved in \cite{HM}). 
We note that $\mu_\Phi$
may be interpreted as the isomorphism $\on{F}_{2}(\kk)\to 
\on{exp}(\hat\f_{2})$ corresponding to $\bullet(\bullet\bullet)$, so 
$\mu_\Phi = \mu_{\bullet(\bullet\bullet)}$ (we write $\mu_{O}$
instead of $\mu^{O}_{\Phi}$ when no confusion is possible).
We then show the identity 
\begin{equation} \label{muO:muOnew}
\mu_{O^{(i)}} = \mu_{O}^{1,2,...,ii+1,...,n}\circ 
\mu_{\bullet(\bullet\bullet)}^{i,i+1}, 
\end{equation}
where $O$ is a parenthesized word of length $n$ and $O^{(i)}$ is 
the parenthesized word obtained from it by replacing the $(i+1)$th
letter $\bullet$ by $(\bullet\bullet)$. 
Applying this identity to $O=\bullet(\bullet\bullet)$ with $i=1,2$
and using the identity $\mu_{\Phi}^{O'} = \on{Ad}(\Phi_{O,O'})\circ 
\mu_{\Phi}^{O}$
relating the various $\mu_\Phi^O$, we obtain (\ref{Phimumu}). 

We then study the torsor aspects of the map $\Phi\mapsto\mu_{\Phi}$. While 
$M_{1}(\kk)$ is a torsor under the commuting actions of the groups 
$\on{GT}_{1}(\kk)$ and $\on{GRT}_{1}(\kk)$, $\on{SolKV}(\kk)$ is a torsor
under the actions of groups $\on{KV}(\kk)$ and $\on{KRV}(\kk)$. 
We prove that $\Phi\mapsto\mu_{\Phi}$ is a morphism of torsors, i.e., there exist 
group morphisms $\on{GT}_{1}(\kk)\to \on{KV}(\kk)$, $f\mapsto \alpha_{f}$ 
and $\on{GRT}_{1}(\kk) \to \on{KRV}(\kk)$, $g\mapsto a_{g}$, compatible 
with the actions (the Lie algebra version of the latter morphism was already constructed in 
\cite{AT}). We give a direct proof of these facts, based on the nonemptiness of 
$M_{1}(\kk)$ (a result in \cite{Dr:Gal}); we also sketch an independent proof of 
$\alpha_{f}\in \on{KV}(\kk)$; its main ingredient is the identity 
\begin{equation} \label{eq:GT}
\on{Ad}f(x_{12},x_{23})\circ \alpha_{f}^{\widetilde{12},3}
\circ \alpha_{f}^{1,2} = \alpha_{f}^{1,\widetilde{23}}\circ \alpha_{f}^{2,3}.  
\end{equation}
A similar independent proof of $a_{g}\in \on{KRV}(\kk)$ may be 
given based on 
$$
\on{Ad}g(t_{12},t_{23})\circ a_{g}^{12,3}
\circ a_{g}^{1,2} = a_{g}^{1,{23}}\circ a_{g}^{2,3}.  
$$ 
It can be proved using the techniques of \cite{AT} that the 
sets of solutions of both equations are affine lines, and our result
gives explicit formulas for these solutions. We also observe that 
(\ref{eq:GT}) can be generalized to the profinite and pro-$l$ 
setups (i.e., we have morphisms $\widehat{\on{GT}}
\to \on{Aut}(\widehat{\on{F}}_{2})$
and $\on{GT}_{l}\to \on{Aut}((\on{F}_{2})_{l})$, $f\mapsto \alpha_{f}$, 
and (\ref{eq:GT}) takes place in $\on{Aut}(\widehat{\on{F}}_{3})$ or 
$\on{Aut}((\on{F}_{3})_{l})$). 
 
Formula (\ref{muO:muOnew}) and its analogue (\ref{eq:GT}) 
then enable us to compute the Jacobians of 
$\mu_{\Phi}^{O}:\on{F}_{n-1}(\kk)\to 
\on{exp}(\f_{n-1})$ and $\alpha_{f}^{O}\in \on{Aut}(\on{F}_{n-1}(\kk))$, where
$O$ is an arbitrary parenthesized word, $\Phi\in M_{1}(\kk)$, 
$f\in \on{GT}_{1}(\kk)$, in terms of 
in terms of $\Gamma_{\Phi}$ and of the `$\Gamma$-function' of $f$. 

Finally, we show that specializing our construction to the 
Knizhnik--Zamolodchikov (KZ) associator yields an explicit 
solution of the original KV conjecture, where the Lie 
series are required to converge for any finite dimensional 
Lie algebra and the Duflo series is required to coincide with the generating 
series of Bernoulli numbers. 

\medskip 
{\bf Acknowledgements.} We are grateful to V.G. Drinfeld and to 
D. Bar-Natan who posed the question of how to construct 
explicit solutions of the KV problem in terms of associators. 
The formula for $\mu_{\Phi}$ as well as the proof 
of equation (\ref{damien}) were suggested to us by D. Calaque. 
We would like to thank G. Massuyeau for discussions and for pointing
out reference \cite{HM}.
The research of A.A. was supported in part by the grants 200020-120042 and
200020-121675 of the Swiss National Science Foundation, and the research 
of C.T. was supported by CNRS.  

\tableofcontents

\section{Preliminary results} \label{sec:back}

In this section, we recall the notions of tangential derivations and 
automorphisms of free Lie algebras, their divergence and Jacobian cocycles, the actions
of pure braid groups (resp., infinitesimal braid Lie algebras) on free groups
Lie algebras by tangential automorphisms (resp., derivations), and 
simplicial morphisms between these objects. 

\subsection{Tangential automorphisms, the Jacobian cocycle, and complexes}
\label{sect:tang}

Let $\f_{n}$ be the free Lie algebra with generators $x_{1},...,x_{n}$, $\hat\f_{n}$
its degree completion (where the generators $x_{k}$ have degree $1$). For 
$u_{1},...,u_{n}\in\f_{n}$, we denote by $\lbr u_{1},...,u_{n}\rbr$
the derivation of $\f_{n}$ given by $x_{k}\mapsto [u_{k},x_{k}]$. 
In this way, we define a linear map $(\f_{n})^{n}\to \on{Der}(\f_{n})$. 
Its image is a (positively) graded Lie subalgebra $\tder_{n}$ of $\on{Der}(\f_{n})$; 
its elements are called the tangential derivations of $\f_{n}$. 
We similarly define $\tder_{n}^{\wedge} \subset \on{Der}(\hat\f_{n})$
as the degree completion of $\tder_{n}$; it is a pronilpotent Lie algebra.

If $U_{1},...,U_{n}\in \on{exp}(\hat\f_{n})$, we similarly define 
$\lbr U_{1},...,U_{n}\rbr$ as the automorphism of $\hat\f_{n}$ given by 
$x_{k}\mapsto U_{k}x_{k}U_{k}^{-1}$. This defines a map 
$\on{exp}(\hat\f_{n})^{n}\to \on{Aut}(\hat\f_{n})$, whose 
image is the subgroup of tangential automorphisms
$\on{TAut}_{n}\subset \on{Aut}(\hat\f_{n})$. The exponential sets up an 
isomorphism $\on{exp}:\tder_{n}^{\wedge}\to \on{TAut}_{n}$. 

Define $\Tr_{n}:= A_{n}/[A_{n},A_{n}]$  as the quotient of the  
free associative algebra $A_{n} \simeq U(\f_{n})$ by its subspace of 
commutators; this is the vector space spanned by the set of cyclic words 
in $x_1,...,x_n$. $\Tr_{n}$ is equipped with an action of $\on{Der}(\f_{n})$, 
induced by the action of $\on{Der}(\f_{n})$ on $A_{n}$. We denote by 
$x\mapsto \langle x\rangle$ the canonical projection map $A_{n}\to\Tr_{n}$. 
$\Tr_{n}$ is positively graded and we denote by $\hat\Tr_{n}$ its degree 
completion; it is equipped with actions of $\on{Der}(\hat\f_{n})$ and 
$\on{Aut}(\hat\f_{n})$. 

One shows that any $u\in \tder_{n}$ can be written as
$u=\lbr u_{1},...,u_{n}\rbr$, where $(u_{1},...,u_{n})$ 
is uniquely determined by 
the condition $p_{1}(u_{1})=...=p_{n}(u_{n})=0$, where 
$p_{k} : \f_{n}\to\kk$ is the linear map 
such that $u=\sum_{k}p_{k}(u)x_{k}$ modulo $[\f_{n},\f_{n}]$. 

We define simplicial group morphisms $\on{TAut}_n\to\on{TAut}_m$
as follows. Let\footnote{We set $[n]:=\{1,...,n\}$.} 
$\phi : [m]\supset D_\phi\to [n]$ be a partially
defined map, and let $(a_1,...,a_n)\in (\f_n)^n$ be such that each $a_k$
has vanishing linear term in $x_k$. We set 
$\lbr a_1,...,a_n\rbr^\phi := \lbr b_1,...,b_m\rbr$, where 
$b_\ell(x_1,...,x_m):= a_{\phi(\ell)}(\sum_{k\in \phi^{-1}(1)}x_k,...,
\sum_{k\in \phi^{-1}(n)}x_k)$. This formula defines a Lie algebra 
morphism $\tder_n\to \tder_m$, which induces a group morphism 
$\on{TAut}_n \to \on{TAut}_m$, also denoted $x\mapsto x^\phi$. 
We will also use the notation $x^\phi = x^{\phi^{-1}(1),...,\phi^{-1}(n)}$. 
For example, $\lbr a_1,a_2\rbr^{12,3}=\lbr a_1(x_1+x_2,x_3),a_1(x_1+x_2,x_3),
a_2(x_1+x_2,x_3)\rbr$. 

We also define noncommutative variants of these morphisms 
as follows. Let $\tilde\phi$ be a pair consisting of a partially defined 
map $\phi : [m]\supset D_\phi\to [n]$ as above and of total orders on 
each of the sets $\phi^{-1}(1),...,\phi^{-1}(n)$. We define  
$\lbr a_1,...,a_n\rbr^{\tilde\phi}:= \lbr \tilde b_1,...,\tilde b_m\rbr$, 
where $\tilde b_\ell(x_1,...,x_m):= a_{\phi(\ell)}(\on{cbh}(x_k|k\in
\phi^{-1}(1)),...,\on{cbh}(x_k|k\in
\phi^{-1}(n)))$; here $\on{cbh}(a_1,...,a_p) = \on{log}(e^{a_1}...e^{a_p})$
and $\on{cbh}(a_s|s\in S)$ is defined similarly, for $S$ a finite 
ordered set. We use the notation $x^{\tilde\phi} = x^{\widetilde{\phi^{-1}(1)},
...,\widetilde{\phi^{-1}(n)}}$ (where the elements of $\phi^{-1}(k)$
are written in increasing order). 

We then define a `divergence' map 
$$
j : \tder_{n}\to \Tr_{n}
$$
as follows. Let 
$\partial_{k}:A_{n}\to A_{n}$ be the linear maps defined 
by the identity $x = \epsilon(x)1+
\sum_{k=1}^{n}x_{k}\partial_{k}(x)$ (where $\epsilon:A_{n}\to\kk$
is the counit map). We then set 
$$
j(u):= \langle \sum_{k=1}^{n}x_{k}\partial_{k}(u_{k})\rangle.
$$
One can show that $j$ satisfies the cocycle identity 
$$
j([u,v]) = u \cdot j(v) - v\cdot j(u), 
$$ 
where the action of $\tder_{n}$ on $\Tr_{n}$ is understood in the r.h.s.; 
$j$ is graded, so it extends to a cocycle $\tder_{n}^{\wedge}\to\hat\Tr_{n}$. 
The Lie algebra cocycle $j$ gives rise to the `Jacobian' group cocycle  
 $$
 J: \on{TAut}_{n}\to \hat\Tr_{n}. 
 $$
$J$ is uniquely defined by the conditions $J(\on{id})=0$
and $(d/dt)J(e^{tx}g)_{|t=0} = j(x) + x\cdot J(g)$; as a consequence, 
$J$ satisfies the cocycle identity $J(h\circ g) = J(h) + h\cdot J(g)$. 

The compatibility of $j,J$ with simplicial maps can be described 
as follows. Any partially defined $[m]\supset D_\phi 
\stackrel{\phi}{\to}[n]$ gives rise to a Lie algebra
morphism $\f_n\to \f_m$, $x^\mapsto x^\phi$, with 
$x_k^\phi := \sum_{\ell\in \phi^{-1}(k)}x_\ell$, and any $\tilde\phi$
gives rise to a morphism $\hat\f_n\to\hat\f_m$, $x\mapsto x^{\tilde\phi}$, 
with $x_k^{\tilde\phi} = \on{cbh}(x_\ell|\ell\in\phi^{-1}(k))$. 
These morphisms give rise to linear maps $\Tr_n\to\Tr_m$
and $\hat\Tr_n\to \hat\Tr_m$. Then one can show that 
$j(u^\phi) = j(u)^\phi$, $J(g^\phi) = J(g)^\phi$, 
$j(u^{\tilde\phi}) = j(u)^{\tilde\phi}$, $J(g^{\tilde\phi}) 
= J(g)^{\tilde\phi}$.  

We define a complex $\Tr_{1}\stackrel{\delta}{\to}\Tr_{2}
\stackrel{\delta}{\to}\Tr_{3}...$ by $f(x_{1})\mapsto 
f(x_{1}+x_{2})-f(x_{1})-f(x_{2}) = f^{12}-f^1-f^2$, $f(x_{1},x_{2})\mapsto 
f(x_{1}+x_{2},x_{3})-f(x_{1},x_{2}+x_{3})-f(x_{2},x_{3})+f(x_{1},x_{2})
=f^{12,3}-f^{1,23}-f^{2,3}+f^{1,2}$, 
etc. It is proved in \cite{AT} that this complex is acyclic 
in degree 2 (the degree of $\Tr_i$ is $i$).  
The kernel of $\Tr_{1}\stackrel{\delta}{\to}\Tr_{2}$ is 1-dimensional, 
spanned by the class of $x_{1}\in A_{1}\simeq \Tr_{1}$. 

We similarly define a complex 
$\hat\Tr_{1}\stackrel{\tilde\delta}{\to}\hat\Tr_{2}
\stackrel{\tilde\delta}{\to}\hat\Tr_{3}...$ by $f(x_{1})\mapsto 
f(\on{log}(e^{x_{1}}e^{x_{2}}))-f(x_{1})-f(x_{2}) = f^{\widetilde{12}}-f^1-f^2$, 
$f(x_{1},x_{2})\mapsto 
f(\on{log}(e^{x_{1}}e^{x_{2}}),x_{3})-f(x_{1},
\on{log}(e^{x_{2}}e^{x_{3}}))
-f(x_{2},x_{3})+f(x_{1},x_{2})$. It has a decreasing filtration 
by the degree, and its associated
graded is the above complex, so the complex $\hat\Tr_{1}
\stackrel{\tilde\delta}{\to}...$ is again acyclic in degree 2. Since 
$\on{log}(e^{x_{1}}e^{x_{2}})-x_{1}-x_{2}$ is a sum of brackets, 
$\on{Ker}(\hat\Tr_{1}\stackrel{\tilde\delta}{\to}\hat\Tr_{2})$ 
is again 1-dimensional, spanned by the class of 
$x_{1}\in A_{1}^{\wedge}\simeq \hat\Tr_{1}$. 

\subsection{Braid groups and Lie algebras of infinitesimal braids}

Let $\on{B}_{n}$ be the braid group of order $n$. $\on{B}_{n}$ may be viewed as 
$\pi_{1}(X_{n}/S_{n},S_{n}p)$, where $X_{n} = \{(z_{1},...,z_{n})\in \CC^{n}|
z_{i}\neq z_{i}$ if $i\neq j\}$  and $S_{n}p$ is the $S_{n}$-orbit of the set 
$p=\{(z_{1},...,z_{n})|z_{i}\in\RR, z_{1}<...<z_{n}\}$. 
The fibration $X_{n}\to X_{n}/S_{n}$ gives rise to the morphism $\on{B}_{n}\to S_{n}$, 
and the pure braid group $\on{PB}_{n}$ is defined as $\on{Ker}(\on{B}_{n}\to S_{n})$, so 
we have an exact sequence $1\to \on{PB}_{n}\to \on{B}_{n}\to S_{n}\to 1$; 
also $\on{PB}_{n} = \pi_{1}(X_{n},p)$.  

We recall the Artin presentation of $\on{B}_{n}$: generators are 
$\sigma_{1},...,\sigma_{n-1}$, and relations are given by 
$$\sigma_{i}\sigma_{i+1}
\sigma_{i}=\sigma_{i+1}\sigma_{i}\sigma_{i+1} \quad (i=1,...,n-2),
\quad  
\sigma_{i}\sigma_{j} =\sigma_{j}\sigma_{i} \on{\ for\ } |i-j|>1.$$ 
We also recall the Coxeter presentation of $S_{n}$: generators
are $s_{1},...,s_{n-1}$ ($s_{i}$ is the permutation $(i,i+1)$)
and relations are the same as those between the 
$\sigma_{i}$, with the additional relations $s_{i}^{2}=1$ ($i=1,...,n-1$). 
The morphism $\on{B}_{n}\to S_{n}$ is then given by $\sigma_{i}\mapsto s_{i}$. 

\begin{figure}[h!]
\begin{center}
\includegraphics[width=14cm]{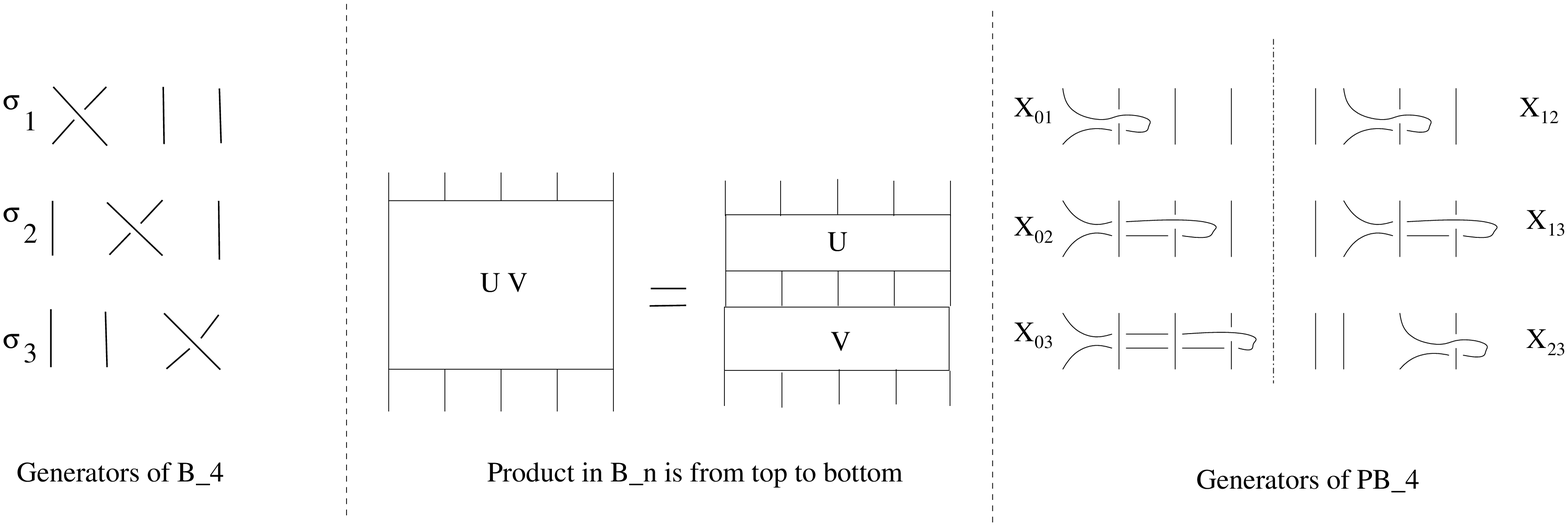}
\caption{\footnotesize }\label{AETfig1.eps}
\end{center}
\end{figure}

The group $\on{PB}_{n}$ admits the following presentation. For $i<j$
($i,j\in [n]$), set 
$$
x_{ij}:= (\sigma_{j-2}...\sigma_{i})^{-1}\sigma_{j-1}^{2} 
(\sigma_{j-2}...\sigma_{i}). 
$$
The generators $x_{ij}$ belong to $\on{PB}_{n}$, and\footnote{We set
$(g,h):= ghg^{-1}h^{-1}$.} 
$$
(x_{ij}x_{ik}x_{jk},x_{ij}) = 
(x_{ij}x_{ik}x_{jk},x_{ik}) = 
(x_{ij}x_{ik}x_{jk},x_{jk}) = 1\on{\ for\ }i<j<k, 
$$
and 
$$
(x_{ij},x_{kl})=(x_{il},x_{jk}) = (x_{ik},x_{jk}x_{jl}x_{jk}^{-1})=1 
\on{\ for\ }i<j<k<l. $$
One proves that this constitutes a presentation of $\on{PB}_{n}$, see 
Figure \ref{AETfig1.eps}. 

For any sequence $(k_{1},...,k_{n})$ of integers $\geq 0$, 
there exists a unique morphism $\on{PB}_{n}\to \on{PB}_{k_{1}+...+k_{n}}$
consisting in replacing the first strand by $k_{1}$ consecutive 
strands, ..., the $n$th strand by $k_{n}$ consecutive strands. 
If we set $m:= k_{1}+...+k_{n}$ and $\phi:[m]\to [n]$ is the map 
such that $\phi(k_{1}+...+k_{i-1}+[k_{i}])=i$, we denote this morphism 
by $x\mapsto x^{\tilde\phi} = x^{\widetilde{1...k_{1}},...,
\widetilde{k_{1}+...+k_{n-1}+1...m}}$. This morphism in explicitly given by 
$$
x_{ij}\mapsto \prod_{i'\in \phi^{-1}(i)}^{\nearrow} (\prod_{j'\in 
\phi^{-1}(j)}^{\searrow} x_{i'j'}), 
$$
where $\prod^{\nearrow}, \prod^{\searrow}$ mean the product 
in increasing and decreasing order of the indices.  

The Lie algebra $\t_{n}$ of infinitesimal braids is presented by generators
$t_{ij}$, $i\neq j\in [n]$ and relations $t_{ji}=t_{ij}$, 
$[t_{ij},t_{ik}+t_{jk}]=0$ for $i,j,k$ distinct and 
$[t_{ij},t_{kl}]=0$ for $i,j,k,l$ distinct. For each partially defined map 
$[m]\supset D_{\phi}\stackrel{\phi}{\to} [n]$, there is a unique 
Lie algebra morphism $\t_{n}\to \t_{m}$, $x\mapsto x^{\phi}$
given by $t_{ij}^{\phi}:= \sum_{i'\in\phi^{-1}(i),j'\in\phi^{-1}(j)}t_{i'j'}$
(in particular, we have an action of $S_{n}$ on $\t_{n}$). 
We often write $x^{\phi^{-1}(1),...,\phi^{-1}(n)}$ instead of $x^{\phi}$. 
We attribute degree $1$ to each of the generators $t_{ij}$, so the
Lie algebra $\t_{n}$ is positively graded; we denote by $\hat\t_{n}$
its degree completion.  

\subsection{The morphism $\t_{n+1}\to \tder_{n}$}

Let us reindex $t_{ij}$, $i\neq j\in \{0,...,n\}$ the generators of $\t_{n+1}$. 
One checks that there is a unique morphism $\on{ad}:\t_{n+1}\to \tder_{n}$, 
defined by $t_{0i}\mapsto (x_{j}\mapsto [x_{i},x_{j}])$ and 
$t_{ij}\mapsto (x_{i}\mapsto [x_{i},x_{j}], x_{j}\mapsto 
[x_{j},x_{i}], x_{k}\mapsto 0$ for $k\neq i,j)$ if $i,j\neq 0$. 
It exponentiates to 
$\on{Ad}:\on{exp}(\hat\t_{n+1})\to \on{TAut}_{n}$. One checks that 
$j(\on{ad}t_{ij})=0$, so the cocycle property implies $j(\on{ad}x)=
J(\on{Ad}X)=0$
for any $x\in\t_{n+1}$ and $X\in \on{exp}(\hat\t_{n+1})$. 

The morphism $\on{ad}:\t_{n+1}\to \tder_{n}$ may be interpreted as follows. 
The Lie subalgebra of $\t_{n+1}$ generated by the elements 
$t_{0i}$, $i\in [n]$
identifies with $\f_{n}$ under $x_{i}\mapsto t_{0i}$; it is a Lie ideal of
$\t_{n+1}$. Then $\on{ad}:\t_{n+1}\to \on{Der}(\f_{n})$
can be viewed as the adjoint action of $\t_{n+1}$ on its Lie ideal 
$\f_{n}\subset \t_{n+1}$. 

Note that the morphism $\t_{n}\to \t_{n+1}$, $t_{ij}\mapsto t_{ij}$
is injective, so $\t_{n}$ may be viewed as a Lie subalgebra of $\t_{n+1}$; 
then $\t_{n+1}$ identifies with the semidirect product 
$\f_{n}\rtimes_{\on{ad}}\t_{n}$.

\subsection{The morphism $\on{PB}_{n+1}\to \on{TAut}_{n}$}

Reindex the generators of $\on{PB}_{n+1}$ as $x_{ij}$, $i<j\in \{0,...,n\}$. 
Let $\on{F}_{n}$ be the free group with generators $X_{i}$ ($i\in [n]$). 
Then: (a) the morphism $\on{F}_{n}\to \on{PB}_{n+1}$, $X_{i}\mapsto 
x_{0i}$, is injective; (b) $\on{F}_{n}$ is a normal subgroup in $\on{PB}_{n+1}$. 
This implies that we have an action $\on{Ad}:\on{PB}_{n+1}\to 
\on{Aut}(\on{F}_{n})$ of $\on{PB}_{n}$ by automorphisms of $\on{F}_{n}$. 

This action can be made explicit as follows: 
if $i>0$, then 
$$
\on{Ad}(x_{0i})(X_{j})=X_{i}X_{j}X_{i}^{-1}, 
$$
and if $0<i<j$, then 
\begin{equation*}
\begin{split}
& \on{Ad}(x_{ij})(X_{i})=
X_{j}^{-1}X_{i}X_{j}, \quad \on{Ad}(x_{ij})(X_{j})=
(X_{i}X_{j})^{-1}X_{j}(X_{i}X_{j}), \\
 &  
\on{Ad}(x_{ij})(X_{k})=X_{k} \quad 
\on{for} \quad k<i\quad  \on{or} \quad k>j, \\
 &  
\on{Ad}(x_{ij})(X_{k})=(X_{j}^{-1}X_{i}^{-1}X_{j}X_{i})
X_{k}(X_{j}^{-1}X_{i}^{-1}X_{j}X_{i})^{-1}\quad  \on{for} \quad i<k<j. 
\end{split}\end{equation*}
This extends to an action of $\on{PB}_{n+1}$ by automorphisms of 
$\on{F}_{n}(\kk)$. Using the isomorphism $\on{F}_{n}(\kk)\simeq 
\on{exp}(\hat\f_{n})$ given by $X_{i}\mapsto e^{x_{i}}$, we therefore
obtain a morphism $\on{PB}_{n}\to \on{Aut}(\hat\f_{n})$. Its image is 
contained in $\on{TAut}_{n}$ (since $\on{Ad}x_{ij}$ belongs to 
this subgroup and the elements $x_{ij}$ generate $\on{PB}_{n}$), and 
since $\on{TAut}_{n}$ is prounipotent, the universal property of 
Malcev completions implies that $\on{Ad}$ extends to a morphism
$\on{Ad}:\on{PB}_{n}(\kk)\to \on{TAut}_{n}$. 

\begin{lemma}
For any $g\in \on{PB}_{n+1}(\kk)$, $J(\on{Ad}g)=0$. 
\end{lemma}

{\em Proof.} It suffices to show that $J(\on{Ad}x_{ij})=0$. For any 
$u\in \on{F}_{n}(\kk)$, $J(\on{Inn}u)=0$ (where $\on{Inn}u$ is
$v\mapsto uvu^{-1}$) and $\on{Ad}(x_{0i}) = \on{Inn}(X_{i})$, 
so it suffices to prove that 
$J(\on{Inn}X_{j} \circ \on{Ad}x_{ij})=0$ for $0<i<j$. Let $\theta_{ij}:=
\on{Inn}X_{j}\circ\on{Ad}x_{ij}$, then $\theta_{ij} : 
X_{i}\mapsto X_{i}$, $X_{j}\mapsto X_{i}^{-1}X_{j}X_{i}$, 
$X_{k}\mapsto X_{j}X_{k}X_{j}^{-1}$ for $k<i$ or $k>j$, 
$X_{k}\mapsto (X_{i}^{-1}X_{j}X_{i})X_{k} (X_{i}^{-1}X_{j}X_{i})^{-1}$
for $i<k<j$.

Let $\u\subset \tder_n^\wedge$ be the subspace 
of all elements $\lbr a_1,...,a_n\rbr$, where
$a_j\in \kk x_i$, $a_i=0$, and for $k\neq i,j$, 
$a_k\in \hat\f_n$ has the form $a_k(x_i,x_j)$ ($a_k\in \hat\f_2$). 
This is a Lie subalgebra in $\tder_n^\wedge$, so $\on{exp}$ maps
it bijectively to a subgroup of $\on{TAut}_n$. One checks that
$\on{exp}(\u)\subset U$, where $U\subset \on{TAut}_n$ is the 
subspace of all $\lbr U_1,...,U_n\rbr$, where $U_j\in \{e^{\lambda x_i},
\lambda\in\kk\}$, $U_i=1$, and for $k\neq i,j$, $U_k$ has the form
$U_k(x_i,x_j)$ ($U_k\in\on{exp}(\hat\f_2)$), and that $U$ is an 
algebraic subgroup of $\on{TAut}_n$. Therefore $\u\subset \on{Lie}(U)$. 
On the other hand, one checks that $\u$ coincides with the 
tangent subspace of $U$ at the origin, so $\u = \on{Lie}(U)$. 
It follows that $\on{log}$ takes $U$ to $\u$.   

All this implies that $\log\theta_{ij}$ has the form
$\lbr a_{1},...,a_{n}\rbr$, where $a_{i}=0$, $a_{j}=-x_{i}$ and for 
$k\neq i,j$, $a_{k}\in\hat\f_{n}$ has the form $a_{k}(x_{i},x_{j})$. Then 
$j(\on{log}\theta_{ij})=0$, hence $J(\theta_{ij})=0$, as wanted. 
\hfill \qed\medskip 

Note that the quotient group $\on{PB}_{n+1}/\on{F}_{n}$ identifies with 
$\on{PB}_{n}$ under $x_{ij}\mapsto x_{ij}$ for $0<i<j$, 
$x_{0i}\mapsto 1$. We then have an exact sequence $1\to \on{F}_{n}\to 
\on{PB}_{n+1}\to \on{PB}_{n}\to 1$. Moreover, this exact sequence 
admits the splitting $\on{PB}_{n}\to \on{PB}_{n+1}$, $x_{ij}\mapsto x_{ij}$. 
It follows that $\on{PB}_{n+1}$ identifies with the semidirect product 
$\on{F}_{n}\rtimes_{\on{Ad}}\on{PB}_{n}$. 

\begin{remark}
We will rename $x,y$ (resp., $x,y,z$, $X,Y$, $X,Y,Z$) the generators 
$x_{1},x_{2}$ (resp., $x_{1},x_{2},x_{3}$, $X,Y$, $X,Y,Z$) of $\hat\f_{2}$
(resp., $\hat\f_{3}$, $\on{F}_{2}$, $\on{F}_{3}$). 
\end{remark}

\section{The main results} \label{sec:result}

\subsection{The map $M_{1}(\kk)\to \on{SolKV}(\kk)$}

Let $\hat\f_{2}$ be the topologically free Lie algebra generated by $x,y$. 
Let $\on{F_{2}}$ be the free group with generators $X,Y$ and let 
$\on{F}_{2}(\kk)$ be its prounipotent completion; we have an identification 
$\on{F}_{2}(\kk) \simeq \on{exp}(\hat\f_{2})$, induced by the morphism 
$\on{F}_{2}\to \on{exp}(\hat\f_{2})$ given by $X\mapsto e^{x},Y\mapsto e^{y}$. 

The set of solutions of the Kashiwara--Vergne equations 
is (see \cite{KV,AT})\footnote{For $g,h$ in a prounipotent group 
$G$ or its Lie algebra, we use the notation $g\sim h$ for
`$g$ is conjugated to $h$', i.e., $g=khk^{-1}$ for some $k\in G$.} 
\footnote{If $\Gamma$ is a finitely generated
group, we denote by $\Gamma(\kk)$ its prounipotent (of Malcev) completion. 
There is a group morphism $\Gamma\to\Gamma(\kk)$ with the universal property that 
any group morphism $\Gamma\to U$, with $U$ prounipotent, extends uniquely to a 
morphism $\Gamma(\kk)\to U$ of algebraic groups.} 
\footnote{The definition given here is equivalent to that of \cite{AT} 
as $\Tr_{1}\to\Tr_{2}\to\Tr_{3}$ is acyclic.}
\begin{equation*}\begin{split}
\on{SolKV}(\kk) := & \{\mu\in\on{Iso}(\on{F}_{2}(\kk),\on{exp}(\hat\f_{2})) | 
\mu(X)\sim e^{x}, \mu(Y)\sim e^{y},\mu(XY)=e^{x+y}, \\
 & \on{and\ }\exists r\in u^{2}\kk[[u]] | J(\mu) = 
 \langle r(x+y) - r(x) - r(y)\rangle\}. 
\end{split}\end{equation*}
Here $\mu$ gives rise to an element of $\on{TAut}_{2}$ (using 
$\on{F}_{2}(\kk) \simeq \on{exp}(\hat\f_{2})$) and $J(\mu)$ is its Jacobian. 
As the kernel of $\Tr_{1}\to\Tr_{2}$ is equal to $\kk u$, $r$ is uniquely 
determined by $\mu\in\on{SolKV}(\kk)$, so we define a map 
$\on{Duf} : \on{SolKV}(\kk)\to u^{2}\kk[[u]]$, $\mu\mapsto 
r = \on{Duf}(\mu)$; we will call $r$ the Duflo formal series of $\mu$. 

The set of associators with coupling constant 1 is 
\begin{equation*} \begin{split}
& M_{1}(\kk):= \{\Phi(t_{12},t_{23})\in \on{exp}(\hat\t_{3})|
\Phi^{3,2,1}=\Phi^{-1}, 
e^{t_{23}/2}\Phi^{1,2,3}e^{t_{12}/2}\Phi^{3,1,2}e^{t_{31}/2}\Phi^{2,3,1}=
e^{(t_{12}+t_{23}+t_{31})/2}, \\
& \Phi^{2,3,4}\Phi^{1,23,4}\Phi^{1,2,3}=\Phi^{1,2,34}
\Phi^{12,3,4}\}.
\end{split}\end{equation*}

\begin{theorem} \label{thm:main}
There is a unique map $M_{1}(\kk)\to \on{SolKV}(\kk)$, $\Phi\mapsto \mu_{\Phi}$, 
such that\footnote{If $G$ is a prounipotent group, we use the 
notation $g\cdot h\cdot (same)^{-1}$ for $ghg^{-1}$ for if $g\in G$
and $h\in G$ or $\on{Lie}(G)$.} 
$$
\mu_{\Phi}(X) = \Phi(x,-x-y) e^{x}\Phi(x,-x-y)^{-1}, 
\mu_{\Phi}(Y)= 
e^{-(x+y)/2}\Phi(y,-x-y)\cdot e^{y} \cdot (same)^{-1}. 
$$
\end{theorem}

The Jacobian of $\mu_{\Phi}$ can be computed as follows. 
In \cite{DT,E} (see also \cite{Ih}), it was proved\footnote{The key 
ingredient in the proof of this result is the statement that 
the image of $\grt_{1}$ in $\f_{2}'/\f_{2}''$ is spanned 
by the classes of the Drinfeld generators. This statement  
also follows from 
Theorem 4.1 in \cite{AT}; indeed, one sees easily that
the diagram 
$$
\begin{matrix}
& \f_{2} & \stackrel{\psi\mapsto \langle a\partial_{a}\psi\rangle}{\to} & 
\Tr_{2} & \stackrel{\phi\mapsto \phi^{ab}}{\to} &  \kk[\bar a,\bar b] \\
& \uparrow & & & & \uparrow \\
  & \f_{2}' & &\to & & \f_{2}'/\f_{2}'' = \bar a\bar b\kk[\bar a,\bar b]\\
& \uparrow & & & & \uparrow \\
 & \grt_{1} & & \to & & \grt_{1}/\grt_{1}' 
\end{matrix}$$
commutes (the upper part follows from the fact that $\f_{2}'$
is freely generated by the $(\on{ad}a)^{k}(\on{ad}b)^{l}([a,b])$ and the
bottom part from $\grt_{1}'\subset \f_{2}''$); 
Theorem 4.1 in \cite{AT} implies that the image of $\grt_{1}\to\Tr_{2}$
is spanned by the images of the Drinfeld generators; it follows that the same is
true of the image of $\grt_{1}\to\f_{2}'/\f_{2}''$.} 
 that for any $\Phi(a,b)
\in M_{1}(\kk)$, 
there exists a formal series $\Gamma_{\Phi}(u) = e^{\sum_{n\geq 2}
(-1)^n\zeta_{\Phi}(n)u^{n}/n}$, such that 
\begin{equation} \label{Phi:Gamma}
(1+b\partial_{b}\Phi(a,b))^{ab} = {\Gamma_{\Phi}(\overline a + \overline b)
\over {\Gamma_{\Phi}(\overline a)\Gamma_{\Phi}(\overline b)}}, 
\end{equation}
where $\partial_{b}\Phi(a,b)$ is defined as above and $x\mapsto x^{ab}$
is the abelianization morphism $\kk\langle\langle a,b\rangle\rangle
\to \kk[[\overline a,\overline b]]$ 

The values of the $\zeta_{\Phi}(n)$
for $n$ even are independent of $\Phi$, given by 
$-{1\over 2}({u\over{e^u-1}}-1+{u\over 2}) = \sum_{n\geq 1}
\zeta_\Phi(2n)u^{2n}$, so they are related to Bernoulli numbers
by $\zeta_\Phi(2n) = -{1\over 2}{{B_{2n}}\over{(2n)!}}$ for $n\geq 1$ 
(we have $\zeta_{\Phi}(2)=-1/24$, $\zeta_{\Phi}(4)=1/1440$, etc.)

\begin{proposition} \label{form:J}
$J(\mu_{\Phi}) = \langle \on{log}\Gamma_{\Phi}(x) 
+ \on{log}\Gamma_{\Phi}(y)- \on{log}\Gamma_{\Phi}(x+y) \rangle$, so 
$\on{Duf}(\mu_{\Phi}) = -\on{log}\Gamma_{\Phi}$. 
\end{proposition}

\subsection{Torsor aspects}

We set 
\begin{equation*} \begin{split}
\on{KV}(\kk) := \{ & \alpha\in\on{Aut}(\on{F}_{2}(\kk)) | \alpha(X)\sim X,\alpha(Y)
\sim Y, \alpha(XY)=XY, \\
 & \on{and\ }\exists \sigma\in u^{2}\kk[[u]] | 
 J(\alpha) = \langle \sigma(\on{log}(e^{x}e^{y})) - \sigma(x)-\sigma(y)\rangle\}
\end{split}\end{equation*}
and 
\begin{equation*} \begin{split}
\on{KRV}(\kk) := \{ & a\in\on{Aut}(\hat\f_{2}) | a(x)\sim x, a(y)\sim y, a(x+y)
= x+y, \\
 & \on{and\ }\exists s\in u^{2}\kk[[u]] | 
 J(a) = \langle s(x+y) - s(x)-s(y)\rangle\}. 
\end{split}\end{equation*}
Here $\alpha,a$ give rise to elements of $\on{TAut}_{2}$ (using 
$\on{F}_{2}(\kk) \simeq \on{exp}(\hat\f_{2})$) and $J(\alpha),J(a)$ 
are their Jacobians. As before, we will denote $\on{Duf}:\on{KV}(\kk)
\to u^{2}\kk[[u]]$, $\on{KRV}(\kk)\to u^{2}\kk[[u]]$ the maps 
$\alpha\mapsto\sigma$, $a\mapsto s$. 

\begin{proposition} \label{SolKV:torsor}
$\on{KV}(\kk)$ and $\on{KRV}(\kk)$ are groups. 
$\on{SolKV}(\kk)$ is a torsor under the commuting left action of $\on{KV}(\kk)$
and right action of $\on{KRV}(\kk)$ given by $(\alpha,\mu)\mapsto \mu\circ
\alpha^{-1}$ and $(\mu,a)\mapsto a^{-1}\circ\mu$. 
\end{proposition}

In particular, any element of $\on{SolKV}(\kk)$ gives rise to an isomorphism
${{\mathfrak{kv}}}\to {{\mathfrak{krv}}}$ between the Lie algebras of these groups, 
whose associated graded is the canonical identification 
$\on{gr}({{\mathfrak{kv}}})
\simeq{{\mathfrak{krv}}}$. 

The prounipotent radical of the Grothendieck-Teichm\"uller group is 
\begin{equation*} \begin{split}
& \on{GT}_{1}(\kk) = \{f(X,Y)\in \on{F}_{2}(\kk) | f(Y,X)=f(X,Y)^{-1}, 
f(X,Y)f(Y^{-1}X^{-1},X)f(Y,Y^{-1}X^{-1})=1, \\
& f(x_{23},x_{34})f(x_{12}x_{13},x_{24}x_{34})f(x_{12},x_{23})=
f(x_{12},x_{23}x_{24})f(x_{13}x_{23},x_{34})\}
\end{split}\end{equation*}
(the last equation is in $\on{PB}_{4}(\kk)$) with product
$(f_{1}*f_{2})(X,Y)=f_{1}(f_{2}(X,Y)Xf_{2}(X,Y)^{-1},Y)f_{2}(X,Y)$. 
Its graded version is
$$
\on{GRT}_{1}(\kk) = \{g(t_{12},t_{23})\in \on{exp}(\hat\t_{3})|g^{3,2,1}
=g^{-1}, 
g(A,C)A g(A,C)^{-1}+g(B,C)Bg(B,C)^{-1}+C=0 $$ 
$$\on{\ for\ } A+B+C=0, 
g^{1,2,3}g^{3,1,2}g^{2,3,1}=1,g^{2,3,4}g^{1,23,4}g^{1,2,3}
=g^{1,2,34}g^{12,3,4}\}$$
with product $(g_{1}*g_{2})(a,b)=g_{1}(g_{2}(a,b)ag_{2}(a,b)^{-1},b)g_{2}(a,b)$
(we set $a:=t_{12},b:= t_{23}$). 

\begin{proposition} (see \cite{Dr:Gal})
$M_{1}(\kk)$ is a torsor under the commuting left action of $\on{GT}_{1}(\kk)$
and right action of $\on{GRT}_{1}(\kk)$ by $(f,\Phi)\mapsto 
(f*\Phi)(a,b):=f(\Phi(a,b)e^{a}\Phi(a,b)^{-1},e^{b}) \Phi(a,b)$
and $(\Phi,g)\mapsto (\Phi*g)(a,b):= 
\Phi(g(a,b)a g(a,b)^{-1},b)g(a,b)$. 
\end{proposition}

The following Theorem \ref{thm:morph:tors} and Proposition 
\ref{form:J:torsor} express torsor properties of the map 
$\Phi\mapsto \mu_{\Phi}$. 

\begin{theorem} \label{thm:morph:tors}
There are unique group morphisms $\on{GT}_{1}(\kk)\to \on{KV}(\kk)$, 
$f(X,Y)\mapsto \alpha_{f}^{-1}$, where
$$
\alpha_{f}(X) = f(X,Y^{-1}X^{-1})X f(X,Y^{-1}X^{-1})^{-1}, \alpha_{f}(Y)=
 f(Y,Y^{-1}X^{-1})Y f(Y,Y^{-1}X^{-1})^{-1}, $$
and $\on{GRT}_{1}(\kk)\to \on{KRV}_{1}(\kk)$, $g(a,b)\mapsto a_{g}^{-1}$, where 
$$
a_{g}(x)=g(x,-x-y)xg(x,-x-y)^{-1},a_{g}(y)=g(y,-x-y)yg(y,-x-y)^{-1}. $$
These group morphisms are compatible with the map 
$M_{1}(\kk)\to \on{SolKV}(\kk)$, 
which is therefore a morphism of torsors. 
\end{theorem}

\begin{proposition} \label{form:J:torsor}
We have a commuting diagram of torsors 
$$
\begin{matrix}
M_{1}(\kk) &\stackrel{\Phi\mapsto\mu_{\Phi}}{\to} & \on{SolKV(\kk)}\\
\scriptstyle{\Phi\mapsto \on{log}\Gamma_{\Phi}} 
\downarrow& & \downarrow\scriptstyle{\on{Duf}} \\
\{r\in u^{2}\kk[[u]] | r_{ev}(u) = -{u^{2}\over{24}}
+{u^{4}\over{1440}}...\} & 
\stackrel{(-1)\times -}{\hookrightarrow}& 
u^{2}\kk[[u]]
\end{matrix}$$
where $r_{ev}(u)$ is the even part of $r(u)$, and the 
spaces in the lower line are viewed as affine spaces. 
\end{proposition}

\subsection{Analytic aspects}

Let us recall the original form of the KV conjecture. Let $\kk = \RR$ or $\CC$. 

\begin{conjecture} (\cite{KV})
For any finite dimensional $\kk$-Lie algebra $\G$, there exists a 
pair of Lie series $A(x,y),B(x,y)\in \hat\f_{2}$, such that: 

(KV1) $x+y-\on{log}e^{y}e^{x} = (1-e^{-\on{ad}x})(A(x,y)) 
+ (e^{\on{ad}y}-1)(B(x,y))$; 

(KV2) $A,B$ give convergent power series at the neighborhood 
of $(0,0)\in\G^{2}$; 

(KV3) $\on{tr}_{\G}((\on{ad}x)\partial_{x}A + 
(\on{ad}y)\partial_{y}B) = {1\over 2}\on{tr}_{\G}
({{\on{ad}x}\over{e^{\on{ad}x}-1}} 
+ {{\on{ad}y}\over{e^{\on{ad}y}-1}}
- {{\on{ad}z}\over{e^{\on{ad}z}-1}}-1)$
(identity of analytic functions on $\G^{2}$ near the origin), 
where $z = \on{log}e^{x}e^{y}$ and for $(x,y)\in \G^{2}$, 
$(\partial_{x}A)(x,y)\in \on{End}(\G)$ is 
$a\mapsto {d\over{dt}}_{|t=0}A(x+ta,y)$, 
$(\partial_{y}B )(x,y)(a) =  {d\over{dt}}_{|t=0}B(x,y+ta)$. 
\end{conjecture}

According to \cite{AT}, there is a unique map 
$\kappa:\on{TAut}_{2}\to \tder_{2}$, where $\kappa(g):= 
\ell - g\ell g^{-1}$, and $\ell\in \on{Der}(\hat\f_{2})$ is the `grading'
derivation $\ell(x_{i})=x_{i}$. It is proved in \cite{AT}
that if $\mu\in\on{SolKV}(\kk)$, and $(A,B)$ are such that 
$-\kappa(\mu^{-1}) = \lbr A,B\rbr$, then (KV1) and (KV3)
hold as identities between formal series for any $\G$, where in (KV3) the formal 
series ${1\over 2}{t\over{e^{t}-1}}$ is replaced by $r_{\mu}(t)$. 

Let $\Phi_{\on{KZ}}(a,b)\in \on{exp}(\hat\f_{2})$ be the 
KZ associator, and $\tilde\Phi_{\on{KZ}}(a,b):= \Phi_{\on{KZ}}
({a\over{2\pi\i}},{b\over{2\pi\i}})$; recall that $\tilde\Phi_{\on{KZ}}$
is the renormalized holonomy from $0$ to $1$ of $G'(t) = 
{1\over{2\pi\i}}({a\over t} + {b\over {t-1}})G(t)$, and 
$\tilde\Phi_{\on{KZ}}\in M_{1}(\CC)$. Set 
$\mu_{\on{KZ}}:= \mu_{\tilde\Phi_{\on{KZ}}}$ and 
$u_{\on{KZ}} = \lbr A_{\on{KZ}},B_{\on{KZ}}\rbr
:= -\kappa(\mu_{\on{KZ}}^{-1})$. 

Let $(A_{\RR},B_{\RR})$ be defined as the real parts of $(A_{\on{KZ}},B_{\on{KZ}})$
(w.r.t. the natural real structure of $\hat\f_{2}$). 
Then: 

\begin{theorem} \label{thm:an}
1) $(A_{\RR},B_{\RR})$ satisfies (KV1), (KV2)
and (KV3) for any finite dimensional Lie algebra $\G$ and is therefore 
a universal solution of the KV conjecture. 

2) For any $s\in\RR$, $(A_{s},B_{s}):= (A_{\RR}+s(\on{log}(e^{x}e^{y})-x),
B_{\RR}+s(\on{log}(e^{x}e^{y})-y))$ is a universal solution of the KV conjecture. 

3) When $s=-1/4$, we have $(A_{s}(x,y),B_{s}(x,y))=(B_{s}(-y,-x),A_{s}(-y,-x))$. 
\end{theorem}

Of course, the main new result here is the analyticity statement (KV2). 

\subsection{Organization of the proofs}
We construct the isomorphisms $\tilde\mu_{\Phi}^{O}$ and 
$\mu_{\Phi}^{O}$ in Section \ref{sec:assoc}.
In Section \ref{sec:mu:mu}, we prove the identity relating $\mu_{O}$ and 
$\mu_{O^{(i)}}$. We then prove Theorem \ref{thm:main} and Proposition 
\ref{form:J} in Section \ref{pf:main}. In Section \ref{sec:torsor}, we 
prove Proposition \ref{SolKV:torsor}, Theorem \ref{thm:morph:tors} 
and Proposition \ref{form:J:torsor}. 
Section \ref{act:gt} is devoted to a direct proof of the properties of $\alpha_{f}$. 
In Section \ref{sec:jac}, we compute the Jacobians of $\mu_{\Phi}^{O}$ 
and $\alpha_{f}^{O}$ and in Section \ref{sec:anal}, we prove the analytic 
Theorem \ref{thm:an}. Appendix \ref{app} is devoted to  
results on centralizers in $\t_{n}$ and $\on{PB}_{n}(\kk)$. 

\section{Associators and isomorphisms of free groups} \label{sec:assoc} 

\subsection{The categories ${\bf PaB},{\bf PaCD}$}

In \cite{B}, Bar-Natan introduced the category ${\bf{PaB}}$
of parenthesized braids. Its set of objects is the set of pairs $O = (n,P)$, 
where $n$ in an integer $\geq 0$ and $P$ is a parenthesization of the word
$\bullet...\bullet$ ($n$ letters); alternatively, $P$ is a planar binary tree
with $n$ leaves (we will set $|O|=n$). The object with $n=0$ is denoted 
${\bf 1}$. The morphisms 
are defined by ${\bf{PaB}}(O,O') = \emptyset$ if $|O|\neq |O'|$, and 
$=\on{B}_{n}$ if $|O|=|O'|=n$; the composition is then defined using the product in 
$\on{B}_{n}$. 

\begin{figure}[h!]
\begin{center}
\includegraphics[width=10cm]{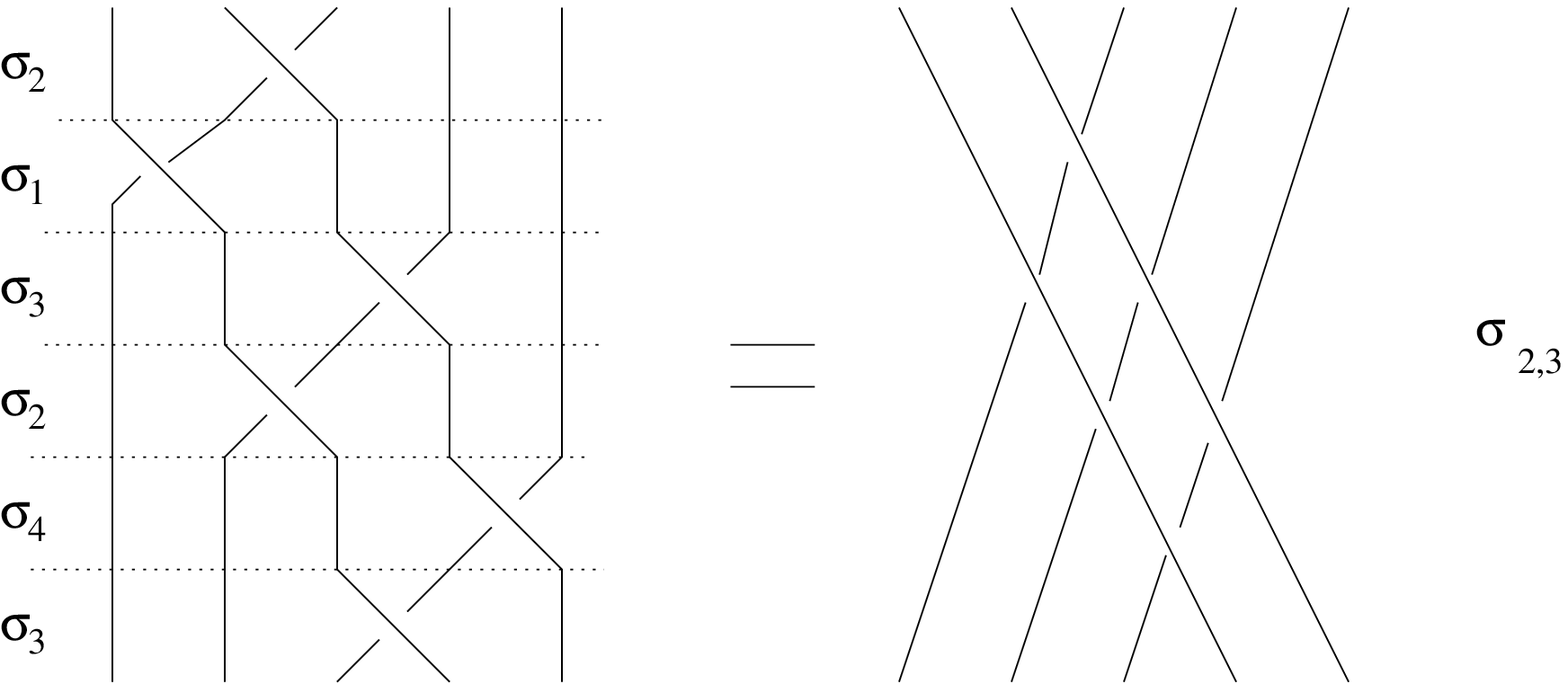}
\caption{\footnotesize }\label{AETfig2.eps}
\end{center}
\end{figure}

${\bf{PaB}}$ is a braided monoidal category (see e.g. \cite{CE}), where the 
tensor product of objects is $(n,P)\otimes (n',P'):= 
(n+n',P*P')$ (where $P*P'$ is the concatenation of parenthesized words, e.g.
for $P=\bullet\bullet$ and $P'=(\bullet\bullet)\bullet$, 
$P*P'=(\bullet\bullet)((\bullet\bullet)\bullet)$). The tensor 
product of morphisms ${\bf{PaB}}(O_{1},O'_{1}) \times 
{\bf{PaB}}(O_{2},O'_{2}) \to
{\bf{PaB}}(O_{1}\otimes O_{2},O'_{1}\otimes O'_{2})$ is induced by the 
juxtaposition of braids $\on{B}_{|O_{1}|}\times \on{B}_{|O_{2}|}
\to \on{B}_{|O_{1}|+|O_{2}|}$ (the group morphism 
$(\sigma_{i},e)\mapsto \sigma_{i}$, $(e,\sigma_{j})\mapsto 
\sigma_{j+|O_{1}|}$). The braiding $\beta_{O,O'}
\in {\bf{PaB}}(O\otimes O',O'\otimes O)$
is the braid $\sigma_{n,n'}\in \on{B}_{n+n'}$ where the 
$n$ first strands are globally exchanged with the $n'$ last strands
(see Figure \ref{AETfig2.eps}); we have $\sigma_{n,n'} = (\sigma_{n}...\sigma_{1})
(\sigma_{n+1}...\sigma_{2})...(\sigma_{n+n'-1}...\sigma_{n'})$ (where $n=|O|$,
$n'=|O'|$). 
Finally, the associativity constraint $a_{O,O',O''}\in 
{\bf{PaB}}((O\otimes O')\otimes O'',O\otimes (O'\otimes O''))$ 
corresponds to the trivial braid $e\in B_{|O|+|O'|+|O''|}$. 

Moreover, the pair $({\bf{PaB}},\bullet)$ is universal 
for pairs $(\cC,M)$ of a braided monoidal 
category and an object, i.e., for each such pair, there exists a unique tensor functor
${\bf{PaB}}\to \cC$ taking $\bullet$ to $M$. 

Bar-Natan introduced another category ${\bf PaCD}$, 
which we will describe as follows. Its set of objects is the same as that of 
${\bf PaB}$, and ${\bf PaB}(O,O')=\emptyset$ if $|O|\neq |O'|$, 
$=\on{exp}(\hat\t_{n})\rtimes S_{n}$ if $|O|=|O'|=n$. We define 
the tensor product as above at the level of objects, and by the juxtaposition 
map 
$(\on{exp}\hat\t_{n}\rtimes S_{n}) \times (\on{exp}\hat\t_{n}\rtimes S_{n}) \to
\on{exp}\hat\t_{n+n'}\rtimes S_{n+n'}$ (the group morphism induced by 
$((t_{ij},1),1)\mapsto t_{ij}$, $((1,s_{i}),1)\mapsto s_{i}$, 
$(1,(t_{ij},1))\mapsto t_{n+i,n+j}$, $(1,(1,s_{i}))\mapsto s_{n+i}$)
at the level of morphisms. 

Any $\Phi\in M_{1}(\kk)$ gives rise to a structure of braided monoidal 
category on ${\bf PaCD}$ (and therefore to a tensor functor ${\bf PaB}\to
{\bf PaCD}$, which is the identity at the level of objects) as follows: 
$\beta_{O,O'} = e^{\sum_{i=1}^{n}\sum_{j=n+1}^{n+n'}t_{ij}/2}s_{n,n'}$, 
where $n=|O|,n'=|O'|$, and $s_{n,n'}\in S_{n+n'}$ is given by 
$s_{n,n'}(i)=n'+i$ for $i\in [n]$, $s_{n,n'}(n+i)=i$
for $i\in [n']$, and $a_{O,O',O''} = \Phi(t_{12},t_{23})^{1...n,n+1...n+n',
n+n'+1...n+n'+n''}$ for $n=|O|$, $n'=|O'|$, $n''=|O''|$.

\subsection{Morphisms $\on{B}_{n}\to \on{exp}(\hat\t_{n})\rtimes S_{n}$, 
$\on{PB}_{n}\to \on{exp}(\hat\t_{n})$}

Fix $\Phi\in M_{1}(\kk)$. It gives rise to a functor $F_{\Phi}:
{\bf PaB}\to {\bf PaCD}$, 
so for any $n\geq 1$ and any $O\in \on{Ob}({\bf PaB})$, $|O|=n$, we get a 
group morphism 
$$
F_{\Phi}(O) = \tilde\mu_{O} : \on{B}_{n}\simeq {\bf PaB}(O)
\to {\bf PaCD}(O)=\on{exp}(\hat\t_{n})\rtimes S_{n},  
$$
such that  
$\begin{matrix}
\on{B}_{n} & \stackrel{\tilde\mu_{O}}{\to} & \on{exp}(\hat\t_{n})\rtimes S_{n}\\
 \searrow & & \swarrow \\
  & S_{n} & \end{matrix}$ commutes. It follows that $\tilde\mu_{O}$
restricts to a morphism 
$$
\tilde\mu_{O} : \on{PB}_{n} \to \on{exp}(\hat\t_{n}). 
$$
  
Let us show that the various $\tilde\mu_{O}$ are are all conjugated to each other. 
Let $\on{can}_{O,O'}
\in {\bf PaB}(O,O')$ correspond to $e\in \on{B}_{n}$. Then $\on{can}_{O',O''}
\circ \on{can}_{O,O'}=\on{can}_{O,O''}$. Moreover, if we denote by 
$\sigma_{O}:\on{B}_{n}\to {\bf PaB}(O)$ the canonical identification, then 
$\sigma_{O'}(b) = \on{can}_{O,O'} \circ \sigma_{O}(b)\circ 
\on{can}_{O,O'}^{-1}$. Let us set $\Phi_{O,O'}:= F_{\Phi}(\on{can}_{O,O'})$. 
Then: 

1) $\Phi_{O,O'}\in \on{exp}(\hat\t_{n})$, 
$\Phi_{O',O''}\Phi_{O,O'}=\Phi_{O,O'}$;

2) $\tilde\mu_{O'}(b) = \Phi_{O,O'}\tilde\mu_{O}(b)\Phi_{O,O'}^{-1}$. 

If $O = \bullet(...(\bullet\bullet))$ is the `right parenthesization', 
the explicit formula for $\tilde\mu_{O}$ is 
$$
\tilde\mu_{O}(\sigma_{i}) = \Phi^{i,i+1,i+2...n} 
e^{t_{i,i+1}/2}s_{i} (\Phi^{i,i+1,i+2...n})^{-1}, \quad i =0,...,n-1. 
$$

The morphisms $\tilde\mu_{O}$ extend to isomorphisms between prounipotent 
completions as follows. The prounipotent completion of 
$\on{B}_{n}$ relative to $\on{B}_{n}\to S_{n}$ will be denoted 
$\on{B}_{n}(\kk,S_{n})$; it may be constructed as follows: $\on{B}_{n}$
acts by automorphisms of $\on{PB}_{n}$, hence of $\on{PB}_{n}(\kk)$; 
$\on{B}_{n}(\kk,S_{n})$ fits in an exact sequence $1\to \on{PB}_{n}(\kk)
\to \on{B}_{n}(\kk,S_n)\to S_{n}\to 1$ and identifies with the quotient of the 
semidirect product $\on{PB}_{n}(\kk)\rtimes \on{B}_{n}$
by the image of the morphism $\on{PB}_{n}\to 
\on{PB}_{n}(\kk)\rtimes \on{B}_{n}$, $g\mapsto (g^{-1},g)$ 
(which is a normal subgroup). Then the morphisms $\tilde\mu_{O}$ give rise to 
isomorphisms 
$$\begin{matrix}
\on{PB}_{n}(\kk) & \stackrel{\sim}{\to}& \on{exp}(\hat\t_{n})\\
 \downarrow & & \downarrow\\
 \on{B}_{n}(\kk,S_{n})  & \stackrel{\sim}{\to}& \on{exp}(\hat\t_{n})
 \rtimes S_{n}\end{matrix}$$ 
When $\Phi$ is the KZ associator (with coupling constant $2\pi\on{i}$), 
these isomorphisms are given by Sullivan's theory of minimal models 
applied to the configuration 
space of $n$ points in the complex plane (which computes all the rational homotopy 
groups of a simply-connected Kaehler manifold, but only the Malcev completion of 
its fundamental group in the non-simply-connected case, whence the name 
`$1$-formality').  

\subsection{Restriction to free groups}

Renumber $x_{ij}$, $i<j\in\{0,...,n\} $ and $t_{ij}$, $i\neq j\in 
\{0,...,n\}$ the generators for $\on{PB}_{n+1}$ and $\t_{n+1}$. 
Recall that $\on{PB}_{n+1}$ contains the free group with $n$ generators 
$\on{F}_{n} = \langle x_{01},...,x_{0,n}\rangle$ as a normal subgroup. 
Similarly, $\t_{n+1}$ contains the free Lie algebra with $n$ generators 
$\f_{n} = \on{Lie}(t_{01},...,t_{0,n})$. For coherence of notation with 
the previous sections, we will set $X_{i}:= x_{0i}$, $x_{i}:= t_{0i}$. 

\begin{proposition}
For any $O\in \on{Ob}({\bf PaB})$ with $|O|=n+1$, the morphism 
$\tilde\mu_{O}$ restricts to a morphism $\mu_{O} : \on{F}_{n}\to 
\on{exp}(\hat\f_{n})$, which extends to an isomorphism 
$\mu_{O}:\on{F}_{n}(\kk)\to \on{exp}(\hat\f_{n})$.
The composition of $\mu_{O}$ with the isomorphism $\on{exp}(\hat\f_{n})\to 
\on{F}_{n}(\kk)$, $\on{exp}(x_{i})\mapsto X_{i}$, is a tangential automorphism
of $\on{exp}(\hat\f_{n})$, i.e., an element of $\on{TAut}_{n}$. 
\end{proposition}

{\em Proof.} Let us first treat the case of $\mu_{n}:= 
\mu_{\bullet(...(\bullet\bullet))}$. As 
$x_{0i}=(\sigma_{i-2}...\sigma_{0})^{-1}\sigma_{i-1}^{2}(\sigma_{i-2}...
\sigma_{0})$, we have 
$\mu_{n}(x_{0i}) = \mu_{n}(\sigma_{i-2}...\sigma_{0})^{-1}\Phi^{i-1,i,i+1...n}
\cdot e^{t_{i-1,i}}\cdot 
(\Phi^{i-1,i,i+1...n})^{-1}\mu_{n}(\sigma_{i-2}...\sigma_{0})$. 
There exists $y_{i}\in \hat\t_{n+1}$ such that $\mu_{n}(\sigma_{i-2}...\sigma_{0})
=e^{y_{i}}s_{i-2}...s_{0}$, so for some $\tilde y_{i}\in \hat\t_{n+1}$, 
$$(\Phi^{i-1,i,i+1...n})^{-1}\mu_{n}(\sigma_{i-2}...\sigma_{0}) =
s_{i-2}...s_{0}e^{\tilde y_{i}}.
$$ 
Then $\mu_{n}(x_{0i}) = 
e^{-\tilde y_{i}}(s_{i-2}...s_{0})^{-1}e^{t_{i-1,i}} s_{i-2}...s_{0}
e^{\tilde y_{i}} = e^{-\tilde y_{i}}e^{t_{0i}}e^{\tilde y_{i}}$. 
As the action of $\hat\t_{n+1}$ on $\hat\f_{n}$ is by tangential automorphisms, 
we have $e^{-\tilde y_{i}}e^{t_{0i}}e^{\tilde y_{i}} = e^{z_{i}}e^{t_{0i}}
e^{-z_{i}}$ for some $z_{i}\in\hat\f_{n}$. So $\mu_{n}\circ (e^{t_{0i}}
\mapsto x_{0i})\in \on{TAut}_{n}$. The general case follows from the 
identity $\mu_{O'}(b) = \Phi_{O,O'}\mu_{O}(b)\Phi_{O,O'}^{-1}$ and the 
fact that for any $\Psi\in \on{exp}(\hat\t_{n+1})$, $x\mapsto \Psi x\Psi^{-1}$
induces a tangential automorphism of $\on{exp}(\hat\f_{n})$.  
\hfill \qed\medskip 

\begin{proposition} \label{prop:3:2}
If moreover $O = \bullet\otimes \bar O$, where $\bar O\in\on{Ob}({\bf PaB})$
has length $n$, then $\mu_{O}(X_{1}...X_{n}) = e^{x_{1}+...+x_{n}}$. 
\end{proposition}

{\em Proof.} The map $\on{PB}_{2}\to \on{PB}_{n+1}$, $p\mapsto 
p^{0,\widetilde{1...n}}$ takes $x_{01}$ to $x_{01}...x_{0n}=X_{1}...X_{n}$. 
Similarly to (\ref{diag2}), one proves that the diagram 
$$
\begin{matrix}
\on{PB}_{2} & \stackrel{p\mapsto 
p^{0,\widetilde{1...n}}}{\to} & \on{PB}_{n+1}\\
\scriptstyle{\tilde\mu_{\bullet\bullet}}
\downarrow & & \downarrow\scriptstyle{\tilde\mu_{\bullet\otimes \bar O}}\\
\on{exp}(\hat\t_{2}) & \stackrel{x\mapsto x^{0,1...n}}{\to}
& \on{exp}(\hat\t_{n+1}) 
\end{matrix}$$
commutes.  \hfill \qed \medskip 

The various isomorphisms $\mu_{O}$ are related by the identities
\begin{equation} \label{muO:muO'}
\mu_{O'} = \on{Ad}(\Phi_{O,O'}) \circ \mu_{O}; 
\end{equation}
the automorphisms $\on{Ad}(\Phi_{O,O'})$ are no longer necessarily inner.

\section{The identity $\mu_{O^{(i)}} = \mu_{O}^{1,2,...,ii+1,...,n}\circ 
\mu_{\bullet(\bullet\bullet)}^{i,i+1}$} \label{sec:mu:mu}

Let $O\in \on{Ob}({\bf PaB})$ be a parenthesized word of length $n$; 
its letters are numbered $0,...,n-1$. Let $i\in \{1,...,n-1\}$, let $O^{(i)}$
be the object obtained by replacing the letter $\bullet$ numbered $i$
by $(\bullet\bullet)$ (e.g., if $O=\bullet(\bullet\bullet)$, then 
$O^{(1)} = \bullet((\bullet\bullet)\bullet)$).  The purpose of this section 
is to show the identity 
$$
\mu_{O^{(i)}} = \mu_{O}^{1,2,...,ii+1,...,n}\circ 
\mu_{\bullet(\bullet\bullet)}^{i,i+1}. 
$$

\subsection{Free magmas and semigroups}
Recall that a magma is a triple $(M,M\times M\to M,e\in M)$
satisfying $e\times m\mapsto m$ and $m\times e\mapsto m$. 
A semigroup is a magma, where $M\times M\to M$ is associative.

Let $X$ be a finite set. Let $\Mg_X$ be the free magma generated by $X$
and $\Sg_X$ the semigroup generated by $X$.The assignments $X\mapsto \Sg_X$, 
$X\mapsto \Mg_X$ are functorial and we have a natural map $Mg_X \to \Sg_X$;
so we have a commutative diagram 
$$
\begin{matrix}
\Mg_X & \to  & \Sg_X \\
\downarrow & & \downarrow \\
\Mg_{\{\bullet\}}& \to & \Sg_{\{\bullet\}} = \NN
\end{matrix}
$$ 
This diagram is Cartesian, so $\Mg_X$ can be identified
with a fibered product. Explicitly, we have 
$\Sg_X = \sqcup_{n\geq 0} X^n$, $\Mg_{\{\bullet\}}
= \sqcup_{n\geq 0}\{$parenthesizations of the word $\bullet...\bullet$
of length $n\} = \sqcup_{n\geq 0} \{$rooted planar
binary trees with $n$ leaves$\}$, $\Mg_X = \sqcup_{n\geq 0}
\{$parenthesized words of length $n$ in the alphabet $X\}$.  

We denote by $w: \Mg_X\to \Sg_X$ (word),  
$P:\Mg_X\to \Mg_{\{\bullet\}}$ (parenthesization)
the natural maps; the various maps to $\NN$
are denoted by $x\mapsto |x|$ (length). 

Note that $S_n$ acts on $X^n$. For $w,w'\in \Sg_X$, 
with $|w| = |w'| = n$, we then set $S_{w,w'} = 
\{\sigma\in S_n | \sigma\cdot w = w'\}$. 

\subsection{A braided monoidal category ${\bf PaB}_X$}

We denote by BMC the `category' of braided monoidal categories (b.m.c.), 
where morphisms are the tensor functors. 
 
We define a functor $\on{Sets}\to \on{BMC}$, 
$X\mapsto {\bf PaB}_{X}$, adjoint to the 
`objects' functor $\on{BMC}\to \on{Sets}$, $\cC\mapsto \on{Ob}\cC$. 
This means that for any set $X$ and b.m.c. $\cC$, we have a natural 
bijection $\on{Mor}_{\on{Sets}}(X,\on{Ob}\cC) \simeq \on{Mor}_{\on{BMC}}
({\bf PaB}_{X},\cC)$. More precisely, we have an injection $X\subset 
\on{Ob}{\bf PaB}_{X}$, and for any b.m.c. $\cC$ and any map $X\to \on{Ob}\cC$, 
there is attached a tensor functor ${\bf PaB}_X\to \cC$, 
such that $\on{Ob}{\bf PaB}_{X}\to \on{Ob}\cC$ extends
$X\to \on{Ob}\cC$. When $X=\{\bullet\}$, ${\bf PaB}_X$
identifies with Bar-Natan's ${\bf PaB}$. 

We now construct ${\bf PaB}_X$. 
We set $\on{Ob}({\bf PaB}_X):= \Mg_X$. For $O,O'\in \Mg_X$, 
we set ${\bf PaB}_X(O,O') = \emptyset$ if $|O|\neq |O'|$, 
and $= \on{B}_n \times_\pi S_{w(O),w(O')}$ if $|O|=|O'|=n$
($\pi : \on{B}_n\to S_n$ is the canonical projection). So
${\bf PaB}_X(O,O') \subset \on{B}_n$; since $S_{w,w'}S_{w',w''}
\subset S_{w,w''}$, the product in $\on{B}_n$ restricts to a map 
${\bf PaB}_X(O,O') \times  {\bf PaB}_X(O',O'') \to 
{\bf PaB}_X(O,O'')$, which we define as the composition in 
${\bf PaB}_X$.

The tensor product is defined at the level of objects
by the product in $\Mg_X$, and at the level of morphisms 
is induced by the juxtaposition map $\on{B}_n\times 
\on{B}_{m}\to \on{B}_{n+m}$. 

We now construct the braiding and associativity constraints. 
For $O,O',O''\in Mg_X$, $a_{O,O',O''}\in {\bf PaB}_X(
(O\otimes O')\otimes O'',O\otimes (O'\otimes O''))$ 
is defined as the identity element in $\on{B}_{n+n'+n''}$ 
($n=|O|$, $n'=|O'|$, $n'' = |O''|$). 

Then $\beta_{O,O'} \in {\bf PaB}_X(O\otimes O',O'\otimes O)
\simeq \on{B}_{|O|+|O'|}$ corresponds to $\sigma_{n,n'}$
(one checks that the image $s_{n,n'}\in S_{n+n'}$ of $\sigma_{n,n'}$
belongs to the desired $S_{w,w'}$). 

One checks that ${\bf PaB}_X$, equipped with this structure, is a 
b.m.c., and that $X\mapsto{\bf PaB}_{X}$ is adjoint to the `objects' functor. 

\subsection{The category ${\bf PaCD}_X$}

We first define a tensor category $F_X$ as follows. 
$\on{Ob}(F_X) := \Sg_X$, and for $w,w'\in \Sg_X$, 
$F_X(w,w') = \emptyset$ if $|w|\neq |w'|$, 
and $= (\on{exp}(\hat\t_n)\rtimes S_n)\times_\pi S_{w,w'}$
else, where $\pi : \on{exp}(\hat\t_n)\rtimes S_n\to S_n$ is the 
canonical projection. The composition is defined as above, 
using the product in $\on{exp}(\hat\t_n)\rtimes S_n$, again 
using $S_{w,w'}S_{w',w''} \subset S_{w,w''}$. 

The tensor product is defined, at the level of objects,  
by the semigroup law, and at the level of morphisms using 
the juxtaposition $(\on{exp}(\hat\t_n)\rtimes S_n)
\times (\on{exp}(\hat\t_{n'})\rtimes S_{n'})
\to \on{exp}(\hat\t_{n+n'})\rtimes S_{n+n'}$. 

Let $\Phi\in M_{1}(\kk)$. For $X = \{\bullet\}$, 
$\Sg_{X} = \NN$ (we then have $n\otimes m = n+m$). 
For $n,n',n''\in\NN$, we then set 
$$
a_{n,n',n''}:= 
\Phi^{1...n,n+1...n+n',n+n'+1...n+n'+n''}
\in \on{exp}(\hat\t_{n+n'+n''})
\in F_{\{\bullet\}}(n\otimes n'\otimes n''); 
$$
$s_{n,n'}\in S_{n+n'}$ is the block permutation 
$i\mapsto n'+i$ ($i\in [n]$), $n+i\mapsto i$ ($i\in [n']$)
and 
$$
\beta_{n,n'}:= 
(e^{t_{12}/2})^{1...n,n+1...n+n'}s_{n,n'}
\in \on{exp}(\hat\t_{n+n'}) \rtimes S_{n+n'} = 
F_{\{\bullet\}}(n\otimes n',n'\otimes n). 
$$
We note that if $X$ is arbitrary and $w,w',w''\in \Sg_{X}$, then 
$a_{|w|,|w'|,|w''|}\in F_{X}(w\otimes w'\otimes w'')$ and 
$\beta_{|w|,|w'|}\in F_{X}(w\otimes w',w'\otimes w)$. 

We define the category ${\bf PaCD}_X$ by 
$\on{Ob}({\bf PaCD}_X):= \Mg_X$, and for $O,O'\in \Mg_X$, 
we set ${\bf PaCD}_X(O,O') := F_X(w(O),w(O'))$.
The tensor product is defined at the level of objects 
as the product in $\Mg_2$; as $w : \Mg_2\to \Sg_2$ is 
compatible with products, a tensor product is defined at the level 
of morphisms by 
${\bf PaCD}_X(O_1,O_2) \otimes {\bf PaCD}_X(O'_1,O'_2) 
= F_X(w(O_1),w(O_2)) \otimes  F_X(w(O'_1),w(O'_2))
\to F_X(w(O_1)\otimes w(O'_1),w(O_2) \otimes w(O'_2)) = 
F_X(w(O_1\otimes O'_1),w(O_2\otimes O'_2)) = 
{\bf PaCD}_X(O_1\otimes O'_1,O_2\otimes O'_2)$. 

Let $\Phi\in M_1(\kk)$. Then $\Phi$ gives rise to a 
b.m.c. structure on ${\bf PaCD}_X$ by 
$a_{O,O',O''}:= a_{|O|,|O'|,|O''|} \in 
F_{X}(w(O)\otimes w(O')\otimes w(O''))
= {\bf PaCD}_X((O\otimes O') \otimes O'',
O\otimes (O'\otimes O''))$ and $\beta_{O,O'}:= \beta_{|O|,|O'|}
\in F_{X}(w(O)\otimes w(O'),w(O')\otimes w(O)) = 
{\bf PaCD}_X(O\otimes O',O'\otimes O)$ for $O,O',O''\in \Mg_{X}$. 

We denote by ${\bf PaCD}_X^\Phi$ the resulting b.m.c.

\subsection{Tensor functors}

When $X = X_1 := \{\bullet\}$,  ${\bf PaB}_X$ coincides with 
${\bf PaB}$; we also denote $\Mg_{X}$, ${\bf PaCD}_X^\Phi$ 
by $\Mg$, ${\bf PaCD}^\Phi$. For $X = X_2 := \{\bullet,\circ\}$, we denote 
${\bf PaB}_X$, ${\bf PaCD}_X^\Phi$, 
$\Mg_{X}$, $\Sg_{X}$ by ${\bf PaB}_2$, ${\bf PaCD}_2^\Phi$, $\Mg_{2}$, 
$\Sg_{2}$. 

We define ${\bf PaB}_2\to {\bf PaB}$ as the tensor functor  
induced by the map $X_2\to \Mg_1$, 
$\bullet\mapsto \bullet$, $\circ\mapsto \bullet\bullet$. 

We denote by ${\bf PaB}\to {\bf PaCD}^\Phi$ the tensor 
functor induced by the canonical injection 
$X_1\to \on{Ob}({\bf PaCD}^\Phi) = \Mg_1$. 

Similarly, we denote by ${\bf PaB}_2\to {\bf PaCD}_2^\Phi$ the tensor 
functor induced by the canonical injection 
$X_2\to \on{Ob}({\bf PaCD}^\Phi_2) = \Mg_2$. 

Let us now construct a functor $F_{X_2} \to F_{X_1}$. 
At the level of objects, this is the semigroup morphism
$\Sg_2\to \Sg_1$ induced by the map 
$l:X_2\to \Sg_1 \simeq \NN$, $w\mapsto \tilde w$ given by 
$\bullet\mapsto 1$ and $\circ \mapsto 2$. 
So for $w=(w_1,...,w_n)\in \sqcup_{n\geq 0} X_2^n$, 
$\tilde w = \sum_{i=1}^n l(w_i)$, where $l(\bullet)=1$
and $l(\bullet\bullet)=2$. 
Let us now define the functor at the level of morphisms, i.e. 
the maps $F_{X_2}(w,w') \to F_{X_1}(\tilde w,\tilde w')$. 
As $F_{X_2}(w,w') = \emptyset$ unless 
$(\on{card}\{i|w_i = \bullet\}, \on{card}\{i|w_i = \circ\}) 
= (\on{card}\{i|w'_i = \bullet\}, 
\on{card}\{i|w'_i = \circ\})$, we will assume that these
pairs of integers are equal (in particular $|w| = |w'|$); 
we denote this pair by $(n_1,n_2)$. Note that 
$|w|=|w'|=n_1+n_2$, while $\tilde w=\tilde w'=n_1+2n_2$. 

There is a unique non-decreasing map 
$\phi_w : [n_1+2n_2]\to [n_1+n_2]$, such that 
$i$ has one preimage by $\phi_w$ if $w_i=\bullet$
and two preimages if $w_i = \circ$; 
for example, if $w=(\bullet,\bullet,\circ,\circ,\bullet)$, 
then $\phi_w : [7]\to [5]$ is 
$(1,...,7)\mapsto (1,2,3,3,4,4,5)$. 

Moreover, for any $\sigma\in S_{n_1+n_2}$, 
there is a unique $\sigma^w\in S_{n_1+2n_2}$
such that: (a) $\sigma\circ \phi_w = \phi_{w'} \circ
\sigma^w$, where $w'= \sigma \cdot w$, so that $\sigma^w$ restricts to 
bijections $\phi_w^{-1}(i)\to \phi_{w'}^{-1}(i)$; 
(b) these bijections are increasing (this 
condition is nonempty only if $\on{card}\phi_w^{-1}(i)>1$). 
The map $\sigma \mapsto \sigma'$ is a group morphism 
$S_{n_1+n_2}\to S_{n_1+2n_2}$ (it maps a permutation to 
a block permutation); for example, if $w=(\circ,\bullet,\bullet)$, 
this map is $S_3\to S_4$, $
 \bigl(\begin{smallmatrix}
 1& 2& 3\\
 2& 1& 3\end{smallmatrix}\bigl)
\mapsto \bigl(\begin{smallmatrix}
 1& 2& 3& 4\\
 3& 1& 2& 4\end{smallmatrix}\bigl)$, 
$\bigl(\begin{smallmatrix}
 1& 2& 3\\
 1& 3& 2\end{smallmatrix}\bigl)
\mapsto 
 \bigl(\begin{smallmatrix}
 1& 2& 3& 4\\
 1& 2& 4& 3\end{smallmatrix}\bigl)$. 

The morphisms $\t_{n_1+n_2}\to \t_{n_1+2n_2}$, $x\mapsto 
x^{\phi_w}$ and $S_{n_1+n_2}\to S_{n_1+2n_2}$, 
$\sigma\mapsto \sigma^w$ are compatible, so we obtain a 
group morphism $\on{exp}(\hat\t_{n_1+n_2}) \rtimes S_{n_1+n_2}
\to \on{exp}(\hat\t_{n_1+2n_2}) \rtimes S_{n_1+2n_2}$. 
We then define $F_{X_2}(w,w')\to F_{X_1}(\tilde w,\tilde w')$
as the restriction of this group morphism. 
One checks that this map is compatible with tensor 
products, so we have defined a tensor functor $F_{X_2}\to F_{X_1}$. 

The tensor functor $F_{X_{2}}\to F_{X_{1}}$ extends to a tensor functor 
${\bf PaCD}_2^\Phi \to {\bf PaCD}^\Phi$ as follows. 
There is a unique magma morphism 
$\Mg_2\to \Mg_1$, $O\mapsto \tilde O$, extending the map 
$X_2\to \Mg_1$, $\bullet\mapsto\bullet$, $\circ\mapsto 
\bullet\bullet$. It is such that the diagram 
$$
\begin{matrix}
\Mg_2 &\to & \Mg_1 \\ 
\downarrow & & \downarrow \\
 \Sg_2 & \to & \Sg_1
\end{matrix}
$$
commutes. The functor ${\bf PaCD}_2^\Phi \to {\bf PaCD}^\Phi$
is defined, at the level of objects, as the map
$\Mg_2\to \Mg_1$ and at the level of morphisms by 
${\bf PaCD}_2^\Phi(O,O') = F_{X_2}(w(O),w(O'))
\to F_{X_1}(\widetilde{w(O)},\widetilde{w(O')})
= F_{X_1}(w(\tilde O),w(\tilde O')) = {\bf PaCD}(O,O')$. 

It remains to show that it takes braidings and
associativity constraints to their analogues. 
Namely:

(a) it takes $\beta_{O,O'}\in 
{\bf PaCD}_{2}(O\otimes O',O'\otimes O)$
to $\beta_{\tilde O,\tilde O'}\in {\bf PaCD}
(\tilde O\otimes \tilde O',\tilde O'\otimes \tilde O)$. 

(b) it takes $a_{O,O',O''}\in {\bf PaCD}_{2}
((O\otimes O')\otimes O'',O\otimes (O'\otimes O''))$
to $a_{\tilde O,\tilde O',\tilde O''}\in {\bf PaCD}
((\tilde O\otimes \tilde O')\otimes \tilde O'',
\tilde O\otimes (\tilde O'\otimes \tilde O''))$. 

To prove (a), let $w,w' = w(O),w(O')$, 
$(\on{card}\{i|w_{i}=\bullet\}, 
\on{card}\{i|w_{i}=\circ\})=(n_{1},n_{2})$, 
$(\on{card}\{i|w'_{i}=\bullet\}, 
\on{card}\{i|w'_{i}=\circ\})=(n'_{1},n'_{2})$. 
Then $\beta_{O,O'} = \beta_{n_{1}+n_{2},n'_{1}+n'_{2}}
\in F_{X_{2}}(w\otimes w',w'\otimes w)$. 
Similarly, $\beta_{\tilde O,\tilde O'} = 
\beta_{n_{1}+2n_{2},n'_{1}+2n'_{2}}
\in F_{X_{1}}(\tilde w\otimes \tilde w',
\tilde w'\otimes \tilde w)$. 

Now note that: 
$$((t_{12})^{1...n,n+1...n+n'})^{\phi_{w\otimes w'}}
= (t_{12})^{1...n_{1}+2n_{2},n_{1}+2n_{2}+1...n_{1}+2n_{2}
+n'_{1}+2n'_{2}},$$
and 
$$ (s_{n,n'})^{w\otimes w'} = s_{n_{1}+2n_{2},
n'_{1}+2n'_{2}}.
$$  
So the map $F_{X_{2}}(w\otimes w',w'\otimes w)
\to F_{X_{1}}(\tilde w\otimes \tilde w',
\tilde w'\otimes \tilde w)$
takes $\beta_{n,n'}$ to $\beta_{n_{1}+2n_{2},n'_{1}+2n'_{2}}$. 
The proof of (b) is similar. 

Then the diagram of functors 
$$
\begin{matrix}
{\bf PaB}_2 & \to & {\bf PaB}\\
\downarrow & & \downarrow \\
{\bf PaCD}_2^\Phi & \to & {\bf PaCD}^\Phi
\end{matrix}
$$
commutes by universal properties 
(the two composed functors ${\bf PaB}_2\to {\bf PaCD}^{\Phi}$
coincides as their restrictions to the elements of $X_{2}\subset 
\on{Ob}({\bf PaB}_{2})$ do). 

\begin{remark}
More generally, to any map $X\to \Mg_1$, one associates a 
tensor functor ${\bf PaCD}_X^\Phi\to {\bf PaCD}^\Phi$, defined 
at the level of objects by the extension of this map to a morphism 
$\Mg_X\to \Mg_1$ and at the level of morphisms by suitable iterations 
of cobrackets, and it is such that 
$$
\begin{matrix}
{\bf PaB}_X & \to  & {\bf PaB}\\
\downarrow & & \downarrow \\
{\bf PaCD}_X & \to & {\bf PaCD} 
\end{matrix}
$$
commutes. 
\end{remark}

\subsection{Relation between braid groups representations} \label{rel:braid}

Let $n\geq 1$, let $i\in [n]$, let $w^{i} = (\bullet,...,\bullet,\circ,
\bullet,...,\bullet)\in \Sg_{2}$ be given by $w_{i}=\circ$
and $w_{j} = \bullet$ for $j\in [n]-\{i\}$. Let $O\in \Mg_{2}$
be such that $w(O)=w^{i}$. We have proved that the diagram 
$$
\begin{matrix}
{\bf PaB}_{2}(O) & \to & {\bf PaB}(\tilde O) \\
\downarrow & & \downarrow \\
{\bf PaCD}_{2}(O) & \to & {\bf PaCD}(\tilde O)
\end{matrix}
$$
commutes. 

We have isomorphisms: 

 ${\bf PaB}_{2}(O)\simeq 
\on{B}_{n}\times_{\pi}S_{n-1}$, where $S_{n-1}\subset S_{n}$
identifies with $\{\sigma\in S_{n}|\sigma(i)=i\}$; 

${\bf PaCD}_{2}(O)\simeq (\on{exp}(\hat\t_{n})\rtimes S_{n})
\times_{\pi} S_{n-1}$; 

${\bf PaB}(\tilde O)\simeq \on{B}_{n+1}$; 

${\bf PaCD}(\tilde O) \simeq \on{exp}(\hat\t_{n+1})\rtimes S_{n+1}$. 

For $O\in \Mg_{X_{1}}$, $|O|=n$, the morphism 
${\bf PaB}(O)\to {\bf PaCD}(O)$ is a morphism 
$\mu_{O} : \on{B}_{n}\to \on{exp}(\hat\t_{n})\rtimes S_{n}$. 
Note that if $O_{X}\in \Mg_{X}$ and $O := P(O_{X})$, then we have a 
commutative diagram 
$$ 
\begin{matrix}
{\bf PaB}_{X}(O_{X}) & \to & {\bf PaCD}_{X}(O_{X})\\
\downarrow   & & \downarrow \\
 \on{B}_{n}  & \stackrel{\mu_{O}}{\to}& 
\on{exp}(\hat\t_{n})\rtimes S_{n} \end{matrix}
$$
where the vertical maps are injective. 

The above commutative diagram therefore inserts in a diagram 
\begin{equation} \label{diag1}
\begin{matrix}
\on{B}_{n}& \leftarrow& \on{B}_{n}\times_{\pi}S_{n-1} & \stackrel{
1,2,...,\widetilde{ii+1},...,n}{\to} 
& \on{B}_{n+1}\\
\scriptstyle{\mu_{O}}\downarrow 
& & \downarrow & & \downarrow\scriptstyle{\mu_{O^{(i)}}} \\
\on{exp}(\hat\t_{n})\rtimes S_{n}
& \leftarrow & (\on{exp}(\hat\t_{n})\rtimes S_{n})\times_{\pi}S_{n-1} & 
\stackrel{1,2,...,ii+1,...,n}{\to} & \on{exp}(\hat\t_{n+1})\rtimes S_{n+1} 
\end{matrix}
\end{equation}
Restricting to pure braid groups, we obtain the commutative diagram 
\begin{equation} \label{diag2}
\begin{matrix}
\on{PB}_{n} & \stackrel{1,2,...,\widetilde{ii+1},...,n}{\to} & \on{PB}_{n+1}\\
\scriptstyle{\mu_{O}}\downarrow
 & & \downarrow\scriptstyle{\mu_{O^{(i)}}}\\
 \on{exp}(\hat\t_{n}) & \stackrel{1,2,...,ii+1,i+2,...,n+1}{\to}& 
 \on{exp}(\hat\t_{n+1})
\end{matrix}
\end{equation}

\subsection{Relation between $\mu_O$ and $\mu_{O^{(i)}}$}

Let $O\in \on{Ob}({\bf PaB})$, $|O|=n$. 
We index letters in $O$ by $0,...,n-1$, fix an index $i\neq 0$ and construct 
$O^{(i)}$ by doubling inside $O$ the $\bullet$ with index $i$. 

$O$ gives rise to a morphism $\tilde\mu_O : \on{B}_n\to
\on{exp}(\hat\t_n) \rtimes S_n$, which restricts to 
$\mu_O : \on{F}_{n-1}\to \on{exp}(\hat\f_{n-1})$. 
Similarly, $\tilde\mu_{O^{(i)}} : \on{B}_{n+1}\to
\on{exp}(\hat\t_{n+1}) \rtimes S_{n+1}$ restricts to 
$\mu_{O^{(i)}} : \on{F}_{n}\to \on{exp}(\hat\f_{n})$.

We want to prove that 
\begin{equation} \label{main:in:section}
\mu_{O^{(i)}} = \mu_O^{1,2,...,ii+1,...,n} \circ 
\mu_{\bullet(\bullet\bullet)}^{i,i+1}. 
\end{equation}

We first show that there are uniquely determined elements 
$g_1,...,g_{n-1}\in \on{exp}(\hat\f_{n-1})$ and 
$g,h\in \on{exp}(\hat\f_2)$ such that: 

(a) $\mu_O = \lbr g_1(x_1,...,x_{n-1}),...,
g_{n-1}(x_1,...,x_{n-1})\rbr$, $\on{log} g_i(x_1,...,x_{n-1}) 
= -{\frac{1}{2}}(x_1+...+x_{i-1}) + O(x^2)$, and\footnote{
$O(x^2)$ means an element of $\hat\f_{n-1}$ of valuation  
$\geq 2$.}

(b) $\mu_{\bullet(\bullet\bullet)} = \lbr g(x_1,x_2),h(x_1,x_2) \rbr$, 
$\on{log}g(x_1,x_2) = O(x^2)$, $\on{log}h(x_1,x_2) = 
-{\frac{1}{2}}x_1 + O(x^2)$.

Let us prove the first statement (it actually contains the second 
statement as a particular case). The elements 
$g_i(x_1,...,x_{n-1})$ are uniquely
determined by the equality $\mu_O = \lbr g_1,...,g_{n-1}\rbr$, 
together with the condition that the coefficient of $x_i$ in the expansion 
of $\on{log}g_i$ vanishes. We should then prove that 
$\on{log} g_i = -{\frac{1}{2}}(x_1+...+x_{i-1}) + O(x^2)$. 
We have 
$$\tilde\mu_O(\sigma_j)=e^{a_j} \cdot e^{t_{j-1,j}/2}s_j 
\cdot e^{-a_j},
$$ 
where $a_j\in \hat\t_n$ has valuation $\geq 2$
(we write this as $a_j\in O(t^2)$), and 
$$
\mu_O(X_i) = \tilde \mu_O(\sigma_1)^{-1}...\tilde 
\mu_O(\sigma_{i-1})^{-1}\tilde\mu_O(\sigma_i)^2
\tilde\mu_O(\sigma_{i-1})...\tilde\mu_O(\sigma_1).
$$ 
Now 
$$
\tilde\mu_O(\sigma_{i-1})...\tilde\mu_O(\sigma_1)
= s_{i-1}...s_1 e^{{1\over 2}(x_1+...+x_{i-1})+O(t^2)}
$$
and $\tilde\mu_O(\sigma_i^2) = e^{a_i}e^{t_{i-1,i}}e^{-a_i}$.
It follows that 
$$\mu_O(X_i) = e^{-{1\over 2}(x_1+...+x_{i-1})+O(t^2)}
e^{\tilde a_i}\cdot e^{x_i} \cdot (same)^{-1},
$$ 
where $\tilde a_i = 
s_1...s_{i-1}\cdot a_i\cdot s_{i-1}...s_1 \in O(t^2)$, so  
$\mu_O(X_i) = e^{-{1\over 2}(x_1+...+x_{i-1})+O(t^2)}
\cdot e^{x_i} \cdot (same)^{-1}$, which implies that $g_i$
has the announced form.

To prove (\ref{main:in:section}), we need to prove the equality   
\begin{align} \label{second:in:section}
& \mu_{O^{(i)}} = 
\lbr g_1(x_1,...,x_i+x_{i+1},...,x_{n}),...,
g_i(x_1,...,x_i+x_{i+1},...,x_{n})g(x_i,x_{i+1}),
\\ & \nonumber 
g_i(x_1,...,x_i+x_{i+1},...,x_{n})h(x_i,x_{i+1}),...,
g_{n-1}(x_1,...,x_i+x_{i+1},...,x_{n})\rbr. 
\end{align}

(\ref{diag2}) implies that the diagram
$$
\begin{matrix}
\on{F}_{n-1} &\to  & \on{F}_n \\
\scriptstyle{\mu_O}\downarrow 
& & \downarrow \scriptstyle{\mu_{O^{(i)}}}\\
\on{exp}(\hat\f_{n-1}) & \to & \on{exp}(\hat\f_{n})
\end{matrix}
$$
commutes, where the upper morphism takes $X_j$ ($j\in [n-1]$)
to: $X_j$ if $j<i$, $X_iX_{i+1}$ if $j=i$, $X_{j+1}$ if $j>i+1$
and the lower morphism is similarly defined (replacing products 
by sums and $X_k$'s by $x_k$'s).  
Specializing to the generators $X_j$ ($j\neq i$) of $\on{F}_{n-1}$, 
this yields 
$$
\mu_{O^{(i)}}(X_j) = g_j^{0,1,...,ii+1,...,n}\cdot e^{x_j}\cdot (same)^{-1}
$$
for $j<i$ and 
$$
\mu_{O^{(i)}}(X_j) = g_{j-1}^{0,1,...,ii+1,...,n}\cdot e^{x_j}\cdot
 (same)^{-1}
$$
for $j>i+1$, which implies that (\ref{second:in:section})
holds when applied to the generators $X_j$, $j\neq i,i+1$. 

We now prove that (\ref{second:in:section})
also holds when applied to $X_i$ and $X_{i+1}$. 

The morphism $X_i\in \on{B}_n = {\bf PaB}(O,O)$
can be decomposed as 
$$
O \stackrel{(\sigma_{i-2}...\sigma_0)^{-1}}{\to} 
(O_1\otimes (\bullet\bullet))\otimes O_2
\stackrel{\sigma_{i-1}^2}{\to} (O_1\otimes (\bullet\bullet))\otimes O_2 
\stackrel{\sigma_{i-2}...\sigma_0}{\to} O.  
$$
Here the braid group elements indicate the morphisms. 
Let $\gamma\in \on{exp}(\hat\t_n)\rtimes S_n$
be the image of the morphism $O \stackrel{(\sigma_{i-2}...\sigma_0)^{-1}}{\to} 
(O_1\otimes (\bullet\bullet))\otimes O_2$ under ${\bf PaB}\to {\bf PaCD}$; 
its image in $S_n$ is the permutation $s_0...s_{i-2}$, i.e., 
$(0,...,n-1)\mapsto (i-1,0,1,...,i-2,i,i+1,...,n-1)$. 
The image of $(O_1\otimes (\bullet\bullet))\otimes O_2
\stackrel{\sigma_{i-1}^2}{\to} (O_1\otimes (\bullet\bullet))\otimes O_2$
is $e^{t_{i-1,i}}$, therefore the image of $X_i$ is 
$$
\mu_O(X_i) = \gamma e^{t_{i-1,i}} \gamma^{-1}.  
$$
We have $\gamma = \gamma_0 s_0...s_{i-2}$, where $\gamma_0\in 
\on{exp}(\hat\t_n)$. As $s_0...s_{i-2}\cdot t_{i-1,i} = x_i$, 
we have 
$$
\mu_O(X_i) = \gamma_0 e^{x_i} \gamma_0^{-1}. 
$$
As this image is also $g_i(x_1,...,x_{n-1}) \cdot e^{x_i} \cdot (same)^{-1}$, 
we derive from this that $g_i^{-1}\gamma_0$ commutes with $x_i$, 
hence by Proposition \ref{prop:comm}
has the form $e^{\lambda x_i}\alpha^{0i,1,2,...,i-1,i+1,...,n-1}$, 
where $\alpha\in \on{exp}(\hat\t_{n-1})$. 

Since $\mu_O(\sigma_j)=s_j e^{t_{j,j+1}/2}$, we get 
$\on{log}\gamma_0 = -{1\over 2}(x_1+...+x_{i-1}) + O(x^2)$. 
Comparing linear terms in $x_i$, we get $\lambda=0$. 

Let us now compute $\mu_{O^{(i)}}(X_i)$. 
The morphism $X_i \in \on{B}_{n+1} = {\bf PaB}(O^{(i)},O^{(i)})$
can be decomposed as 
$$
O^{(i)} \stackrel{(\sigma_{i-2}...\sigma_0)^{-1}}{\to}
(O_1 \otimes (\bullet(\bullet\bullet))) \otimes O_2
\stackrel{\sigma_{i-1}^2}{\to}
(O_1 \otimes (\bullet(\bullet\bullet))) \otimes O_2
\stackrel{\sigma_{i-2}...\sigma_0}{\to}
O^{(i)} 
$$
(here $\sigma_{i-1}^2$ involves the two first $\bullet$ of 
$\bullet(\bullet\bullet)$). 
The morphism $O^{(i)} \stackrel{(\sigma_{i-2}...\sigma_0)^{-1}}{\to}
(O_1 \otimes (\bullet(\bullet\bullet))) \otimes O_2$
is obtained from $O^{(i)} \stackrel{(\sigma_{i-2}...\sigma_0)^{-1}}{\to}
(O_1 \otimes (\bullet\bullet)) \otimes O_2$
by the operation of doubling of the $i$th strand, 
so its image is $\gamma^{0,1,2,...,ii+1,...,n}
= \gamma_0^{0,1,2,...,ii+1,...,n}(s_0...s_{i-2})$. 
The image of $\bullet(\bullet\bullet) \stackrel{\sigma_1^2}{\to}
\bullet(\bullet\bullet)$ is $g(x_1,x_2)\cdot e^{x_1}\cdot (same)^{-1}$, 
so the image of 
$$
(O_1 \otimes (\bullet(\bullet\bullet))) \otimes O_2
\stackrel{\sigma_{i-1}^2}{\to}
(O_1 \otimes (\bullet(\bullet\bullet))) \otimes O_2
$$
is $g(t_{i-1,i},t_{i-1,i+1})e^{t_{i-1,i}}(same)^{-1}$. 
It follows that 
$$
\mu_{O^{(i)}}(X_i) = \gamma^{0,1,2,...,ii+1,...,n}
g(t_{i-1,i},t_{i-1,i+1})\cdot e^{t_{i-1,i}} \cdot (same)^{-1}
= \gamma_0^{0,1,2,...,ii+1,...,n}
g(x_i,x_{i+1})\cdot e^{x_i} \cdot (same)^{-1}. 
$$
Now we claim that 
$$
\gamma_0^{0,1,2,...,ii+1,...,n}
g(x_i,x_{i+1})e^{x_i} (same)^{-1} = 
g_i^{0,1,2,...,ii+1,...,n}
g(x_i,x_{i+1})\cdot e^{x_i}\cdot  (same)^{-1}.
$$
Indeed, 
\begin{align*}
& (g_i^{-1}\gamma_0)^{0,1,2,...,ii+1,...,n}g(x_i,x_{i+1})
\cdot e^{x_i}\cdot (same)^{-1} 
\\ & 
= 
(\alpha^{0i,1,2,...,i-1,i+1,...,n-1}
)^{0,1,2,...,ii+1,...,n}g(x_i,x_{i+1})
\cdot e^{x_i}\cdot (same)^{-1} 
\\ & = 
\alpha^{0ii+1,2,3,...,i-1,i+2,...,} 
g(x_i,x_{i+1})\cdot e^{x_i}\cdot (same)^{-1}. 
\end{align*}
Now $x_i$ and $x_{i+1}$ commute with any $\alpha^{0ii+1,...}$, 
so this is $g(x_i,x_{i+1})\cdot e^{x_i}\cdot (same)^{-1}$. 

So we get  
$$
\mu_{O^{(i)}}(X_i) = g_i^{0,1,2,...,ii+1,...,n}
g(x_i,x_{i+1})\cdot e^{x_i}\cdot (same)^{-1}.
$$
The same argument shows that 
$$
\mu_{O^{(i)}}(X_{i+1}) = g_i^{0,1,2,...,ii+1,...,n}
h(x_i,x_{i+1})\cdot e^{x_{i+1}}\cdot  (same)^{-1}, 
$$
as wanted.

\section{The map $M_{1}(\kk)\to \on{SolKV}(\kk)$} \label{pf:main}

We show that for $\Phi\in M_{1}(\kk)$, $\mu_{\Phi}\in \on{SolKV}(\kk)$. 
By construction of $\mu_\Phi$, we have $\mu_{\Phi}(X)\sim e^{x}$, 
$\mu_{\Phi}(Y)\sim e^{y}$, so $\mu_{\Phi}\in \on{TAut}_{2}$.

\subsection{Proof of $\on{Ad}\Phi(t_{12},t_{23})\circ
\mu_{\Phi}^{12,3}\circ\mu_{\Phi}^{1,2} = 
\mu_{\Phi}^{1,23}\circ\mu_{\Phi}^{2,3}$} \label{pf:Phimumu}
 
We first prove: 

\begin{proposition} \label{prop:mu}
1)  $\mu_{\bullet(\bullet\bullet)} = \mu_{\Phi}$. 

2) $\Phi_{\bullet((\bullet\bullet)\bullet),\bullet(\bullet(\bullet\bullet))}
= \Phi(t_{12},t_{23})$. 
\end{proposition}

{\em Proof.} Let us prove 1). $x_{01}\in \on{B}_{3} 
= {\bf PaB}(\bullet(\bullet\bullet))$ corresponds to 
$a_{\bullet,\bullet,\bullet} \circ (\beta_{\bullet,\bullet}^{2}
\otimes\on{id}_{\bullet}) \circ a_{\bullet,\bullet,\bullet}^{-1}$. 
The image of this element in $\on{exp}(\hat\t_{3})\rtimes S_{3}$
is $\mu_{\bullet(\bullet\bullet)}(X) = 
\Phi(t_{01},t_{12}) e^{t_{01}}\Phi(t_{01},t_{12})^{-1}$. Since 
$t_{01}+t_{12}+t_{02}$ is central in $\t_{3}$ and since $\Phi$
is group-like, this is 
$\Phi(t_{01},-t_{01}-t_{02}) e^{t_{01}}\Phi(t_{01},-t_{01}-t_{02})^{-1}
= \Phi(x,-x-y)e^{x}\Phi(x,-x-y)^{-1} = \mu_{\Phi}(X)$.
Similarly, $x_{02}$ corresponds to 
$(\on{id}_{\bullet}\otimes \beta_{\bullet,\bullet})
\circ a_{\bullet,\bullet,\bullet} \circ (
\beta_{\bullet,\bullet}^{2}\otimes \on{id}_{\bullet}) 
\circ a_{\bullet,\bullet,\bullet}^{-1}
\circ (\on{id}_{\bullet}\otimes \beta_{\bullet,\bullet}^{-1})$. 
The image of this
element in $\on{exp}(\hat\t_{3})\rtimes S_{3}$ is 
 \begin{equation*} \begin{split}
 & \mu_{\bullet(\bullet\bullet)}(Y) \\
  & = e^{t_{12}/2}(12)\Phi(t_{01},t_{12}) e^{t_{01}}\Phi(t_{01},t_{12})^{-1}
 (12)e^{-t_{12}/2} = 
 e^{t_{12}/2}\Phi(t_{02},t_{12}) e^{t_{02}}\Phi(t_{02},t_{12})^{-1}
 e^{-t_{12}/2} \\
  & = e^{-(t_{01}+t_{02})/2}\Phi(t_{02},-t_{01}-t_{02}) 
  e^{t_{02}}\Phi(t_{02},-t_{01}-t_{02})^{-1}
 e^{(t_{01}+t_{02})/2}  \\
  & = 
 e^{-(x+y)/2}\Phi(y,-x-y)e^{y}\Phi(y,-x-y)^{-1}e^{(x+y)/2} = \mu_{\Phi}(Y). 
\end{split}\end{equation*}
So $\mu_{\bullet(\bullet\bullet)}=\mu_{\Phi}$. 

Let us now prove 2). Let $O:= \bullet((\bullet\bullet)\bullet)$, 
$O':= \bullet(\bullet(\bullet\bullet))$. Then $\on{can}_{O,O'} = 
\on{id}_{\bullet}\otimes a_{\bullet,\bullet,\bullet}\in
{\bf PaB}(O,O')$, whose image in ${\bf PaCD}(O,O') = 
\on{exp}(\hat\t_{4})\rtimes S_{4}$ is $\Phi(t_{12},t_{23})
= \Phi_{O,O'}$. 
\hfill \qed\medskip 

We now prove (\ref{Phimumu}). 
Applying (\ref{muO:muOnew}) to $O=\bullet(\bullet\bullet)$ and $i=1,2$, 
and using $\mu_{\bullet(\bullet\bullet)} = \mu_{\Phi}$, we get 
$$
\mu_{\bullet((\bullet\bullet)\bullet)} = 
\mu_{\Phi}^{12,3}\circ\mu_{\Phi}^{1,2},
\quad 
\mu_{\bullet(\bullet(\bullet\bullet))} = 
\mu_{\Phi}^{1,23}\circ\mu_{\Phi}^{2,3}. 
$$
Moreover, (\ref{muO:muO'}) implies 
$$
\on{Ad}\Phi_{\bullet((\bullet\bullet)\bullet),
\bullet(\bullet(\bullet\bullet))}\circ 
\mu_{\bullet((\bullet\bullet)\bullet)}= 
\mu_{\bullet(\bullet(\bullet\bullet))} . $$ 
As $\Phi_{\bullet((\bullet\bullet)\bullet),
\bullet(\bullet(\bullet\bullet))}
= \Phi(t_{12},t_{23})$, we get  (\ref{Phimumu}). 

\subsection{Proof of $\mu_{\Phi}(XY)=e^{x+y}$} \label{X:Y}

We will give three proofs: 

{\it First proof.}  We have 
\begin{equation*} \begin{split}
&  \mu_{\Phi}(XY)=\mu_{\Phi}(X) \mu_{\Phi}(Y)  \\ 
& =\Phi(x,-x-y)e^{x}\Phi(-x-y,x) e^{-(x+y)/2}\Phi(y,-x-y)e^{y}
\Phi(-x-y,x)e^{(x+y)/2}  \\ 
& =\Phi(x,-x-y)e^{x/2}\Phi(y,x) e^{y/2}\Phi(-x-y,x)e^{(x+y)/2}  =e^{x+y},   
\end{split}\end{equation*}
where the second equality follows from the duality identity and the 
third and fourth equalities both follow from the hexagon identity. 

{\it Second proof.} Let us set $\nu:=\mu_\Phi^{-1}$. 
Since $\mu_\Phi$ satisfies (\ref{Phimumu}), we have 
\begin{equation} \label{eq:mu}
\nu^{2,3} \circ \nu^{1,23} = \nu^{1,2} \circ \nu^{12,3} 
\circ \on{Ad}(\Phi(t_{12},t_{23})).  
\end{equation}
Let us set $C(x,y):= \nu(x+y)$, and apply (\ref{eq:mu}) 
to $x+y+z$ to obtain $C(x, C(y,z))=C(C(x,y), z)$.
According to \cite{AT}, this implies 
$C(x,y)=s^{-1}\on{log}(e^{sx}e^{sy})$ for some $s \in \kk^\times$. 
Checking degree 1 and 2 terms in $\nu$, we get $s=1$.

{\it Third proof.} As $\tilde\mu_{\bullet\bullet}(x_{01}) = e^{t_{01}}$, and using 
Proposition \ref{prop:3:2}, we get 
$\mu_{\bullet\otimes \bar O}(X_{1}...X_{n}) 
= \tilde\mu_{\bullet\otimes \bar O}(X_{1}...X_{n}) = (e^{t_{01}})^{0,1...n}
= e^{x_{1}+...+x_{n}}$. 
This implies $\mu_{\Phi}(XY)=e^{x+y}$ since $\mu_{\Phi} 
= \mu_{\bullet(\bullet\bullet)}$. 

\subsection{Proof that $J(\mu_\Phi)$ is a $\delta$-coboundary (end of proof of 
Theorem \ref{thm:main})} 
\label{delta:J:mu}

Since $J(\on{Ad}\Phi(t_{12},t_{23}))=0$, and $J(\mu_{\Phi}^{12,3}) = 
J(\mu_{\Phi})^{12,3}$, etc., we get by applying $J$ to (\ref{Phimumu}),
$$
\Phi(t_{12},t_{23})\cdot
J(\mu_{\Phi})^{12,3} + \Phi(t_{12},t_{23})\circ \mu_{\Phi}^{12,3}
\cdot J(\mu_{\Phi})^{1,2} = J(\mu_{\Phi})^{1,23} + \mu_{\Phi}^{1,23}\cdot 
J(\mu_{\Phi})^{2,3}.
$$
 
Applying the inverse of  (\ref{Phimumu}), we get 
$$
(\mu_{\Phi}^{1,2})^{-1} \circ (\mu_{\Phi}^{12,3})^{-1}\cdot
J(\mu_{\Phi})^{12,3} + (\mu_{\Phi}^{1,2})^{-1}
\cdot J(\mu_{\Phi})^{1,2} = 
(\mu_{\Phi}^{2,3})^{-1} \circ (\mu_{\Phi}^{1,23})^{-1}\cdot
J(\mu_{\Phi})^{1,23} + (\mu_{\Phi}^{2,3})^{-1} \cdot 
J(\mu_{\Phi})^{2,3},
$$ 
and since $a^{12,3}\cdot t^{12,3} = (a\cdot t)^{12,3}$, etc., 
$$
(\mu_{\Phi}^{1,2})^{-1}\cdot (\mu_{\Phi}^{-1}\cdot J(\mu_{\Phi}))^{12,3}
+ (\mu_{\Phi}^{-1}\cdot J(\mu_{\Phi}))^{1,2}
=  (\mu_{\Phi}^{2,3})^{-1}\cdot (\mu_{\Phi}^{-1}\cdot J(\mu_{\Phi}))^{1,23}
+ (\mu_{\Phi}^{-1}\cdot J(\mu_{\Phi}))^{2,3}.
$$ 

Now $\mu_{\Phi}^{-1}(x+y) = \on{log}(e^{x}e^{y})$ implies that 
$(\mu_{\Phi}^{1,2})^{-1}\cdot t^{12,3} = t^{\tilde{12},3}$, and similarly
with $1,23$, so $\tilde\delta(\mu_{\Phi}^{-1}\cdot J(\mu_{\Phi}))=0$. 
So there exists $\gamma\in\hat\Tr_{1}$ with valuation $\geq 2$
such that $\mu_{\Phi}^{-1}\cdot J(\mu_{\Phi})
=\tilde\delta(\gamma)$. Now $\mu_{\Phi}\cdot \gamma^{\tilde{12}} = \gamma^{12}$, 
and $\mu_{\Phi}\cdot \gamma^{1} = \gamma^{1}$, $\mu_{\Phi}\cdot \gamma^{2} 
= \gamma^{2}$
as $\mu_{\Phi}(x)\sim x$, $\mu_{\Phi}(y)\sim y$, therefore
$\mu_{\Phi}\cdot \tilde\delta(\gamma) = \delta(\gamma)$.
So $J(\mu_{\Phi}) = \delta(\gamma)$. 
It follows that for a suitable $\gamma\in u^{2}\kk[[u]]$, 
we have $J(\mu_{\Phi}) = \delta(\gamma) = \langle \gamma(x+y)
-\gamma(x)-\gamma(y)\rangle$. 

All this ends the construction of the map 
$M_{1}(\kk)\to \on{SolKV}(\kk)$, hence the proof of 
Theorem \ref{thm:main}. 

\subsection{Computation of $J(\mu_{\Phi})$ (proof of Proposition \ref{form:J})}

Let $U:= \lbr 1,A(x,y)\rbr\in \on{TAut}_{2}$, where 
$$\on{log}A(x,y) 
= \sum_{k\geq 1} \alpha_{k}(\on{ad}x)^{k}(y)+O(y^{2})
$$ 
(here $O(y^{2})$ means a series of elements with $y$-degree $\geq 2$). Then 
$\on{log}U=\lbr 0,\sum_{k\geq 1}\alpha_{k}(\on{ad}x)^{k}(y)+O(y^{2})\rbr$, 
and $J(U) = j(\on{log}U)+O(y^{2})$. Now $j(\on{log}U) = 
\langle \sum_{k\geq 1}\alpha_{k}y(-x)^{k}+O(y^{2})\rangle$. So  
$$
J(U) = \langle \sum_{k\geq 1}\alpha_{k}(-x)^{k}y\rangle + O(y^{2}). 
$$
On the other hand, the hexagon identity implies that 
$\mu_{\Phi} = \on{Inn}(\Phi(x,-x-y)e^{-x/2}) \circ \bar\mu_{\Phi}$, 
where $\bar\mu_{\Phi}=\lbr 1,\Phi(x,y)^{-1}\rbr$, and we then have 
$J(\bar\mu_{\Phi})=J(\mu_{\Phi})$. 

We have $\on{log}\Phi(x,y) = -\sum_{k\geq 1}\zeta_{\Phi}(k+1)
(\on{ad}x)^{k}(y) + O(y^{2})$, therefore 
$$
J(\mu_{\Phi}) = J(\bar\mu_{\Phi}) = \langle \sum_{k\geq 1}(-1)^k
\zeta_{\Phi}(k+1)x^{k}y\rangle+O(y^{2}). $$

As we have $J(\mu_{\Phi}) = \langle f(x)+f(y)-f(x+y)\rangle$ for some series 
$f(x)$, we get 
\begin{equation} \label{interm}
J(\mu_{\Phi}) = \langle (-1)^k{{\zeta_{\Phi}(k+1)}\over{k+1}}
((x+y)^{k+1}-x^{k+1}-y^{k+1})\rangle = \langle\on{log}\Gamma_{\Phi}(x)
+\on{log}\Gamma_{\Phi}(y)-\on{log}\Gamma_{\Phi}(x+y)\rangle. 
\end{equation}
This proves Proposition \ref{form:J}. 

\section{Group and torsor aspects} \label{sec:torsor}

\subsection{Group structures of $\on{KV}(\kk)$ and $\on{KRV}(\kk)$}

It is proved in \cite{AT} that $\on{KRV}(\kk)$ is a group, acting 
freely and transitively on $\on{SolKV}(\kk)$. 
 
Let us prove that $\on{KV}(\kk)$ is a group. For $\alpha\in \on{KV}(\kk)$, 
let $\sigma_{\alpha}:= \on{Duf}(\alpha)$, so $\sigma_{\alpha}\in u^{2}\kk[[u]]$, 
and $J(\alpha) = \tilde\delta(\sigma_\alpha)$. If $\alpha,\alpha'
\in \on{KV}(\kk)$, we have clearly $\alpha'\circ\alpha(X) \sim X$, 
$\alpha'\circ\alpha(Y) \sim Y$, $\alpha'\circ\alpha(XY) = XY$. 
Moreover, $J(\alpha'\circ\alpha) = J(\alpha') + \alpha'\cdot J(\alpha)
= \tilde\delta(\sigma_{\alpha'}) + \alpha'\cdot \tilde\delta(\sigma_\alpha)
= \tilde\delta(\sigma_{\alpha} + \sigma_{\alpha'})$, where the last equality 
follows from $\alpha'(X)\sim X$, $\alpha'(Y)\sim Y$, 
$\alpha'(XY)=XY$, which implies 
$\tilde\delta(\alpha'\cdot t) = \tilde\delta(t)$ for 
$t\in\hat\Tr_{1}$. So $\alpha'\circ
\alpha\in \on{KV}(\kk)$. One proves similarly that 
$\alpha^{-1}\in \on{KV}(\kk)$. We have also proved that 
$\sigma_{\alpha'\circ \alpha} = \sigma_\alpha+\sigma_{\alpha'}$, i.e., 
$\on{Duf}:\on{KV}(\kk)\to u^2\kk[[u]]$ is a group morphism. 

\subsection{The torsor structure of $\on{SolKV}(\kk)$ (proof of 
Proposition \ref{SolKV:torsor})} \label{sec:prop}

Let us prove that $\on{KV}(\kk)$ acts on $\on{SolKV}(\kk)$. 
For $\mu\in \on{SolKV}(\kk)$, let $r_{\mu}:= \on{Duf}(\mu)$, so 
$r_{\mu}\in u^{2}\kk[[u]]$, and $J(\mu) = \delta(r_\mu)$. 
For 
$\mu\in\on{SolKV}(\kk)$, $\alpha\in\on{KV}(\kk)$, we have 
$\mu\circ\alpha(X) \sim \mu(X)\sim e^x$, 
$\mu\circ\alpha(Y) \sim \mu(Y)\sim e^y$, 
$\mu\circ\alpha(XY) = \mu(XY) = e^{x+y}$. Moreover, 
$J(\mu\circ\alpha) = J(\mu)  + \mu\cdot J(\alpha) = 
\delta(r_\mu) + \mu\cdot \tilde\delta(\sigma_\alpha) = 
\delta(r_\mu+\sigma_\alpha)$, where the last equality uses the identity 
$\delta(t) = \mu\cdot \tilde\delta(t)$ for $t\in\hat\Tr_{2}$, 
which follows from $\mu(XY)=e^{x+y}$, $\mu(X)\sim e^{x}$, $\mu(Y)
\sim e^{y}$. So $\mu\circ\alpha
\in \on{SolKV}(\kk)$. We have also proved that $r_{\mu\circ\alpha} = 
r_\mu + \sigma_\alpha$, so $\on{Duf}:\on{SolKV}(\kk)\to u^2\kk[[u]]$
is a morphism of torsors.  

Let us now prove that the action of $\on{KV}(\kk)$ on $\on{SolKV}(\kk)$
is free and transitive. For $\mu,\mu'\in \on{SolKV}(\kk)$, set $\alpha:= 
\mu^{-1}\circ\mu'$; 
then $\alpha(X)\sim X$, $\alpha(Y)\sim Y$, $\alpha(XY)=XY$, 
and $J(\alpha) = J(\mu^{-1}) + \mu^{-1}\cdot J(\mu') = 
\mu^{-1}\cdot (J(\mu')-J(\mu))$
as $J(\mu^{-1}) = -\mu^{-1}\cdot J(\mu)$. Then $J(\alpha) = 
\mu^{-1}\cdot (\delta(r_{\mu'} - r_{\mu})) = 
\tilde\delta(r_{\mu'}-r_{\mu})$, where the last equality uses 
$\mu^{-1}\cdot \delta(t) = \tilde\delta(t)$ for $t\in\hat\Tr_{1}$. 
So $\alpha\in\on{KV}(\kk)$. 

\subsection{Compatibilities of morphisms with group structures and actions
(proof of Theorem \ref{thm:morph:tors})}
\label{sec:end}

We now show that: (a) $f\mapsto \alpha_{f}^{-1}$ is a group 
morphism  $\on{GT}_{1}(\kk)\to \on{KV}(\kk)$, (b) $g\mapsto 
a_{g}^{-1}$ is a group morphism $\on{GRT}_{1}(\kk)\to \on{KRV}(\kk)$, 
(c) the map $\Phi\mapsto\mu_{\Phi}$ is compatible with the actions of 
these groups. 

For this, we will show that
\begin{equation} \label{eq:torsors}
\mu_{f*\Phi} = \mu_{\Phi}\circ\alpha_{f}, \quad  
\mu_{\Phi*g}=a_{g}\circ \mu_{\Phi}.
\end{equation}
We will check these identities 
on the first generator ($X$ or $x$), the proofs in the second case being 
similar. 

The proofs go as follows: 
\begin{equation*} \begin{split}
& \mu_{f*\Phi}(X) = (f*\Phi)(x,-x-y) \cdot e^{x}\cdot (same)^{-1} \\
& =f(\Phi(x,-x-y)e^{x}\Phi(x,-x-y)^{-1},e^{-x-y})\Phi(x,-x-y)\cdot e^{x}
\cdot (same)^{-1} \\
& =f(\mu_{\Phi}(X),\mu_{\Phi}(Y^{-1}X^{-1}))\cdot \mu_{\Phi}(X)
\cdot (same)^{-1} \\ 
& = \mu_{\Phi}(f(X,Y^{-1}X^{-1})\cdot X\cdot (same)^{-1})
= \mu_{\Phi}\circ \alpha_{f}(X)
\end{split}\end{equation*} 
and 
\begin{equation*} \begin{split}
& \mu_{\Phi*g}(X) = (\Phi*g)(x,-x-y)\cdot e^{x}\cdot (same)^{-1} \\
& =\Phi(g(x,-x-y)x g(x,-x-y)^{-1},-x-y)g(x,-x-y)\cdot e^{x}\cdot (same)^{-1}\\
& =\Phi(a_{g}(x),a_{g}(-x-y))\cdot a_{g}(x)\cdot (same)^{-1}\\
& =a_{g}(\Phi(x,-x-y)x \Phi(x,-x-y)^{-1}) = a_{g}\circ \mu_{\Phi}(X). 
\end{split}\end{equation*}

The first part of (\ref{eq:torsors}) implies the following: 
(a) if $f\in \on{GT}_{1}(\kk)$, then $\alpha_{f}\in \on{KV}(\kk)$; 
(b) $\alpha_{f_{1}*f_{2}} = \alpha_{f_{2}}\circ \alpha_{f_{1}}$; 
(c) $M_{1}(\kk)\to \on{SolKV}(\kk)$ is compatible with the 
group morphism $f\mapsto \alpha_{f}^{-1}$. 

Indeed, using the nonemptinesss of $M_{1}(\kk)$ (see \cite{Dr:Gal})
we get $\alpha_{f} = \mu_{\Phi}^{-1}\circ \mu_{f*\Phi}$, which implies 
$\alpha_{f}\in \on{KV}(\kk)$ according to Subsection \ref{sec:prop}, i.e., (a).
Again using the nonemptiness of  $M_{1}(\kk)$, we get 
$\alpha_{f_{1}*f_{2}} = \mu_{\Phi}^{-1} \circ \mu_{(f_{1}*f_{2})*\Phi}
= (\mu_{\Phi}^{-1} \circ \mu_{f_{2}*\Phi}) \circ 
( \mu_{f_{2}*\Phi}^{-1} \circ \mu_{f_{1}*(f_{2}*\Phi)}) = 
\alpha_{f_{2}} \circ \alpha_{f_{1}}$ (where we used 
$(f_{1}*f_{2})*\Phi = f_{1}*(f_{2}*\Phi)$), which proves (b).
(c) is then tautological.  

Similarly, the second part of (\ref{eq:torsors}) implies: 
(a) if $g\in \on{GRT}_{1}(\kk)$, then $a_{g}\in \on{KRV}(\kk)$; 
(b) $a_{g_{1}*g_{2}} = a_{g_{2}}\circ a_{g_{1}}$; 
(c) $M_{1}(\kk)\to \on{SolKV}(\kk)$ is compatible with the 
group morphism $g\mapsto a_{g}^{-1}$. 
All this proves Theorem \ref{thm:morph:tors}. 

It is easy to prove the identities $\alpha_{f_{1}*f_{2}} = \alpha_{f_{2}} \circ 
\alpha_{f_{1}}$, $a_{g_{1}*g_{2}} = a_{g_{2}}\circ a_{g_{1}}$ directly 
(i.e., not using the nonemptiness of $M_{1}(\kk)$):  the verifications on the first generators
($X$ and $x$) are 
\begin{equation*} \begin{split}
& \alpha_{f_{1}*f_{2}}(X)=(f_{1}*f_{2})(X,Y^{-1}X^{-1})\cdot X\cdot (same)^{-1}
\\
& =f_{1}(f_{2}(X,Y^{-1}X^{-1})Xf_{2}(X,Y^{-1}X^{-1})^{-1},Y^{-1}X^{-1})
f_{2}(X,Y^{-1}X^{-1})\cdot X\cdot (same)^{-1}\\
& =f_{1}(\alpha_{f_{2}}(X),\alpha_{f_{2}}(Y^{-1}X^{-1}))\cdot 
\alpha_{f_{2}}(X)\cdot (same)^{-1}\\
& =\alpha_{f_{2}}(f_{1}(X,Y^{-1}X^{-1})\cdot X\cdot (same)^{-1})
=\alpha_{f_{2}}\circ \alpha_{f_{1}}(X),  
\end{split}\end{equation*} 
and 
\begin{equation*} \begin{split}
& a_{g_{1}*g_{2}}(x) = (g_{1}*g_{2})(x,-x-y)\cdot x\cdot (same)^{-1}\\
& =g_{1}(g_{2}(x,-x-y)x g_{2}(x,-x-y)^{-1},-x-y)g_{2}(x,-x-y)\cdot 
x\cdot (same)^{-1} \\
& =g_{1}(a_{g_{2}}(x),a_{g_{2}}(-x-y))\cdot a_{g_{2}}(x)\cdot (same)^{-1}\\
& =a_{g_{2}}(g_{1}(x,-x-y)xg_{1}(x,-x-y)^{-1}) = a_{g_{2}}\circ a_{g_{1}}(x). 
\end{split}\end{equation*}

\begin{remark}
The Lie algebra morphism corresponding to $g\mapsto a_{g}^{-1}$
is the morphism $\nu : \grt_{1}\to \mathfrak{krv}$ from \cite{AT}, given by 
$\psi(x,y)\mapsto\lbr \psi(x,-x-y),\psi(y,-x-y)\rbr$. 
\end{remark}

\subsection{Torsor properties of the Duflo formal series (proof of Proposition 
\ref{form:J:torsor})} \label{sec:J}

We have already proved that $M_{1}(\kk)\to \on{SolKV}(\kk)$, 
and $\on{SolKV}(\kk)\stackrel{\on{Duf}}{\to} u^{2}\kk[[u]]$ is a 
morphism of torsors. On the other hand, it follows from \cite{E} that 
$M_{1}(\kk)\stackrel{\Phi\mapsto \on{log}\Gamma_{\Phi}}{\to} 
\{r\in u^{2}\kk[[u]] | r_{ev}(u)= - {{u^{2}}\over 24}+...\}$
is a morphism of torsors and from Proposition \ref{form:J} that the diagram of
Proposition \ref{form:J:torsor}) commutes. 

For later use, let us make the group morphism $\on{GT}_{1}(\kk)
\to u^{3}\kk[[u^{2}]]$ underlying $\Phi\mapsto \on{log}\Gamma_{\Phi}$
explicit. 

\begin{lemma}
For $f\in\on{GT}_{1}(\kk)$, there is a unique $\Gamma_{f}\in 
\on{exp}(u^{3}\kk[[u^{2}]])$ such that 
$$
[\on{log}f(e^{a},e^{b})] = 1 - {{\Gamma_{f}(-\overline a)\Gamma_{f}(-\overline b)}
\over {\Gamma_{f}(-\overline a - \overline b)}};  
$$
here we use the isomorphism $\hat\f_{2}'/\hat\f_{2}''\simeq \overline a
\overline b \kk[[\overline a,\overline b]]$ given by 
(class of $(\on{ad}a)^{k}(\on{ad}b)^{l}([a,b]))\leftrightarrow 
\overline a^{k+1}\overline b^{l+1}$.
The map $\on{GT}_{1}(\kk)\to u^{3}\kk[[u^{2}]]$, 
$f\mapsto \on{log}\Gamma_{f}$ is a group morphism and $\Gamma_{f*\Phi} = 
\Gamma_{f}\Gamma_{\Phi}$ for any $f\in \on{GT}_{1}(\kk)$, $\Phi\in 
M_{1}(\kk)$. 
\end{lemma}

{\em Proof.} The map $\f_{2}\to \kk[\overline a,\overline b]$, 
$\psi\mapsto (b\partial_{b}\psi)^{ab}$ also induces an isomorphism 
$\hat\f_{2}'/\hat\f_{2}''\simeq \overline a\overline b
\kk[[\overline a,\overline b]]$, 
which takes 
the class $(\on{ad}a)^{k}(\on{ad}b)^{l}([a,b])$ to 
$(-1)^{k+l+1}\overline a^{k+1}\overline b^{l+1}$. 
So for $\psi\in\hat\f_{2}'$, we have 
$(b\partial_{b}\psi)^{ab}(\bar a,\bar b) = -[\psi](-\bar a,-\bar b)$
(where $\psi\mapsto [\psi]$ is the map $\hat\f'_{2}\to 
\hat\f'_{2}/\hat\f''_{2}\simeq \overline a\overline b
\kk[[\overline a,\overline b]]$). 

So (\ref{Phi:Gamma}) may be rewritten 
$$
[\on{log}\Phi](\overline a,\overline b) = 1 - 
{{\Gamma_{\Phi}(-\overline a-\overline b)}\over{
\Gamma_{\Phi}(-\overline a)\Gamma_{\Phi}(-\overline b)}}. 
$$

If now $\psi,\alpha\in\hat\f'_{2}$, we have 
$\psi(e^{-\alpha}ae^{\alpha},b)\in \hat\f'_{2}$ and 
$[\psi(e^{-\alpha}ae^{\alpha},b)] = (1- [\alpha(a,b)])[\psi(a,b)]$. 
Indeed, when $\psi(a,b) = (\on{ad}a)^{k}(\on{ad}b)^{l}([a,b])$, 
one checks that the part of $\psi(e^{-\alpha}ae^{\alpha},b)$
containing $\alpha$ more than twice lies in $\hat\f_{2}''$, 
and the part containing it once has the same class as 
$(\on{ad}a)^{k}(\on{ad}b)^{l}([[-\alpha,a],b])$. 

If now $f\in \on{GT}_{1}(\kk)$, we have 
$(f*\Phi)(a,b) = \Phi(a,b)f(\Phi^{-1}(a,b)e^{a}\Phi(a,b),e^{b})$, 
so 
\begin{equation*}
\begin{split}[\on{log} (f*\Phi)(a,b)] & = 
[\on{log}\Phi(a,b)] + [\on{log}
f(\Phi^{-1}(a,b)e^{a}\Phi(a,b),e^{b})]
\\ & 
= [\on{log}\Phi(a,b)] + [\on{log}f(e^{a},e^{b})]
- [\on{log}\Phi(a,b)] [\on{log}f(e^{a},e^{b})]. 
\end{split}\end{equation*}
so 
\begin{equation} \label{id:log:Phi}
1 - [\on{log}(f*\Phi)(a,b)] = (1 - [\on{log}\Phi(a,b)])
(1 - [\on{log}f(e^{a},e^{b})]). 
\end{equation}
If fix $\Phi_{0}\in M_{1}(\kk)$ and set 
$\Gamma_{f}(u):= \Gamma_{f*\Phi_{0}}(u)/\Gamma_{\Phi_{0}}(u)$, 
then we get 
$$
1 - [\on{log}f(e^{a},e^{b})] = {{\Gamma_{f}(-\overline a)\Gamma_{f}(-\overline b)}
\over {\Gamma_{f}(-\overline a - \overline b)}} 
$$
as wanted. Moreover, (\ref{id:log:Phi}) implies that $\Gamma_{f*\Phi} = 
\Gamma_{f}\Gamma_{\Phi}$, which also implies that $f\mapsto \Gamma_{f}$
is a group morphism. 
\hfill \qed\medskip 

\section{Direct construction of the map $\on{GT}_{1}(\kk) \to \on{KV}(\kk)$}
\label{act:gt}

We will now sketch a proof  of $(f\in \on{GT}_{1}(\kk)) \Rightarrow 
(\alpha_{f}\in \on{KV}(\kk))$, independent of the nonemptiness of  
$M_{1}(\kk)$. 

\subsection {Action of $\on{GT}_{1}(\kk)$ on completed braid groups}

Let $\cC$ be a b.m.c. We denote by $\beta_{X,Y}:X\otimes Y\to Y\otimes X$
and $a_{X,Y,Z}:(X\otimes Y)\otimes Z\to X\otimes (Y\otimes Z)$
the braiding and associativity constraints. For $O\in \on{Ob}({\bf PaB})$ of length $n$ 
and any $X_{1},...,X_{n}\in \on{Ob}(\cC)$, we construct 
the tensor product $O(X_{1},...,X_{n})$ of 
$X_{1},...,X_{n}$ with parenthesization $O$. We say that $\cC$ is 
prounipotent if for any $X_{1},...,X_{n}$ and any $O$, the image of 
$\on{PB}_{n}\to \on{Aut}_{\cC}(O(X_{1},...,X_{n}))$ is prounipotent
(it suffices to require this for a given $O$). 
If $\cC$ is a prounipotent b.m.c. and $f\in \on{GT}_{1}(\kk)$, we construct 
a new b.m.c. $^{f}\!\cC$ as follows: $^{f}\!\cC$ is the same as $\cC$ 
at the level of objects and morphisms, the composition and the 
tensor product of morphisms are not modified,  
but the braiding and associativity constraints are modified as follows: 
$$
\beta'_{X,Y} = \beta_{X,Y}, \quad a'_{X,Y,Z} = a_{X,Y,Z} \circ 
f(\beta_{YX}\beta_{XY},a_{X,Y,Z}^{-1}
\circ \beta_{ZY}\beta_{YZ}\circ a_{X,Y,Z}). 
$$
We then have $^{f_{1}}(^{f_{2}}\cC) = \ ^{f_{1}*f_{2}}\cC$. 
Moreover, the action of $\on{GT}_{1}(\kk)$ on $\on{BMC}$ 
is functorial, so a tensor functor $\phi:\cC\to\cD$ and $f\in \on{GT}_{1}(\kk)$
give rise to $^{f}\!\phi:\ ^{f}\!\cC\to \ ^{f}\!\cD$. 
Note that for $O,O'\in\on{Ob}(\cC)$, 
and under the identifications $^{f}\!\cC(O,O')=\cC(O,O')$, 
$^{f}\!\cD(\phi(O),\phi(O')) = \cD(\phi(O),\phi(O'))$, the map 
$^{f}\!\phi(O,O'):\ ^{f}\!\cC(O,O')\to \ ^{f}\!\cD(\phi(O),\phi(O'))$ coincides with 
$\phi(O,O'):\cC(O,O')\to \cD(\phi(O),\phi(O'))$.

Let ${\bf PaB}_{\kk}$ be the completion of 
${\bf PaB}$ obtained by replacing each group 
$\on{B}_{n}$ by its completion $\on{B}_{n}(S_{n},\kk)$ relative to 
the morphism $\on{B}_{n}\to S_{n}$. By universal properties, 
we have a unique morphism $\phi_{f}:{\bf PaB}\to \ ^{f}\!{\bf PaB}$
which is the identity on objects. If then $f_{1},f_{2}\in \on{GT}_{1}(\kk)$, 
we have 
\begin{equation} \label{smallid}
^{f_{1}}\!\phi_{f_{2}}\circ \phi_{f_{1}} = \phi_{f_{1}*f_{2}}; 
\end{equation}
indeed, both terms are tensor functors ${\bf PaB}_{\kk}\to 
\ ^{f_{1}*f_{2}}\!{\bf PaB}_{\kk}$ which are the identity on objects. 

If now $O\in \on{Ob}({\bf PaB})$ has length $n$, $\phi_{f}$ gives rise
to a group morphism $\phi_{f}(O):{\bf PaB}_{\kk}(O)\to 
\ ^{f}\!{\bf PaB}_{\kk}(O)$. We denote by 
$$
\tilde\alpha_{f}^{O} : \on{B}_{n}(S_{n},\kk) \to \on{B}_{n}(S_{n},\kk) 
$$
the group endomorphism derived from $\phi_{f}(O)$ and the identifications
${\bf PaB}_{\kk}(O) = \ ^{f}\!{\bf PaB}_{\kk}(O) = \on{B}_{n}(S_{n},\kk)$. 
Identity (\ref{smallid}) and the identification of $^{f_{1}}\!\tilde\alpha_{f_{2}}^{O}$
with $\tilde\alpha_{f_{2}}^{O}$ imply 
$$
\tilde\alpha_{f_{2}}^{O} \circ \tilde\alpha_{f_{1}}^{O} = 
\tilde\alpha_{f_{1}*f_{2}}^{O}, 
$$
so we have  a group antimorphism $\on{GT}_{1}(\kk)\to 
\on{Aut}(\on{B}_{n}(S_{n},\kk))$, $f\mapsto \tilde\alpha_{f}^{O}$. 

It is easy to see that we have a commutative diagram 
$\begin{matrix} \on{B}_{n}(S_{n},\kk) & \stackrel{\tilde\alpha_{f}^{O}}{\to} & 
\on{B}_{n}(S_{n},\kk) \\
\searrow & & \swarrow\\
 & S_{n} & \end{matrix}$
 so $\tilde\alpha_{f}^{O}$ restricts to an automorphism 
 $\tilde\alpha_{f}^{O}\in \on{Aut}(\on{PB}_{n}(\kk))$. 

If now $O,O'\in\on{Ob}({\bf PaB})$ have length $n$, 
then $\on{can}_{O,O'}\in {\bf PaB}_{\kk}(O,O')$ is the morphism 
corresponding to the trivial braid. Then $\phi_{f}(\on{can}_{O,O'})
\circ \on{can}_{O,O'}^{-1}\in {\bf PaB}_{\kk}(O)$. Let 
$f^{O,O'}\in \on{PB}_{n}(\kk)$ be the image of this 
element. Since the diagram 
$\begin{matrix} {\bf PaB}(O) & \stackrel{x\mapsto 
\on{can}_{O,O'}\circ x \circ \on{can}_{O,O'}^{-1}}{\to} & 
{\bf PaB}(O') \\
\searrow & & \swarrow\\
& \on{B}_{n}(S_{n},\kk) & \end{matrix}$
commutes, we have 
\begin{equation} \label{alpha:O:O'}
\tilde\alpha_{f}^{O'} = \on{Inn}(f^{O,O'}) 
\circ  \tilde\alpha_{f}^{O}. 
\end{equation}

\subsection{Actions of  $\on{GT}_{1}(\kk)$ on free groups}

Let us index the generators of $\on{PB}_{n}(\kk)$
by $x_{ij}$, $0\leq i<j\leq n-1$. Recall that the subgroup of $\on{PB}_{n}(\kk)$
generated by $x_{01},...,x_{0,n-1}$ is isomorphic to $\on{F}_{n-1}(\kk)$. 
We set $X_{i}=x_{0i}$ for $i=1,...,n-1$. 

\begin{proposition}
Each $\tilde\alpha_{f}^{O}$ restricts to an automorphism 
$\alpha_{f}^{O}\in \on{Aut}(\on{F}_{n-1}(\kk))$, such that 
for any $i$, $\alpha_{f}^{O}(X_{i})\sim X_{i}$. 
\end{proposition}

{\em Proof.} Let us index the letters of $O$ by $0,...,n-1$. For $i=1,...,n-1$, 
let $O_{i}$ be an object of ${\bf PaB}$ of length $n$, in which the letters 
$i-1$ and $i$ appear as $...(\bullet\bullet)....$. We have $X_{i} = 
(\sigma_{0}...\sigma_{i-2})^{-1}\sigma_{i-1}^{2}\sigma_{0}...\sigma_{i-2}$. 
We have $\tilde\alpha_{f}^{O}(\sigma_{0}...\sigma_{i-2}) = 
\sigma_{0}...\sigma_{i-2}\cdot p_{i}$, where $p_{i}\in 
\on{PB}_{n}(\kk)$. On the other hand, 
$\tilde\alpha_{f}^{O}(\sigma_{i-1}) = f^{O,O_{i}}
\tilde\alpha_{f}^{O}(\sigma_{i-1}) (f^{O,O_{i}})^{-1}$
and $\tilde\alpha_{f}^{O}(\sigma_{i-1})  = \sigma_{i-1}$
as $ \on{B}_{n} \simeq {\bf PaB}(O_{i})$ takes $\sigma_{i-1}$ to
$\on{id}_{\bullet}^{\otimes i-1} \otimes \beta_{\bullet,\bullet}\otimes 
\on{id}_{\bullet}^{\otimes n-i-2}$. So 
$$\tilde\alpha_{f}^{O}(\sigma_{i-1}^{2}) = f^{O,O_{i}}
\sigma_{i-1}^{2} (f^{O,O_{i}})^{-1}, 
$$
with $\alpha^{O,O_{i}}_{f}\in \on{PB}_{n}(\kk)$. Then 
\begin{equation} \begin{split}
\tilde\alpha^{O}_{f}(X_{i}) & = 
(\sigma_{0}...\sigma_{i-2}p_{i})^{-1}f^{O,O_{i}}
\sigma_{i-1}^{2} ( f^{O,O_{i}})^{-1}
\sigma_{0}...\sigma_{i-2}p_{i} \\
 & = 
p_{i}^{-1} (\sigma_{0}...\sigma_{i-2})^{-1}
f^{O,O_{i}} (\sigma_{0}...\sigma_{i-2})
\cdot X_{i}\cdot (same)^{-1}.
\end{split}\end{equation}
As $p_{i}^{-1} (\sigma_{0}...\sigma_{i-2})^{-1}
f^{O,O_{i}} (\sigma_{0}...\sigma_{i-2})$
belongs to $\on{PB}_{n}(\kk)$, and as $\on{PB}_{n}(\kk)$
acts on $\on{F}_{n-1}(\kk)$ by tangential automorphisms, 
we obtain that $\tilde\alpha_{f}^{O}(X_{i})$ lies in 
$\on{F}_{n-1}(\kk)$ and is conjugated in $\on{F}_{n-1}(\kk)$
to $X_{i}$. \hfill \qed\medskip 

Similarly to Proposition \ref{prop:3:2}, one can prove: 

\begin{proposition} \label{prop:add}
If $O=\bullet\otimes \bar O$, where $\bar O\in \on{Ob}({\bf PaB})$, 
then $\alpha_{f}^{O}(X_{1}...X_{n-1})=X_{1}...X_{n-1}$. 
\end{proposition}


We then have 
$$
\alpha_{f}^{O'} = \on{Ad}(f^{O,O'}) \circ \alpha_{f}^{O}; 
$$
this is an identity in $\on{Aut}(\on{F}_{n-1}(\kk))$, where 
$ \on{Ad}(\alpha_{f}^{O,O'})$ is not necessarily inner. 

We also record the identities 
\begin{equation} \label{mu:f:Phi}
\tilde\mu_{f*\Phi}^{O} = \tilde\mu_{\Phi}^{O} \circ \alpha_{f}^{O},
\quad  
\mu_{f*\Phi}^{O} = \mu_{\Phi}^{O} \circ \alpha_{f}^{O}.
\end{equation}

\subsection{The map $\on{GT}_{1}(\kk)\to \on{KV}(\kk)$}
 \label{map:gt}
 
Let us fix an element $f\in \on{GT}_{1}(\kk)$ and 
denote $\tilde\alpha_{f}^{O}$, $\alpha_{f}^{O}$ simply by
$\tilde\alpha_{O}$, $\alpha_{O}$. 

As in Subsection \ref{rel:braid}, one proves that 
\begin{equation} \label{diag:GT}
\begin{matrix} 
\on{PB}_{n}(\kk) & \stackrel{1,2,...,\widetilde{ii+1},...,n}{\to}& 
\on{PB}_{n+1}(\kk)\\
\scriptstyle{\alpha_{O}} \downarrow & & \downarrow
\scriptstyle{\alpha_{O^{(i)}}}\\
\on{PB}_{n}(\kk) & \stackrel{1,2,...,\widetilde{ii+1},...,n}{\to}& 
\on{PB}_{n+1}(\kk)
\end{matrix}
\end{equation}
commutes. Using Proposition \ref{app:centr}, one then proves 
\begin{equation} \label{diag:GT:bis}
\alpha_{O^{(i)}} = 
\alpha_{O}^{1,...,\widetilde{ii+1},...,n} \circ 
\alpha_{\bullet(\bullet\bullet)}^{i,i+1}.
\end{equation}
Similarly to Proposition \ref{prop:mu}, one proves that 

1) $\alpha_{\bullet(\bullet\bullet)} = \alpha_{f}$. 

2) $f^{\bullet((\bullet\bullet)\bullet),\bullet(\bullet(\bullet\bullet))}
=f(x_{12},x_{23})$. 

As in Subsection \ref{pf:Phimumu}, one proves that this implies 
\begin{equation} \label{just:proved}
\on{Ad}f(x_{12},x_{23})\circ \alpha_{f}^{\widetilde{12},3}
\circ \alpha_{f}^{1,2} = \alpha_{f}^{1,\widetilde{23}}\circ \alpha_{f}^{2,3}.
\end{equation}
As in Subsection \ref{X:Y}, one can give three proofs of the fact that 
$\alpha_{f}(XY)=XY$. 
Similarly to Subsection \ref{delta:J:mu}, one then proves that 
identity (\ref{just:proved}) then implies that $J(\alpha_{f})$ is a 
$\tilde\delta$-coboundary. 

Let us explain this proof in some detail. Since $J(\on{Ad}
f(x_{12},x_{23}))=0$
and $J(\alpha_{f}^{\widetilde{12},3})=J(\alpha_{f})^{\widetilde{12},3}$, 
we get by applying $J$ to (\ref{just:proved}) 
$$
\on{Ad}f(x_{12},x_{23}) \cdot J(\alpha_{f})^{\widetilde{12},3} + 
\big( \on{Ad}f(x_{12},x_{23})\circ \alpha_{f}^{\widetilde{12},3} \big)\cdot
J(\alpha_{f})^{1,2} = J(\alpha_{f})^{\widetilde{12},3} + J(\alpha_{f})^{2,3}.
$$ 
Applying the inverse of (\ref{just:proved}), we get 
$$
(\alpha_{f}^{1,2})^{-1}\cdot 
(\alpha_{f}^{-1}\cdot J(\alpha_{f}))^{\widetilde{12},3} + 
(\alpha_{f}^{-1}\cdot J(\alpha_{f}))^{1,2} = 
(\alpha_{f}^{2,3})^{-1}\cdot 
(\alpha_{f}^{-1}\cdot J(\alpha_{f}))^{1,\widetilde{23}} + 
(\alpha_{f}^{-1}\cdot J(\alpha_{f}))^{2,3}
$$
Now $\alpha_{f}(XY)=XY$ implies that $\alpha_{f}^{1,2}\cdot t^{\widetilde{12},3}
= t^{\widetilde{12},3}$ and similarly with $1,\tilde{23}$, so 
$\tilde\delta(\alpha_{f}^{-1}\cdot J(\alpha_{f}))=0$. 
As $\hat\Tr_{1}\stackrel{\tilde\delta}{\to}\hat\Tr_{2}\to...$
is acyclic in degree 2, there exists $\beta\in\hat\Tr_{1}$ with valuation $\geq 2$
such that  
$\alpha_{f}^{-1}\cdot J(\alpha_{f}) = \tilde\delta(\beta)$, so 
$J(\alpha_{f}) = \alpha_{f}\cdot \tilde\delta(\beta)$. Now 
$\alpha_{f}(XY)=XY$, $\alpha_{f}(X)\sim X$, $\alpha_{f}(Y)\sim Y$
imply that $\alpha_{f}\cdot \tilde\delta(\beta) = \tilde\delta(\beta)$, 
so $J(\alpha_{f}) = \tilde\delta(\beta)$. 
It follows that $J(\alpha_{f})$ has the form 
$\tilde\delta(\beta)=\langle\beta(\on{log}(e^{x}e^{y})) - \beta(x)
-\beta(y)\rangle$.  

\begin{remark}
(\ref{just:proved}) can also be proved directly, checking the identity on
each of the generators of $\on{F}_{3}(\kk)$ and using only the duality, hexagon and
pentagon relations. This proof then extends to the profinite and pro-$l$ cases. 
\end{remark}

\section{The Jacobians of $\mu_{\Phi,O}$ and $\alpha_f^{O}$}
\label{sec:jac}

\subsection{Telescopic formulas}

If $O\in \on{Ob}({\bf{PaB}})$ has the form $O=\bullet\otimes O'$, 
with $|O'|=n$, then one proves by using (\ref{muO:muOnew}) that 
$\mu_{O}$ expresses directly in terms of $\mu_{\Phi}$, 
for example
$$
\mu_{\bullet((((\bullet\bullet)(\bullet\bullet))(\bullet(\bullet\bullet)))
(\bullet\bullet))} = 
\mu_{\Phi}^{1234567,89}\mu_{\Phi}^{1234,567}\mu_{\Phi}^{8,9}
\mu_{\Phi}^{12,34}\mu_{\Phi}^{5,67}
\mu_{\Phi}^{1,2}\mu_{\Phi}^{3,4}\mu_{\Phi}^{6,7}. $$
The general formula is 
$$
\mu_{\bullet\otimes O'} = \prod_{n\geq 0}
\prod_{\nu \in N(T'), d(\nu)=n}\mu_{\Phi}^{L(\nu),R(\nu)};
$$ 
here $T'$ is the binary planar rooted tree underlying $O'$; $N(T')$ is the set of 
its nodes; $d(\nu)$ is the degree of $\nu$ (distance to the root of the tree); 
$L(\nu)$, $R(\nu)$ is the set of left and right leaves of $\nu$ (these are 
disjoints subsets of $\{1,...,n\}$). The first product is
taken according to increasing values of $n$ (the order in the second product 
does not matter as the arguments of this product commute with each other). 
Here is the tree corresponding to the above example (Figure 
\ref{AETfig9.eps}): 

\begin{figure}[h!]
\begin{center}
\includegraphics[width=10cm]{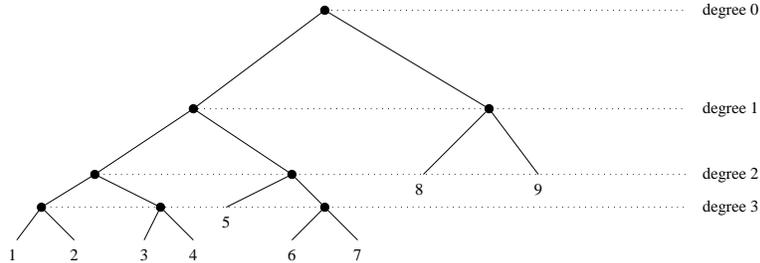}
\caption{\footnotesize There are 8 nodes}\label{AETfig9.eps}
\end{center}
\end{figure}

Similarly, using (\ref{diag:GT:bis}), one proves that for $f\in \on{GT}_{1}(\kk)$, 
we have 
$$
\alpha_{f}^{\bullet\otimes O'} = \prod_{n\geq 0}
\prod_{\nu \in N(T'), d(\nu)=n}\alpha_{f}^{\widetilde{L(\nu)},
\widetilde{R(\nu)}}. 
$$ 

\subsection{Computation of Jacobians}

%
%
%

Let $\mu_{n}:= \mu_{\bullet(\bullet...(\bullet\bullet))}$. Then:  

\begin{proposition}


$J(\mu_{n}) = \langle \sum_{i=1}^{n} \on{log}\Gamma_{\Phi}(x_{i}) 
-\on{log}\Gamma_{\Phi}(\sum_{i=1}^{n}x_{i})\rangle$. 
\end{proposition}

(We identified $\mu_{n}$ with its composition with $e^{x_{i}}\mapsto X_{i}$, 
which belongs to $\on{TAut}_{n}$.)

{\em Proof.} 
We have $\mu_{n} = \mu_{\Phi}^{1,2...n}\circ \mu_{\Phi}^{2,3...n}\circ ... 
\circ \mu_{\Phi}^{n-1,n}$. 
One then proves by descending induction 
on $k$ that $J(\mu_{\Phi}^{k,k+1...n}\circ...\circ\mu_{\Phi}^{n-1,n})
= \langle \sum_{i=k}^{n}\on{log}\Gamma_{\Phi}(x_{i}) 
- \on{log}\Gamma_{\Phi}(\sum_{i=k}^{n}x_{i})\rangle$, using the fact that 
the action of $\mu_{\Phi}^{k,k+1...n}$ on the various $\langle 
\on{log}\Gamma_{\Phi}(x_{i})\rangle$ as well as on 
$\langle \on{log}\Gamma_{\Phi}(\sum_{i=k}^{n}x_{i})\rangle$ 
is trivial. \hfill \qed\medskip 

If now $O\in \on{Ob}({\bf PaB})$ is arbitrary with with $|O|=n+1$, then: 

\begin{proposition}
$J(\mu_{\Phi,O}) = J(\mu_{n}) = \langle \sum_{i=1}^{n} \on{log}\Gamma_{\Phi}(x_{i}) 
-\on{log}\Gamma_{\Phi}(\sum_{i=1}^{n}x_{i})\rangle$. 
\end{proposition}

{\em Proof.} We have $\mu_{O} = \on{Ad}\Phi_{O_{n},O} \circ \mu_{n}$,
where $O_{n}=\bullet(...(\bullet\bullet))$. We then use the
cocycle property of $J$, the above formula for $J(\mu_{n})$, the 
fact that $J(\on{Ad}g)=0$ for $g\in \on{exp}(\hat\t_{n+1})$, and the 
following lemma: 

\begin{lemma} \label{lemma:special}
If $g\in\on{exp}(\hat\t_{n+1})$, then $(\on{Ad}g)(x_{1}+...+x_{n})
\sim x_{1}+...+x_{n}$. 
\end{lemma}

{\em Proof of Lemma.} Decompose $a\in\t_{n+1}$ as $a_{0}+a_{1}^{1,2,...,n}$, with 
$a_{0}\in \f_{n}$ and $a_{1}\in \t_{n}$ (the map $a_{1}\mapsto a_{1}^{1,2,...,n}$
is the injection $\t_{n}\to \t_{n+1}$, $t_{ij}\mapsto t_{ij}$). Then 
$[t_{ij},x_{1}+...+x_{n}]=0$ for $i,j\in\{1,...,n\}$, so $[a_{1}^{1,2,...,n},x_{1}+
...+x_{n}]=0$, so $[a,x_{1}+...+x_{n}] = [a_{0},x_{1}+...+x_{n}]$. 
It follows that if $g\in \on{exp}(\hat\t_{n+1})$, there exists 
$x_{g}\in \on{exp}(\hat\f_{n})$ such that $(\on{Ad}g)(x_{1}+...+x_{n})
= g(x_{1}+...+x_{n})g^{-1}$. \hfill \qed\medskip \hfill \qed\medskip 

We then have: 

\begin{proposition}
$J(\alpha^{O}_{f}) = \langle \sum_{i=1}^{n} \on{log}\Gamma_{f}(\on{log}X_{i})
- \on{log}\Gamma_{f}(\on{log} \prod_{i=1}^{n}X_{i})\rangle$. 
\end{proposition}

{\em Proof.} Fix $\Phi\in M_{1}(\kk)$. We have $\mu_{f*\Phi}^{O} = 
\mu_{\Phi}^{O} \circ \alpha_{f}^{O}$, so $J(\mu_{f*\Phi}^{O})
= J(\mu_{\Phi}^{O}) + \mu_{\Phi}^{O} \circ J(\alpha)_{f}^{O}$. It follows that 
$ \mu_{\Phi}^{O} \circ J(\alpha)_{f}^{O} = 
\langle \sum_{i=1}^{n}\on{log}\Gamma_{f}(x_{i}) 
- \on{log}\Gamma_{f}(\sum_{i=1}^{n}x_{i})\rangle$. 
The result then follows from $\mu_{Phi}^{O}(X_{i})\sim e^{x_{i}}$, 
$\mu_{\Phi}^{O}(X_{1}...X_{n}) \sim e^{x_{1}+...+x_{n}}$.
\hfill \qed\medskip

\begin{remark}
In \cite{AT}, the Lie subalgebra ${\mathfrak{sder}}_{n}\subset \tder_{n}$
of special derivations (normalized special in the terms of Ihara) was introduced: 
${\mathfrak{sder}}_{n}=\{u\in \tder_{n}|u(x_{1}+...+x_{n})=0\}$. 
Let $\tilde{\mathfrak{sder}}_{n}$ be the intermediate Lie algebra
$\tilde{\mathfrak{sder}}_{n} = \{u\in\tder_{n}|\exists u_{0}
\in \f_{n-1} | u(x_{1}+...+x_{n})=[u_{0},x_{1}+...+x_{n}]\}$
(special derivations in Ihara's terms). So ${\mathfrak{sder}}_{n}
\subset \tilde{\mathfrak{sder}}_{n} \subset \tder_{n}$. Then 
Lemma \ref{lemma:special} says that we have a diagram 
$$
\begin{matrix}
\t_{n} &\to &{\mathfrak{sder}}_{n} & & \\
\downarrow & & \downarrow & & \\
\t_{n+1} &\to & \tilde{\mathfrak{sder}}_{n} & \hookrightarrow
& \tder_{n}\end{matrix}$$
\end{remark}

\begin{remark} Set $\on{SolKV}_n(\kk):= \{\mu_n\in \on{TAut}_n | 
\mu_n(e^{x_1}...e^{x_n})=e^{x_1+...+x_n}$ and $\exists r\in u^2\kk[[u]]| 
J(\mu_n) = \langle r(\sum_i x_i) - \sum_i r(x_i)\rangle\}$. 
This is a torsor under the action of the groups 
$\on{KV}_n(\kk):= \{\alpha_n\in \on{TAut}_n | \alpha_n(e^{x_1}...e^{x_n})
=e^{x_1}...e^{x_n}$ and $\exists\sigma\in u^2\kk[[u]] | J(\alpha) = 
\langle \sigma(\on{log}e^{x_1}...e^{x_n}) 
- \sum_i\sigma(x_i)\rangle\}$ and $\on{KRV}_n(\kk)$
defined similarly (replacing $e^{x_1}...e^{x_n}$ by 
$e^{x_1+...+x_n}$). These are prounipotent groups; the Lie 
algebra of $\on{KRV}_n(\kk)$ is $\krv_n :=
\{u\in \tder_n | a(\sum_i x_i)=0$ and $\exists s\in u^2\kk[[u]] | 
j(a) = \langle s(\sum_i x_i) - \sum_i s(x_i)\rangle\}$. 
It contains as a Lie subalgebra $\krv_n^0:= \{a\in \krv_n | s=0\}$, 
which is denoted $\kv_n$ in \cite{AT}. One can prove that if 
$|O'|=n$ and $O = \bullet\otimes O'$, the map $M_1(\kk)\to \on{SolKV}_n(\kk)$, 
$\Phi\mapsto \mu_{\Phi,O}$ is a morphism of torsors. 
\end{remark}

\section{Analytic aspects} \label{sec:anal}

In this section, the base field $\kk$ is $\RR$ or $\CC$. 

\subsection{Analytic germs}

We set
$\RR_{+}\{\!\{x\}\!\} := \{f \in \RR_{+}[[x]] | f$ has posititive
radius of convergence$\}$ and $\RR_{+}\{\!\{x\}\!\}_{0}:= 
\{f\in \RR_{+}\{\!\{x\}\!\} | f(0)=0\}$. If $f,g\in \RR_{+}[[r]]$, 
we write $f\preceq g$ iff $g-f\in \RR_{+}[[r]]$. We define 
$f\preceq g$ similarly when $f,g\in\RR_{+}[[r_{1},...,r_{n}]]$. 

Let $V,E$ be finite dimensional vector spaces and let $|.|_{V},|.|_{E}$ be norms on $V,E$. 
The space of $E$-valued formal series on $V$ is $E[[V]] = \{f =\sum_{n\geq 0} f_{n}, 
f_{n}\in S^{n}(V^{*})\otimes E\}$. For $f_{n}\in S^{n}(V^{*})\otimes E$, viewed 
as an homogeneous polynomial $V\to E$, we set $|f_{n}|:= \on{sup}_{v\neq 0}
(|f_{n}(v)|_{E}/|v|_{V}^{n})$. An analytic germ on $V$ (at the neighborhood of $0$)
is a series $f\in E[[V]]$, such that $|f|(r):= \sum_{n\geq 0}|f_{n}|r^{n}
\in \RR_{+}\{\!\{r\}\!\}$. We denote by $E\{\!\{V\}\!\}
\subset E[[V]]$ the subspace of analytic germs, and by $E\{\!\{V\}\!\}_{0}
\subset E[[V]]_{0}$ the subspace defined by $f_{0}=0$. 

If $f\in E\{\!\{V\}\!\}$ and $\alpha = \sum_{n\geq 0}\alpha_{n}r^{n}
\in \RR_{+}[[r]]_{0}$, we say that 
$\alpha$ is a dominating series for $f$ is $|f_{n}|\leq \alpha_{n}$ for any $n$; 
we write this as $|f(v)|_{E}\preceq \alpha(|v|_{V})$. 

If $V_{1},...,V_{k}$ are finite dimensional vector spaces with norms $|.|_{V_{1}}$, ..., 
$|.|_{V_{k}}$, then we equip $V_{1}\oplus ...\oplus V_{k}$ with the norm 
$|(v_{1},...,v_{k})|:= \on{sup}_{k}|v_{i}|_{V_{i}}$. If $f$ is an analytic germ
$V_{1}\oplus...\oplus V_{k}\to E$, we decompose 
$f = \sum_{\nn\in\NN^{k}} f_{\nn}$, where $f_{\nn}:V_{1}\times...\times
V_{k}\to E$ is the $\nn$-multihomogeneous component of $f$. 
We then set 
$$
|f_{\nn}| := \on{sup}_{(x_{1},...,x_{k})\in \prod_{i}(V_{i} - \{0\})}
|f_{\nn}(x_{1},...,x_{k})|_{E}/|x_{1}|_{V_{1}}^{n_{1}}...|x_{k}|_{V_{k}}^{n_{k}}.
$$ 
Then $f$ is an analytic germ iff $|f|(r_{1},...,r_{n}):= 
\sum_{\nn}|f_{\nn}|r_{1}^{n_{1}}...r_{k}^{n_{k}}\in \RR_{+}[[r_{1},...,r_{k}]]$ 
converges in a polydisc. If $\alpha
= \sum_{n_{1},...,n_{k}\geq 0}\alpha_{n_{1},...,n_{k}}r_{1}^{n_{1}}...
r_{k}^{n_{k}} \in \RR_{+}[[r_{1},...,r_{k}]]$, we write $|f(v_{1},...,v_{k})|_{E}
\preceq \alpha(|v_{1}|_{V_{1}},...,|v_{k}|_{V_{k}})$ if for each $\nn$, 
$|f_{\nn}(v_{1},...,v_{k})|_{E}\leq
\alpha_{\nn}(|v_{1}|_{V_{1}},...,|v_{k}|_{V_k})$. 

Let now $\G$ be a finite dimensional Lie algebra; let $|.|$ be a norm on $\G$; 
let $M> 0$ be such that the identity $|[x,y]|
\leq M|x||y|$ holds. 

The specialization to $\G$ of the Campbell--Baker--Hausdorff series 
is a series $x*y = \on{cbh}(x,y)\in \G[[\G\times \G]]_{0}$. 

\begin{lemma}

1) The CBH series is an analytic germ $\G\times \G\to \G$; we have 
$|x*y|\preceq {1\over M}f(M(|x|+|y|))$, where 
$f(u) = \int_{0}^{u} - {{\on{ln}(2-e^{v})}
\over v}dv$. 

2) $\G\times\G\to \G$, $(x,y)\mapsto e^{\on{ad}x}(y)$ is an analytic germ, 
and $|e^{\on{ad}x}(y)| \preceq e^{M|x|}|y|$. 
\end{lemma}

{\em Proof.} 1) is proved as in \cite{Bk}, not making use of the final majorization 
${1\over{r+s}}\leq 1$. Using Dynkin's formula, one can
prove that 2) follows from $|(\on{ad}x)^{n}(y)|
\leq M^{n}|x|^{n}|y|$. \hfill \qed\medskip 

\subsection{$\on{TAut}_{n}^{an}(\G)$ and $\tder_{n}^{an}(\G)$}

We set $\on{TAut}_{n}(\G):= \{(a_{1},...,a_{n}) | a_{i}\in \G[[\G^{n}]]_{0}\}$
and define on this set a product by $(a_{1},...,a_{n})(b_{1},...,b_{n}):= 
(c_{1},...,c_{n})$, where 
$$c_{i}(x_{1},...,x_{n}):= b_{i}(e^{\on{ad}a_{1}(x_{1},...,x_{n})}(x_{1}),...,
e^{\on{ad}a_{n}(x_{1},...,x_{n})}(x_{n}))*a_{i}(x_{1},...,x_{n}). 
$$
This equips $\on{TAut}_{n}(\G)$ with a group structure. We set 
$\on{TAut}_{n}^{an}(\G):=\{(a_{1},...,a_{n})|a_{i}\in \G\{\!\{\G^{n}\}\!\}_{0}\}$. 

\begin{proposition}
$\on{TAut}_{n}^{an}(\G)$ is a subgroup of $\on{TAut}_{n}(\G)$. 
\end{proposition}

{\em Proof.} Let $(a_{1},..,a_{n})$ and $(b_{1},...,b_{n})$
belong to $\on{TAut}_{n}^{an}(\G)$. Let $\alpha(r),\beta(r)
\in \RR_{+}\{\!\{r\}\!\}_{0}$ be germs such that 
the identities $|a_{i}(x_{1},...,x_{n})| \preceq \alpha(\on{sup}_{i}|x_{i}|)$, 
$|b_{i}(x_{1},...,x_{n})| \preceq \beta(\on{sup}_{i}|x_{i}|)$ hold. 
Then 
\begin{equation*}
\begin{split}
|c_{i}(x_{1},...,x_{n})| & \preceq 
f_{M}(|a_{i}(x_{1},...,x_{n})| + |b_{i}(e^{\on{ad}a_{1}}(x_{1}),...,
e^{\on{ad}a_{n}}(x_{n}))|) \\
 & 
 \preceq f_{M}(\alpha(\on{sup}_{i}|x_{i}|) 
 + \beta(e^{M\alpha(\on{sup}_{i}|x_{i}|)}\on{sup}_{i}|x_{i}|))
 = \gamma(\on{sup}_{i}|x_{i}|), 
\end{split}
\end{equation*}
where $f_{M}(u) = {1\over M}f(Mu)$ and 
$\gamma(r) = f_{M}(\alpha(r) + e^{M\alpha(r)}\beta(r))$ has nonzero radius of 
convergence. Here we use the compatibility of norms with composition: namely, 
if $f\in E[[V_{1}\times..\times V_{n}]]_{0}$ and 
$g_{i}\in V_{i}[[W]]_{0}$, with $|f(v_{1},...,v_{n})|\preceq \alpha(|v_{1}|,...,|v_{n}|)$
and $|g_{i}(w)|\preceq \beta_{i}(|w|)$, then $h:= f\in (g_{1},...,g_{n})\in E[[W]]_{0}$
and $|h(w)|\preceq \alpha \circ (\beta_{1},...,\beta_{n})(|w|)$. We also use the 
non-decreasing properties of elements of 
$\RR_{+}[[r_{1},...,r_{n}]]_{0}$ (i.e., if $F\in \RR_{+}[[u_{1},...,u_{k}]]_{0}$
and $u_{i},u'_{i}\in \RR_{+}[[r_{1},...,r_{l}]]_{0}$ with 
$u_{i}\preceq u'_{i}$, then $F(u_{1},...)\preceq F(u'_{1},...)$. 
So $(a_{1},...,a_{n})(b_{1},...,b_{n})\in \on{TAut}_{n}^{an}(\G)$. 

If now $(a_{1},...,a_{n})\in \on{TAut}_{n}^{an}(\G)$, then its inverse
$(b_{1},...,b_{n})$ in $\on{TAut}_{n}(\G)$ is uniquely determined by the 
identities 
$$
b_{i}(x_{1},...,x_{n}) = -a_{i}(e^{\on{ad}b_{1}(x_{1},...,x_{n})}(x_{1}),....,
e^{\on{ad}b_{n}(x_{1},...,x_{n})}(x_{n})). 
$$
Let us show that each $b_{i}(x_{1},...,x_{n})$ is an analytic germ. 
For this, we define inductively the sequence $b^{(k)} = (b_{1}^{(k)},...,b_{n}^{(k)})$
by $b^{(0)} = (0,..,0)$, and 
$$
b^{(k+1)}_{i}(x_{1},...,x_{n}) = -a_{i}(e^{\on{ad}b^{(k)}_{1}(x_{1},...,x_{n})}(x_{1}),....,
e^{\on{ad}b^{(k)}_{n}(x_{1},...,x_{n})}(x_{n})). 
$$
One checks that $b^{(k)} = b^{(k-1)} + O(x^{k})$, so the sequence $(b^{(k)})_{k\geq 0}$
converges in the formal series topology; the limit $b$ is then the inverse of 
$a = (a_{1},...,a_{n})$. 

Let us now set $\beta_{k}:= \on{sup}_{i}|b_{i}^{(k)}|$ (if $u_{i}(r) = 
\sum_{k\geq 0} u_{i,k}r^{k}\in\RR_{+}[[r]]$ is a finite family, we
set $(\on{sup}_{i}u_{i})(r):= \sum_{k\geq 0}(\on{sup}_{i}u_{i,k})r^{k}$). 
We then have 
$$
|b_{i}^{(k+1)}(x_{1},...,x_{n})| \preceq \alpha(\on{sup}_{i}
|e^{\on{ad}b_{i}^{(k)}(x_{1},...,x_{n})}(x_{i})|) \preceq 
\alpha(e^{M\beta_{k}(\on{sup}_{i}|x_i|)}\on{sup}_{i}|x_{i}|),
$$ 
so $\beta_{k+1}(r)\preceq \alpha(e^{\beta_{k}(r)}r)$.

We now define a sequence $(\gamma_{k})_{k\geq 0}$ of elements of 
$\RR_{+}[[r]]_{0}$ by $\gamma_{0}=0$, 
$$
\gamma_{k+1}(r) = \alpha(e^{M\gamma_{k}(r)}r). 
$$
As the exponential function, mutiplication by $r$ and $\alpha$
are non-decreasing, we have $\beta_{k}\preceq \gamma_{k}$. 
On the other hand, we have $\gamma_{k}(r)=\gamma_{k-1}(r)+O(r^{k})$, 
so the sequence $(\gamma_{k})_{k}$ converges in $\RR_{+}[[r]]_{0}$
(one also checks that this sequence is non-decreasing). 
Its limit $\gamma$ then satisfies 
\begin{equation} \label{eq:gamma}
\gamma(r) = \alpha(e^{M\gamma(r)}r). 
\end{equation}
It is easy to show that (\ref{eq:gamma}) determines $\gamma(r)\in \RR[[r]]_{0}$ 
uniquely. 
On the other hand, the function $(\gamma,r)\mapsto \gamma - \alpha(e^{M\gamma}r)
=:F(\gamma,r)$ is analytic at the neighborhood of $(0,0)$, with differential 
at this point $\partial_{\gamma}F(0,0)d\gamma + \partial_{r}F(0,0)dr=
d\gamma - M\alpha'(0)dr$. We may then apply the implicit function theorem
and use the fact that the $d\gamma$-component of $dF(0,0)$ is nonzero to 
derive the existence of an analytic function $\gamma_{an}(r)$ 
satisfying (\ref{eq:gamma}). By the uniqueness of solutions of (\ref{eq:gamma}), 
we get that the expansion of $\gamma_{an}$ is $\gamma$, so 
$\gamma\in \RR_{+}\{\!\{r\}\!\}_{0}$. 

Now $|b_{i}^{(k)}(x_{1},...,x_{n})|
\preceq \beta_{k}(\on{sup}_{i}|x_{i}|) \preceq \gamma_{k}(\on{sup}_{i}|x_{i}|)
\preceq \gamma(\on{sup}_{i}|x_{i}|)$, so by taking the limit $k\to \infty$, 
$|b_{i}(x_{1},...,x_{k})| \preceq \gamma(\on{sup}_{i}|x_{i}|)$, which implies that 
$b_{i}\in \G\{\!\{\G^{n}\}\!\}_{0}$, as wanted. 
\hfill \qed\medskip 

According to \cite{AT}, we have a bijection
$$
\kappa : \on{TAut}_{n}\to 
\tder_{n}, \quad 
g\mapsto  \ell - g\ell g^{-1}, 
$$
where $\ell$ is the derivation given by $x_{i}\mapsto x_{i}$. 

Set $\tder_{n}(\G):= \{(u_{1},...,u_{n})|u_{i}(x_{1},...,x_{n})
\in \G[[\G^{n}]]_{0}\}$, and $\tder_{n}^{an}(\G):= \{(u_{1},...,u_{n}) |
u_{i}\in \G\{\!\{\G^{n}\}\!\}_{0}\}\subset \tder_{n}(\G)$.
We have maps $\on{TAut}_{n}\to \on{TAut}_{n}(\G)$, 
$\tder_{n}\to \tder_{n}(\G)$ induced by the specialization of formal series. 

\begin{lemma}
1) There exists a map $\kappa_{\G} : \on{TAut}_{n}(\G)\to \tder_{n}(\G)$, 
such that the diagram 
$$\begin{matrix}
\on{TAut}_{n} &\stackrel{\kappa}{\to} &\tder_{n} \\
\downarrow & & \downarrow\\
 \on{TAut}_{n}(\G) &\stackrel{\kappa_{\G}}{\to} & \tder_{n}(\G)
  \end{matrix}$$commutes. 
  
2) This map restricts to a map
$\kappa^{an}_{\G} : \on{TAut}^{an}_{n}(\G)\to \tder^{an}_{n}(\G)$.  
\end{lemma}

{\em Proof.} 1) If $a_{i},b_{i}\in \hat\f_{n}$ are such that 
$g=\lbr e^{b_{1}},...,e^{b_{n}}\rbr$, $g^{-1}=\lbr e^{a_{1}},...,e^{a_{n}}\rbr$, 
then $\kappa(g) = u = \lbr u_{1},...,u_{n}\rbr$, with 
$$
u_{i}(x_{1},...,x_{n}) = 
({{1-e^{\on{ad}a_{i}}}\over{\on{ad}a_{i}}}(\dot a_{i}))
(e^{\on{ad}b_{1}(x_{1},...,x_{n})}(x_{1}),...,
e^{\on{ad}b_{n}(x_{1},...,x_{n})}(x_{n}))
$$
and $\dot a_{i} = \ell(a_{i}) = {d\over {dt}}_{|t=1}a_{i}(tx_{1},...,tx_{n})$. 
So we define $\kappa_{\G}$ by the same formula, where $\dot a_{i}$
is now defined as $ {d\over {dt}}_{|t=1}a_{i}(tx_{1},...,tx_{n})$
(or $\sum_{k\geq 0}ka_{i}^{k}$, where $a_{i}^{k}$ is the degree $n$
part of $a_{i}$). 

2) If the functions $a_{i},b_{i}$ are analytic germs, then so is $\dot a_{i}$
and therefore also each $u_{i}$. 
\hfill \qed\medskip 

Recall also from \cite{AT} that if $\mu\in \on{TAut}_{2}$,
$\mu(x*y)=x+y$ and $J(\mu) = \langle r(x)+r(y)-r(x+y) \rangle$
(i.e., $\mu\in\on{SolKV}(\kk)$), 
then $u:= -\kappa(\mu^{-1}) = \lbr A(x,y),B(x,y)\rbr$ satisfies: 

(KV1) $x+y-y*x = (1-e^{-\on{ad}x})(A(x,y)) + (e^{\on{ad}y}-1)(B(x,y))$, 

(KV3) $j(u) = \langle \phi(x)+\phi(y)-\phi(x*y)\rangle$, where $\phi(t)=tr'(t)$. 

Let $\Phi_{\on{KZ}}$ be the KZ associator, $\tilde\Phi_{\on{KZ}}(a,b):=
\Phi_{\on{KZ}}(a/(2\pi\i),b/(2\pi\i))\in M_{1}(\CC)$ and 
$\mu_{\on{KZ}}:= \mu_{\tilde\Phi_{\on{KZ}}}$. Let $u_{\on{KZ}}:= 
\kappa(\mu_{\on{KZ}}^{-1})$. Then $J(\mu_{\on{KZ}}) = \langle 
r_{\on{KZ}}(x)+r_{\on{KZ}}(y)-r_{\on{KZ}}(x*y)\rangle$, 
where $r_{\on{KZ}}(u)=-\sum_{n\geq 2}(2\pi\i)^{-n}\zeta(n)u^{n}/n$, 
therefore 
$$
j(u_{\on{KZ}}) = \langle \phi_{\on{KZ}}(x)+\phi_{\on{KZ}}(y)
-\phi_{\on{KZ}}(x*y)\rangle, 
$$
where $\phi_{\on{KZ}}(u) = -\sum_{n\geq 2}(2\pi\i)^{-n}\zeta(n)u^{n}$. 
Now the real part of this function (obtained by taking the
real part of the coefficients of $u^{n}$) is 
$$
\phi_{\on{KZ}}^{\RR}(u) = {1\over 2}
({u\over{e^{u}-1}}-1 + {u\over 2}).
$$ 

Let us now set $u_{\RR} := \lbr A_{\RR}(x,y),B_{\RR}(x,y)\rbr$, where the real part is 
taken with respect to the natural real structure on $\f_{2}^{\CC}$. 
Then by the linearity of (KV1), (KV3), we have: 
\begin{equation*}
\on{(KV1)} \quad x+y-y*x 
=  (1-e^{-\on{ad}x})(A_{\RR}(x,y)) 
+ (e^{\on{ad}y}-1)(B_{\RR}(x,y)) 
\end{equation*}
$$
\on{(KV3)}\quad 
j(u_{\RR}) = 
{1\over 2}\langle {{x}\over {e^{x}-1}} +  {{y}\over {e^{y}-1}}
-  {{x*y}\over {e^{x*y}-1}}-1\rangle. 
$$

\subsection{Analytic aspects to the KV conjecture (proof of Theorem 
\ref{thm:an})}

Recall that $\on{log}\tilde\Phi_{{\on{KZ}}}\in \hat\f_{2}$. We denote the
specialization of this series to the Lie algebra $\G$ as 
$(\on{log}\tilde\Phi_{{\on{KZ}}})^{\G}\in \G[[\G^{2}]]_{0}$. 

\begin{proposition} \label{prop:an}
$(\on{log}\tilde\Phi_{{\on{KZ}}})^{\G}$ is an analytic germ, i.e., 
$(\on{log}\tilde\Phi_{{\on{KZ}}})^{\G}\in \G\{\!\{\G^{2}\}\!\}_{0}$. 
\end{proposition}

{\em Proof.} Recall that $A_{2} = U(\f_{2})$ is the free associative algebra in 
$a,b$. For $x\in A_{2}$, set 
$$
|x|_{A_{2}}:= \on{sup}_{N\geq 1} \on{sup}_{m_{1},m_{2}
\in M_{N}(\CC)} || x(m_{1},m_{2})||. 
$$
Here $||.||$ is an algebra norm on $M_{N}(\CC)$. Then 
$|x|_{A_{2}}$ is $\leq \sum_{I\in \sqcup_{n\geq 0}
\{0,1\}^{n}}|x_{I}|$, where $x = \sum_{I}x_{I}e_{I}$, 
and for $I=(i_{1},...,i_{n})$, $e_{I}=e_{i_{1}}...e_{i_{n}}$, 
$e_{0}=a$, $e_{1}=b$. It follows from the Amitsur--Levitsky 
theorem (\cite{AL}) that $(|x|_{A_{2}}=0)\Rightarrow (x=0)$; indeed, 
by this theorem, $x(m_{1},m_{2})=0$ for $m_{1},m_{2}
\in M_{N}(\CC)$ implies: (a) that $x$ is in the 2-sided ideal generated 
by $ab-ba$ if $N=1$; (b) that $x=0$ if $N>1$. It follows that 
$|.|_{A_{2}}$ is an algebra norm\footnote{We will not use 
$(|x|_{A_{2}}=0)\Rightarrow (x=0)$, so our proof is 
independent of the Amitsur--Levitsky theorem.} on $A_{2}$, in particular 
$|xy|_{A_{2}}\leq |x|_{A_{2}}|y|_{A_{2}}$. 

We then define a vector space norm $|.|_{\f_{2}}$ on $\f_{2}$
by $|x|_{\f_{2}}:= |x|_{A_{2}}$; we have $|[x,y]_{\f_{2}}\leq
2|x|_{\f_{2}}|y|_{\f_{2}}$. 

For $\nn = (n_{1},...,n_{d})\in\NN^{d}$, and $f$ a function on 
$(\f_{2})^{d}$ (resp., $\RR^{d}$), we denote by $f(\xi_{1},...,\xi_{d})_{\nn}$ 
(resp., $f(t_{1},...,t_{d})_{\nn}$) the 
$\nn$-multihomogeneous part of $f$, which we view as a multihomogeneous 
polynomial on $(\f_{2})^{d}$ (resp., $\RR^{d}$). 

\begin{lemma} \label{lemma:prev}
For any $\nn$, we have the identity 
$$
|\on{log}(e^{\xi_{1}}...e^{\xi_{d}})_{\nn}|_{\f_{2}}\leq 
((\on{log}(2-e^{t_{1}+...+t_{d}})^{-1})_{\nn})_{
t_{1}=|\xi_{1}|_{\f_{2}},...,t_{d}=|\xi_{d}|_{\f_{2}}}. 
$$
\end{lemma}

{\em Proof of Lemma.} 
We have for any $\nn$, 
$|\xi_{1}^{n_{1}}...\xi_{d}^{n_{d}}|_{A_{2}}\leq 
|\xi_{1}|_{\f_{2}}^{n_{1}}...|\xi_{d}|_{\f_{2}}^{n_{d}}$ so 
$$|(e^{\xi_{1}}...e^{\xi_{d}}-1)_{\nn}|_{A_{2}} \leq 
((e^{t_{1}+...+t_{d}}-1)_{\nn})_{
t_{1}=|\xi_{1}|_{\f_{2}},...,t_{d}=|\xi_{d}|_{\f_{2}}}.$$
Then 
$\on{log}(e^{\xi_{1}}...e^{\xi_{d}})_{\nn} = 
\sum_{k\geq 1} {{(-1)^{k+1}}\over k}
\sum_{(\nn_{1},...,\nn_{k})|\nn_{1}+...+\nn_{k} = \nn}
(e^{\xi_{1}}...e^{\xi_{d}}-1)_{\nn_{1}}...
(e^{\xi_{1}}...e^{\xi_{d}}-1)_{\nn_{k}}$ so 
\begin{equation*}\begin{split}
& |\on{log}(e^{\xi_{1}}...e^{\xi_{d}})_{\nn}|_{A_{2}}  
\leq \big( \sum_{k\geq 1}{1\over k}
\sum_{\nn_{1}+...+\nn_{k}=\nn}
(e^{t_{1}+...+t_{d}}-1)_{\nn_{1}}...
(e^{t_{1}+...+t_{d}}-1)_{\nn_{d}} \big)_{
t_{1}=|\xi_{1}|_{\f_{2}},...,t_{d}=|\xi_{d}|_{\f_{2}}}\\
  & 
  = \big( \sum_{k\geq 1}{1\over k} ((e^{t_{1}+...+t_{d}}-1)^{k})_{\nn}
  \big)_{t_{1}=|\xi_{1}|_{\f_{2}},...,t_{d}=|\xi_{d}|_{\f_{2}}} 
  = ((\on{log}(2-e^{t_{1}+...+t_{d}})^{-1})_{\nn})_{
t_{1}=|\xi_{1}|_{\f_{2}},...,t_{d}=|\xi_{d}|_{\f_{2}}}. 
 \end{split}\end{equation*}
\hfill \qed\medskip 

Let $a(t)$ be an function $[0,1]\to \hat\f_{2}$ of the form 
$a(t) =\sum_{k\geq 1} a_{k}(t)$, where $a_{k}(t)\in \f_{2}[k]$ 
(here $k$ is the total degree in $a,b$) and $\int_{0}^{1}
|a_{k}(t)|_{\f_{2}}dt <\infty$. Let $u_{0},u_{1}$ be solutions of 
$u'(t)=a(t)u(t)$ with $u_{0}(0)=u_{1}(1)=1$, and $U:= u_{1}^{-1}u_{0}$. 


\begin{lemma} \label{lemma:ineq}
For $n\geq 1$, let $(\on{log}U)_{n}$ the degree $n$ (in $a,b$)
part of $\on{log}U$. Then 
$$
\sum_{n\geq 1} |(\on{log}U)_{n}|_{\f_{2}}r^{n} \preceq 
\on{log}(2-e^{\sum_{k\geq 1} r^{k}\int_{0}^{1}|a_{k}(t)|_{\f_{2}}dt})^{-1}.
$$  
\end{lemma}

{\em Proof of Lemma.} Let $\on{Lie}(n)$ be the multilinear part
of $\f_{n}$ in the generators $x_{1},...,x_{n}$. We denote by 
$w_{n}(x_{1},...,x_{n})\in \on{Lie}(n)$ the multilinear part
of $\on{log}(e^{x_{1}}...e^{x_{n}})$. 

Let now $\alpha_{n}$ be the coefficient of $t_{1}...t_{n}$
in the expansion of $\on{log}(2-e^{t_{1}+...+t_{n}})^{-1}$
(this is also the $n$th derivative at $t=0$ of $\on{log}(2-e^{t})^{-1}$). 
Specializing Lemma \ref{lemma:prev} for $\nn=(1,..,1)$, we get the identity
$$
|w_{n}(\xi_{1},...,\xi_{n})|_{\f_{2}} \leq 
\alpha_{n}|\xi_{1}|_{\f_{2}}...|\xi_{n}|_{\f_{2}}
$$
for $\xi_{1},...,\xi_{n}\in\f_{2}$. 

Now $\on{log}U$ expands as 
$$
\on{log}U = \sum_{n\geq 0} \int_{0<t_{1}<...<t_{n}<1}
w_{n}(a(t_{1}),...,a(t_{n}))dt_{1}...dt_{n} 
$$
(see e.g. \cite{EG}). It follows that 
$$
(\on{log}U)_{k} = \sum_{n\geq 0} \sum_{k_{1},...,k_{n}|\sum_{i}k_{i}=k}
\int_{0<t_{1}<...<t_{n}<1}
w_{n}(a_{k_{1}}(t_{1}),...,a_{k_{n}}(t_{n}))dt_{1}...dt_{n} 
$$
and therefore 
$$
|(\on{log}U)_{k}|_{\f_{2}} \leq 
\sum_{n\geq 0} \alpha_{n} \sum_{k_{1},...,k_{n}|\sum_{i}k_{i}=k}
\int_{0<t_{1}<...<t_{n}<1}
|a_{k_{1}}(t_{1})|_{\f_{2}}...|a_{k_{n}}(t_{n})|_{\f_{n}}dt_{1}...dt_{n} . 
$$
Now the generating series for the r.h.s. is 
$\on{log}(2-e^{\sum_{k\geq 1} r^{k}\int_{0}^{1}|a_{k}(t)|_{\f_{2}}dt})^{-1}$, 
proving the result. 
\hfill \qed\medskip 

According to \cite{Dr:Gal}, Section 2, if we set 
$$
a(t):= \sum_{k\geq 0,l\geq 1}{1\over{k!l!(2\pi\i)^{k+l+1}}}
{{(-\on{log}(1-t))^{k}(-\on{log}t)^{l}}\over{t-1}}
(\on{ad}b)^{k}(\on{ad}a)^{l}(b),$$
then $\tilde\Phi_{KZ}=U$. 
We have $|(\on{ad}b)^{k}(\on{ad}a)^{l}(b)|_{\f_{2}}\leq k+l+2\leq
2^{k+l+1}$, so 
$$
|a_{n}(t)|\leq \sum_{k\geq 0,l\geq 1,k+l+1=n}
{1\over{\pi^{k+l+1}k!l!}}
{{(-\on{log}(1-t))^{k}(-\on{log}t)^{l}}\over{1-t}}
$$
Then we have the inequality of formal series in $r$ 
\begin{equation*}\begin{split}
\sum_{n\geq 1} r^{n}\int_{0}^{1}|a_{n}(t)|_{\f_{2}}dt & \preceq 
\int_{0}^{1} \sum_{k\geq 0,l\geq 1}
{{r^{k+l+1}}\over{\pi^{k+l+1}k!l!}}
{{(-\on{log}(1-t))^{k}(-\on{log}t)^{l}}\over{1-t}}dt
\\
& = {r\over \pi}\int_{0}^{1} (1-t)^{-1-{r\over\pi}}(t^{-{r\over\pi}}-1)dt. 
\end{split}\end{equation*}
Now the identity $\int_{0}^{1}t^{a}(1-t)^{b}dt = {{\Gamma(a+1)\Gamma(b+1)}
\over{\Gamma(a+b+2)}}$, valid for $\Re(a), \Re(b)>-1$, implies that 
if $\Re(r)<0$, then 
$$
{r\over \pi}\int_{0}^{1} (1-t)^{-1-{r\over\pi}}(t^{-{r\over\pi}}-1)dt
=  {1\over 2}\big( 1 - {{\Gamma(1-2r)^{2}}\over{\Gamma(1-4r)}}\big). 
$$
This implies that the radius of convergence of 
${r\over \pi}\int_{0}^{1} (1-t)^{-1-{r\over\pi}}(t^{-{r\over\pi}}-1)dt$
is $1/4$, so this series belongs to $\RR_{+}\{\!\{r\}\!\}_{0}$. 
Plugging this in Lemma \ref{lemma:ineq}, we get 
$$
\sum_{n\geq 0}|(\on{log}\tilde\Phi_{\on{KZ}})_{n}|_{\f_{2}}
r^{n}\preceq \on{log}(2-
e^{{1\over 2}\big( 1 - {{\Gamma(1-2r)^{2}}\over{\Gamma(1-4r)}}\big)})^{-1}, 
$$ 
where the series in the r.h.s. lies in $\RR_{+}\{\!\{r\}\!\}_{0}$ (being a
composition of two series in $\RR_{+}\{\!\{r\}\!\}_{0}$). 

Let us now prove that $(\on{log}\tilde\Phi_{\on{KZ}})^{\G}\in 
\G\{\!\{\G^{2}\}\!\}_{0}$ is an analytic germ. By Ado's theorem, 
there exists a injective morphism $\rho:\G\to M_{N}(\kk)$,
where $\kk = \RR$ or $\CC$, hence an injective morphism $\tilde\rho:
\G\to M_{N}(\CC)$. Equip $\G$ with the norm $|x|_{\G}
:= ||\tilde\rho(x)||$. We recall that all the norms on $\G$ are equivalent, so 
it will suffice to prove analyticity w.r.t. $|.|_{\G}$. 

The degree $n$ part of the series $(\on{log}\tilde\Phi_{\on{KZ}})^{\G}$
is the specialization to $\G$ of $(\on{log}\tilde\Phi_{\on{KZ}})_{n}$. 
Now if $\psi\in \f_{2}[n]$ and $\psi^{\G}:\G\times\G\to\G$ is its specialization
to $\G$, we have $|\psi^{\G}(x,y)|_{\G} = ||\psi(\tilde\rho(x),\tilde\rho(y))||
\leq |\psi|_{\f_{2}} \on{sup}(||\tilde\rho(x)||,||\tilde\rho(y)||)^{n}
= |\psi|_{\f_{2}} \on{sup}(|x|_{\G},|y|_{\G})^{n}$, therefore 
$|\psi^{\G}| \leq |\psi|_{\f_{2}}$. We then have 
$$
\sum_{n\geq 0}|(\on{log}\tilde\Phi_{\on{KZ}})^{\G}_{n}|r^{n}
\preceq 
\sum_{n\geq 0}|(\on{log}\tilde\Phi_{\on{KZ}})_{n}|_{\f_{2}}r^{n}
\preceq 
\on{log}(2-
e^{{1\over 2}\big( 1 - {{\Gamma(1-2r)^{2}}\over{\Gamma(1-4r)}}\big)})^{-1};  
$$
together with the fact that the
series in the right has positive radius of convergence, this implies the analyticity of 
the series $(\on{log}\tilde\Phi_{\on{KZ}})^{\G}$. 
\hfill \qed\medskip 

Proposition \ref{prop:an}, together with the local analyticity of 
the CBH series, implies that the specialization of 
$\mu_{\tilde{\Phi}_{\on{KZ}}}$ belongs to $\on{TAut}^{an}_{2}(\G)$. 
It follows that $A(x,y)$, $B(x,y)$ are analytic germs, and so 

\medskip 

(KV2) 
$(A^{\RR},B^{\RR})$ is an analytic germ 
$\G^{2}\to\G^{2}$.

\medskip 

All this implies that 
$(A^{\RR},B^{\RR})$ is a solution of the `original' KV 
conjecture (as formulated in \cite{KV}) and proves 1) in 
Theorem \ref{thm:an}.

Let us now prove Theorem \ref{thm:an}, 2).  One checks easily that 
if $(A,B)$ is a solution of the `original' KV conjecture, then 
$(A_{s},B_{s}):= (A+s(\on{log}(e^{x}e^{y})-x), B+s(\on{log}(e^{x}e^{y})-y))$
is a family of solutions. In fact, if $\mu\in \on{SolKV}(\kk)$
and $\lbr A,B\rbr = -\kappa(\mu^{-1})$, then $\lbr A_{s},B_{s} \rbr
= -\kappa(\mu_{-s}^{-1})$, where $\mu_{s} := \on{Inn}(e^{s(x+y)})\circ \mu$; 
this corresponds to the action of `trivial', degree 1 element of $\mathfrak{krv}$ 
on $\on{SolKV}$ (see\cite{AT}).

Finally, let us prove Theorem \ref{thm:an}, 3). Let $\sigma$ be the 
antilinear automorphism of $\hat\f_{2}$ such that $\sigma(x)=-y$, 
$\sigma(y)=-x$. The series $\Phi_{\on{KZ}}(a,b)$ is real, therefore
$\overline{\tilde\Phi_{\on{KZ}}(a,b)}
= \tilde\Phi_{\on{KZ}}(-a,-b)$ (the bar denotes the complex conjugation). This implies 
that $\mu_{{\on{KZ}}}\circ\sigma = \on{Inn}(e^{-(x+y)/2})\circ
\sigma\circ\mu_{{\on{KZ}}}$. 
Using $\sigma\circ\ell\circ\sigma^{-1}=\ell$ and 
$\ell(x+y)=x+y$, we get 
$$
(\mu_{{\on{KZ}}}\circ\sigma\circ
\mu_{{\on{KZ}}}^{-1}) \circ \ell \circ
(\mu_{{\on{KZ}}}\circ\sigma\circ
\mu_{{\on{KZ}}}^{-1})^{-1} = \ell + 
\on{inn}({1\over 2}(x+y)), 
$$
where $\on{inn}(x+y)$ is the inner derivation $z\mapsto [x+y,z]$
of $\hat\f_{2}$.  Using now $\mu_{{\on{KZ}}}^{-1}(x+y)=
\on{log}(e^{x}e^{y})$, we get 
$$
(\sigma\circ
\mu_{{\on{KZ}}}^{-1}) \circ \ell \circ
(\sigma\circ
\mu_{{\on{KZ}}}^{-1})^{-1} = 
\mu_{{\on{KZ}}}^{-1}\circ \ell 
\circ \mu_{{\on{KZ}}}+ 
\on{inn}({1\over 2}\on{log}(e^{x}e^{y})).  
$$
Since $\sigma\circ\ell\circ\sigma^{-1}=\ell$,  
$\mu^{-1}\circ\ell\circ\mu-\ell = -\lbr A_{\on{KZ}},B_{\on{KZ}}\rbr$
and $\on{inn}({1\over 2}\on{log}(e^{x}e^{y})) = \lbr
{1\over 2}(\on{log}(e^{x}e^{y})-x),{1\over 2}(\on{log}(e^{x}e^{y})-y)\rbr$ 
$$
\sigma \circ \lbr A_{\on{KZ}},B_{\on{KZ}}\rbr\circ\sigma^{-1}
= \lbr A_{\on{KZ}},B_{\on{KZ}} \rbr -  \lbr
{1\over 2}(\on{log}(e^{x}e^{y})-x),{1\over 2}(\on{log}(e^{x}e^{y})-y)\rbr. 
$$
This implies 
$$
(B_{\on{KZ}}(-y,-x),A_{\on{KZ}}(-y,-x)) = (A_{\on{KZ}}(x,y),B_{\on{KZ}}(x,y))
- ({1\over 2}(\on{log}(e^{x}e^{y})-x),
{1\over 2}(\on{log}(e^{x}e^{y})-y)). 
$$
If now
$(A',B') := (A_{\on{KZ}},B_{\on{KZ}})
- {1\over 4}(\on{log}(e^{x}e^{y})-x,\on{log}(e^{x}e^{y})-y)$, this 
implies 
$$(B'(-y,-x),A'(-y,-x))=(A'(x,y),B'(x,y)),$$ which by taking real parts
implies $(B_{-1/4}(-y,-x),A_{-1/4}(-y,-x))=(A_{-1/4}(x,y),B_{-1/4}(x,y))$, 
proving Theorem \ref{thm:an}, 3). 

\begin{appendix}

\section{Results on centralizers} \label{app}

\subsection{The centralizer of $t_{ij}$ in $\t_{n}$}

\begin{proposition} \label{prop:comm}
Let $i< j\in [n]$. If $x\in \t_{n}$ is such that $[x,t_{ij}]=0$, 
then there exists $\lambda\in \kk$ and $y\in \t_{n-1}$ such that 
$x=\lambda t_{ij} + y^{ij,1,2,...,\check i,...,\check j,...,n}$. 
\end{proposition}

{\em Proof.} We may and will assume that $i=1,j=2$. We then 
prove the result by induction on $n$. It is obvious when $n=2$. 
Assume that it has been proved at step $n-1$ and let us prove it 
at step $n$. We have $\t_{n} = \t_{n-1}\oplus \f_{n-1}$, 
where $\t_{n-1}$ is the Lie subalgebra generated by the 
$t_{ij}$, $i\neq j\in \{1,...,n-1\}$ and $\f_{n-1}$ is 
freely generated by the $t_{1n},...,t_{n-1,n}$. Both 
$\t_{n-1}$ and $\f_{n-1}$ are Lie subalgebras of $\t_{n}$, 
stable under the inner derivation $[t_{12},-]$. Then if $x\in \t_{n}$
is such that $[t_{12},x]=0$, we decompose $x=x'+f$, with 
$x'\in \t_{n-1}$, $f\in \f_{n-1}$, $[t_{12},x']=[t_{12},f]=0$. 
By the induction hypothesis, we have $x'=\lambda t_{12} 
+ (y')^{12,3,...,n-1}$, where $y'\in \t_{n-2}$ and $\lambda \in\kk$. 

Let us set $x_{i} = t_{in}$ for $i=1,...,n-1$. The derivation $[t_{12},-]$
of $\f_{n-1}$ is given by $x_{1}\mapsto [x_{1},x_{2}]$, 
$x_{2}\mapsto [x_{2},x_{1}]$, $x_{i}\mapsto 0$ for $i>2$. 
In terms of generators $y_{1} = x_{1}$, $y_{2}=x_{1}+x_{2}$,
$y_{3}=x_{3}$...,
$y_{n-1}=x_{n-1}$, it is given by $y_{1}\mapsto [y_{1},y_{2}]$, 
$y_{i}\mapsto 0$ for $i>1$. 

\begin{lemma}
The kernel of the derivation $y_{1}\mapsto [y_{1},y_{2}]$, 
$y_{i}\mapsto 0$ for $i>1$ of $\f_{n-1}$ coincides with the Lie 
subalgebra $\f_{n-2}\subset \f_{n-1}$ generated by $y_{2},...,y_{n-1}$.
\end{lemma} 

{\em Proof of Lemma}. Let us
prove that the kernel of the induced derivation of $U(\f_{n-1})$
is $U(\f_{n-2})$. We have a linear isomorphism 
$U(\f_{n-1}) \simeq \oplus_{k\geq 1} U(\f_{n-2})^{\otimes k}$, 
whose inverse takes $u_{1}\otimes ...\otimes u_{k}$ to 
$u_{1}y_{1}u_{2}y_{1}...y_{1}u_{k}$. The derivation 
$[t_{12},-]$ of $U(\f_{n-1})$ is then transported to the 
direct sum of the endomorphisms of $U(\f_{n-2})^{\otimes k}$
\begin{equation} \label{end:k}
u\mapsto (y_{2}^{(2)}+...+y_{2}^{(k)})u-u(y_{2}^{(1)}+...+y_{2}^{(k-1)})
\end{equation}
(this is $0$ of $k=1$; $y_{2}^{(i)} = 1^{\otimes i-1}\otimes y_{2}\otimes
1^{\otimes k-i}$; we make use of the algebra structure of $U(\f_{n-2})^{\otimes k}$). 
Each of these endomorphisms has degree 1 for the 
filtration of $U(\f_{n-2})^{\otimes k}$ induced by the PBW 
filtration of $U(\f_{n-2})$ (the part of degree $\leq d$ of $U(\f_{n-2})$
for this filtration consists of combinations of products of $\leq d$ elements of 
$\f_{n-2}$) and the associated graded endomorphism of 
$S(\f_{n-2})^{\otimes k}$ is the multiplication by $y_{2}^{(k)}-y_{2}^{(1)}$, 
which is injective if $k\geq 1$, so (\ref{end:k}) is injective for $k\geq 1$; 
the kernel of the direct sum of maps (\ref{end:k}) therefore coincides with 
the degree $1$ part $U(\f_{n-2})$, which transports to 
$U(\f_{n-2}) \subset U(\f_{n-1})$. So the kernel of the derivation 
$[t_{12},-]$ of $U(\f_{n-1})$ is $U(\f_{n-2})$. 
The kernel of the derivation $[t_{12},-]$ of $\f_{n-1}$ is then 
$\f_{n-1}\cap U(\f_{n-2}) = \f_{n-2}$. 
\hfill \qed\medskip 

{\em End of proof of Proposition \ref{prop:comm}.}
It follows that $f$ expresses as $P(t_{1n}+t_{2n},t_{3n},...,t_{n-1,n})$. 
Then if we set $f':= P(t_{1,n-1},...,t_{n-2,n-1})$, we get 
$f=(f')^{12,3,...,n}$ so $x = x'+f = \lambda t_{12} + 
((y')^{1,2,...,n-1}+f')^{12,3,...,n}$, as wanted. 
\hfill \qed\medskip 

\subsection{The centralizer of $x_{ij}$ in $\on{PB}_{n}$}

\begin{proposition} \label{app:centr}
If $g\in \on{PB}_{n}(\kk)$ commutes with $x_{12}$, then 
there exists $\lambda\in \kk$ and $h\in \on{PB}_{n-1}(\kk)$ such that 
$g= x_{12}^{\lambda} h^{\widetilde{12},3,...,n}$. 
\end{proposition}

Since $x_{ij}$ is conjugated to $x_{12}$, it is easy to derive from this 
the centralizer of $x_{ij}$ in $\on{PB}_{n}(\kk)$.

{\em Proof.} Note that $x_{12}$ commutes with the image of 
$\on{PB}_{n-1}(\kk)\to \on{PB}_{n}(\kk)$, $h\mapsto 
h^{\widetilde{12},3,...,n}$, so that $U_{0}:= 
\{x_{12}^{\lambda} h^{\widetilde{12},3,...,n} | 
h\in \on{PB}_{n-1}(\kk), \lambda\in\kk\}$ is an algebraic
subgroup of $\on{PB}_{n}(\kk)$. 
Let $U\subset \on{PB}_{n}(\kk)$ be the centralizer of $x_{12}$; 
then $U_{0}\subset U$, and we need to prove that $U_{0}=U$. 

We have $U_{0} = \on{exp}(\u_{0})$, 
$U = \on{exp}(\u)$, where $\u_{0} = \kk \on{log}x_{12} \oplus 
\on{Im}(\pb_{n-1}\stackrel{\widetilde{12},3,...,n}{\to} \pb_{n})$ 
and $\u = \{x\in \pb_{n} | [\on{log}x_{12},x]=0\}$, where
$\pb_{n}:= \on{Lie}\on{PB}_{n}(\kk)$. Then the lower central series
defines a complete decreasing filtration of $\pb_{n}$, with 
$F^{1}\pb_{n} = \pb_{n}$ and $F^{i+1}\pb_{n} = [\pb_{n},F^{i}\pb_{n}]$. 
The associated graded Lie algebra is $\t_{n}$, i.e., $\t_{n} = \oplus_{i\geq 1}\t_{n}[i]
= \oplus_{i\geq 1} F^{i}\pb_{n}/F^{i+1}\pb_{n}$. 

Set $F^{i}\u:= \u\cap F^{i}\pb_{n}$, $F^{i}\u_{0}:= \u_{0}\cap F^{i}\pb_{n}$. 
We will prove that the images of $F^{i}\u_{0}$ and $F^{i}\u$ in $\t_{n}[i]$
coincide. Clearly, $\on{Im}(F^{i}\u_{0}\to \t_{n}[i]) \subset 
\on{Im}(F^{i}\u\to \t_{n}[i])$. 

Conversely, projecting the identity 
$[\on{log}x_{12},x]=0$ modulo $F^{i+1}\pb_{n}$, we get 
\begin{equation} \label{interm:1}
\on{Im}(F^{i}\u \to \t_{n}[i]) \subset 
\{x\in \t_{n}[i]|[t_{12},x]=0\}, 
\end{equation} 
and 
since $x\mapsto x^{\widetilde{12},...,n}$
takes $F^{i}\pb_{n-1}$ to $F^{i}\pb_{n}$, we have  
$(F^{i}\pb_{n-1})^{\widetilde{12},...,n}
\subset F^{i}\u_{0}$ if $i>1$ and $(F^{1}\pb_{n-1})^{\widetilde{12},...,n}
\oplus \kk \on{log}x_{12}\subset F^{1}\u_{0}$; projecting these inclusions, 
modulo $F^{i+1}\pb_{n}$, we get  
\begin{equation} \label{interm:2}
\on{Im}(F^{i}\u_{0} \to \t_{n}[i]) \supset 
\t_{n-1}[i]^{12,...,n} \on{\ if\ }i>1 \on{\ and\ }
\on{Im}(F^{1}\u_{0} \to \t_{n}[i]) \supset 
\t_{n-1}[1]^{12,...,n}\oplus \kk t_{12}. 
\end{equation}
Using (\ref{interm:1}), (\ref{interm:2}) and Proposition 
\ref{prop:comm}, we get $\on{Im}(F^{i}\u\to\t_{n}[i]) \subset  
\on{Im}(F^{i}\u_{0}\to\t_{n}[i])$. It follows that these spaces are equal, 
which implies (as both $\u_{0}$ and $\u$ are closed for the topology of $\pb_{n}$)
that $\u_{0}=\u$. So $U_{0} = U$. \hfill \qed\medskip 

\begin{remark} One can also prove Proposition \ref{app:centr} similarly
to Proposition \ref{prop:comm}, by induction on $n$ and using the fact that 
the automorphism $\on{Ad}x_{12}$ of the topologically free group
generated by the $x_{in}$ identifies with the automorphism 
$\on{exp}(\on{ad}t_{12})$ of the topologically free Lie algebra generated by the 
$t_{in}$ (using the identification $(x_{1n},x_{1n}x_{2n},x_{3n},...,x_{n-1,n})
\leftrightarrow (e^{t_{1n}},e^{t_{1n}+t_{2n}},e^{t_{3n}},...,e^{t_{n-1,n}})$). 
\end{remark}

\end{appendix}

\end{document}